\documentclass[10pt,twoside]{article}
\usepackage[latin1]{inputenc}
\usepackage{amsmath}
\input epsf
\def\limiteN{\renewcommand{\arraystretch}{0.5}
\begin{array}[t]{c}\stackrel{}{\longrightarrow} \\
{\scriptstyle N\rightarrow
\infty}\end{array}\renewcommand{\arraystretch}{1}}
\def\limiteloi{\renewcommand{\arraystretch}{0.5}
\begin{array}[t]{c}\stackrel{{\cal D}}{\longrightarrow} \\
{\scriptstyle N\rightarrow
\infty}\end{array}\renewcommand{\arraystretch}{1}}
\def\limiteproba{\renewcommand{\arraystretch}{0.5}
\begin{array}[t]{c}\stackrel{{\cal P}}{\longrightarrow} \\
{\scriptstyle N\rightarrow
\infty}\end{array}\renewcommand{\arraystretch}{1}}
\def\egaleloi{\renewcommand{\arraystretch}{0.5}
\begin{array}[t]{c}\stackrel{{\cal D}}{\sim} \\
{}\end{array}\renewcommand{\arraystretch}{1}}
\newcommand{\be}{\begin{equation}}
\newcommand{\ee}{\end{equation}}
\newcommand{\bd}{\begin{displaymath}}
\newcommand{\ed}{\end{displaymath}}

\newcommand{\ba}{\begin{eqnarray}}
\newcommand{\ea}{\end{eqnarray}}
\newcommand{\ban}{\begin{eqnarray*}}
\newcommand{\ean}{\end{eqnarray*}}

\newcommand{\LL} {I\!\!L}
\newcommand{\R} {I\!\!R}
\newcommand{\E} {E\,}
\newcommand{\N} {I\!\! N}

\newcommand{\cov} {\mbox{cov}\hspace{0.4mm}}

\renewcommand{\arraystretch}{.8}
\renewcommand{\Box}{\hfill\rule{0.25cm}{0.25cm}} %% carre noir + rejet fin=ligne

\newtheorem{Def}{Definition}[section]
\newtheorem{Prop}{Proposition}[section]
\newtheorem{lem}{Lemma}[section]
\newtheorem{Theo}{Theorem}[section]
\newtheorem{cor}{Corollary}[section]
\newtheorem{rem}{Remark}[section]
\newenvironment{dem}{\ \\ {\bf Proof. }}
{\Box\par\medskip\noindent}

\def\1{{\bf 1}}
\def\Pr{I\mskip-7muP}

\begin{document}
%%\date{ }
\title{\bf Identification of the multiscale fractional Brownian motion with biomechanical applications}
\maketitle \vspace{-0.5cm}
~\\
Jean-Marc~BARDET${}^*$ and Pierre~BERTRAND${}^{**}$\\
~\\
${}^*$  {\it SAMOS-MATISSE - UMR CNRS 8595, Universit\'e
Panth\'eon-Sorbonne (Paris I),
90 rue de Tolbiac, 75013 Paris Cedex, France, E-mail: bardet@univ-paris1.fr} \\ \\
${}^{**}$ {\it Laboratoire de Math\'ematiques - UMR CNRS 6620,
Universit\'e Blaise Pascal (Clermont-Ferrand II), 24 Avenue des
Landais, 63117 Aubi\`ere Cedex, France. E-mail:
Pierre.Bertrand@math.univ-bpclermont.fr} \pagestyle{myheadings}
\markboth{Statistic of multi-scale fractional Brownian motion}{
J.M. Bardet and P. Bertrand} ~\\ ~\\
Abstract : In certain applications, for instance biomechanics,
turbulence, finance,  or Internet traffic, it seems suitable to
model the data by  a generalization of a fractional Brownian
motion for which the Hurst parameter $H$ is depending on the
frequency as a piece-wise constant function. These processes are
called multiscale fractional Brownian motions. In this
contribution, we provide a statistical study of the multiscale
fractional Brownian motions. We develop a method based on wavelet
analysis. By using this method, we find initially the frequency
changes, then we estimate the different parameters and afterwards
we test the goodness-of-fit. Lastly, we give the numerical algorithm.
Biomechanical data are then studied with these new tools.~\\ ~\\
{\em Keywords:} Biomechanics; Detection of change; Goodness-of-fit
test; Fractional Brownian motion; Semi-parametric estimation;
Wavelet analysis.
\newpage
\section{Introduction}
Fractional Brownian Motion (F.B.M.) was introduced in 1940 by
Kolmogorov as a way to generate Gaussian "spirals" in a Hilbert
space. But the seminal paper of Mandelbrot and  Van Ness (1968)
emphasizes the relevance of F.B.M. to model natural phenomena:
hydrology, finance...  Formally, a fractional Brownian motion
$B_H=\left ( B_H(t),\, t \in \R_+ \right )$ could be defined as a
real centered Gaussian process with stationary increments such
that $B_H(0)=0$ and $\; \E\left|B_H(s) - B_H(t) \right|^2
\,=\,\sigma^2\, |t-s|^{2H},\;$ for all pair $(s,t) \in \R_+\times \R_+$
where $H \in]0, 1[$ and $\sigma>0$. This process is characterized
by two parameters :  the Hurst index $H$ and the scale parameter
$\sigma$. We lay the emphasis on the fact that  the  same
parameter $H$ is linked to different properties of the F.B.M. as
the smoothness of the sample paths, the long range dependence of
its increments and the self-similarity. \\
~\\
During the decades 1970's and 1980's, the statistical study of
F.B.M. was developed, to look at for instance the historical notes
in Samorodnitsky \& Taqqu (1994), \cite[chap.14]{TaSa:1994} and
the references therein. Modelling by a F.B.M. became more and more
widespread during the last decade (traffic Internet, turbulence,
image processing...). Nevertheless, in many applications the real
data does not fit exactly F.B.M. Thus, the F.B.M. must be regarded
only as an ideal mathematical model. Therefore, various
generalizations of F.B.M. have been proposed these last years to
fill the gap between the mathematical modelling and  real data. In
one hand, Gaussian processes where the Hurst parameter $H$ has
been replaced by a function depending on the time were studied,
see for instance Peltier and L\'evy~Vehel (1996), Benassi, Jaffard
and Roux (1997), Ayache and L\'evy~Vehel (1999). However, this
dependence of time implies the loss of the stationarity of the
increments. In other hand, non Gaussian processes, mainly $\alpha$
stable $(0<\alpha<2)$ infinite variance processes, were
considered, see for example the study of telecom processes in Pipiras and Taqqu (2002).\\
~\\
Here, we are concerned with  Gaussian processes having stationary
increments and a Hurst index changing with the frequencies. To our
knowledge, these kinds of processes were introduced implicitly in
biomechanics by Collins and de Luca (1993), in finance by Rogers
(1997) and Cheridito (2003) and explicitly by Benassi and Deguy
(1999) for image analysis or image synthesis. In any case, the
probabilistic properties of these processes have not been
thoroughly  established   and no rigorous statistical studies have
been done. Both Collins and de Luca (1993) and Benassi and Deguy
(1999) propose a model with two different Hurst indices
corresponding respectively to the   high and the low frequencies
separated by one change point at the frequency $\omega_c$. They use
the log variogram to estimate these two Hurst indices. Indeed, in
this case, the log variogram considered as a function of the
logarithm of the scale presents two asymptotic directions with
slopes being twice the Hurst index at low (respectively high)
frequencies. The change point $\omega_c$ is then estimated as the
abscise of the intersection of the these two straight lines.
Numerically, this method is not robust. Moreover it could not be
adapted in the case of more than one change point. Let us stress
that it is not a question of a theoretical refinement, but one
that corresponds precisely to the true situations. Indeed, in
applications, we consider only finite frequency bands, therefore
we should use a statistical method based on the information
included in finite frequency bands. Wavelet analysis seems the
tool {\em had hoc}, when the Fourier transform of the associated
wavelet is compactly supported.\\
~\\
For these reasons, we put forward in Bardet and Bertrand (2003) a
model of generalized F.B.M. including the cases with more than one
frequency change point. We  called it $(M_K)$ multiscale
fractional Brownian motion where $K$ denote the number of
frequency change points. More precisely, a $(M_K)$ multiscale
fractional Brownian motion is a Gaussian process with stationary
increments where the Hurst parameter $H$ is replaced by a
piecewise constant function of the frequency $\xi \mapsto H(\xi)$
in the harmonizable representation, see Formula
(\ref{def:msm.b.f.}) below. The main probabilistic properties of
this model were studied in Bardet and Bertrand (2003). In this
work, we treat the statistical study of the multiscale F.B.M. and
we focus on its application to biomechanics. \\
~\\
The remainder of the paper is organized as follows: in Section 2,
we describe the biomechanical data and the corresponding
statistical problem. In section 3, we recall the initial
definition of the partial Brownian motion and its principal
probabilistic properties.  Then, we show that the variogram method
is not suitable for the estimation of the various parameters of a
$(M_K)$-F.B.M. We then develop a statistical estimation framework,
based on wavelet analysis. We investigate the discretization of
the wavelet coefficient and we state   a functional Central Limit
Theorem for the empirical wavelet coefficients. In Section 4, we
first estimate the different frequency change  points and Hurst
parameters.  Then, we propose a goodness of fit test and derive an
estimator of the number of frequency changes. The numerical
algorithm is detailed at the end of this section.  Finally, in
Section 5, the biomechanical data are studied with the tools
developed in Section 4. The proof of the results of Sections 3 and
4 are given in appendix.
\section{The Biomechanical Problem}
One of the motivations of this work is to model biomechanical data
corresponding to the regulation of the upright position of the
human being. By using a force platform, the position of the center
of pressure (C.O.P.) during quiet postural stance is determined.
This position is usually measured at a frequency of 100 Hz for the
one minute period, which yields a data set of 6000 observations.
The experimental conditions are formed to the standards of the
Association Fran\c{c}aise de Posturologie (AFP), for instance the feet
position (angle and clearance), the open or closed eyes.
\[
\epsfxsize 7cm \epsfysize 7cm \epsfbox{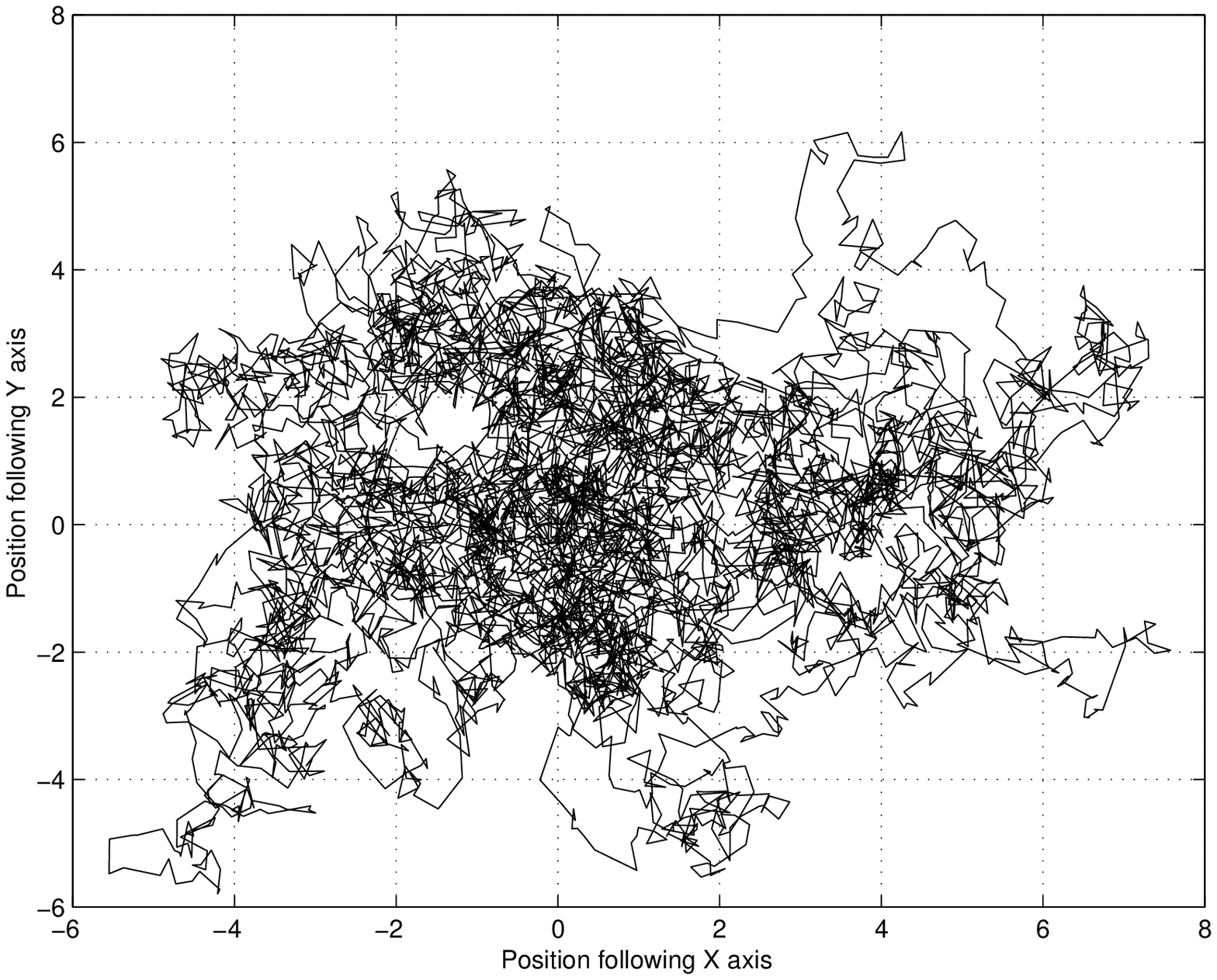}
\]
{\bf Figure 1 :} An example \footnote{these
experimental data were realized by A. Mouzat and are used in
\cite{B2DMV:IEEE}.} of the trajectory of the C.O.P. during 60s at $100$Hz (in mm)\\
~\\
The X axis of the platform corresponds to the fore-aft direction
and the Y axis corresponds to the medio-lateral direction. During
the 1970's, these data were analyzed as a set of points, {\it
i.e.} without taking into account their temporal order. During the
following decade some studies considered them as a process, and
Collins and de Luca (1993) introduced the use of F.B.M. to model
these data. In fact, they used a generalization of F.B.M. More
precisely, let the position  $X_i$ of the C.O.P. be observed at
times $t_i = i \Delta$ for $i=1,\dots,N$ ($\Delta=0.01~s$). The
study of Collins and de Luca is based on the empirical variogram
\be \label{def:V} V_N(\delta ) = \frac{1}{(N- \delta )}
\sum_{i=1}^{N-\delta} \left(X_{(i+\delta)\Delta } -X_{i\Delta}
\right)^2 \ee where $\delta   \in  \N^*$. For a F.B.M., we have
$\E V_N(\delta ) = \sigma ^2\,  \Delta ^{2H}\times \delta ^{2H}$ and
after plotting the log-log graph of the variogram as a function of
the time lag , i.e. $( \log \delta , \log V_N(\delta ))$, a linear
regression provides the slope $2H$. Typically, one gets the
following type of figure (see Figure 2). It is considered by
Collins and de Luca to be a "F.B.M." with two regimes~: with slope
$2H_0$ {\em (short term)} and with slope $2H_1$ {\em (long term)}
separated by a critical time lag $\delta_c$ and these parameters
are estimated graphically~:
\[
\epsfxsize 7cm \epsfysize 7cm \epsfbox{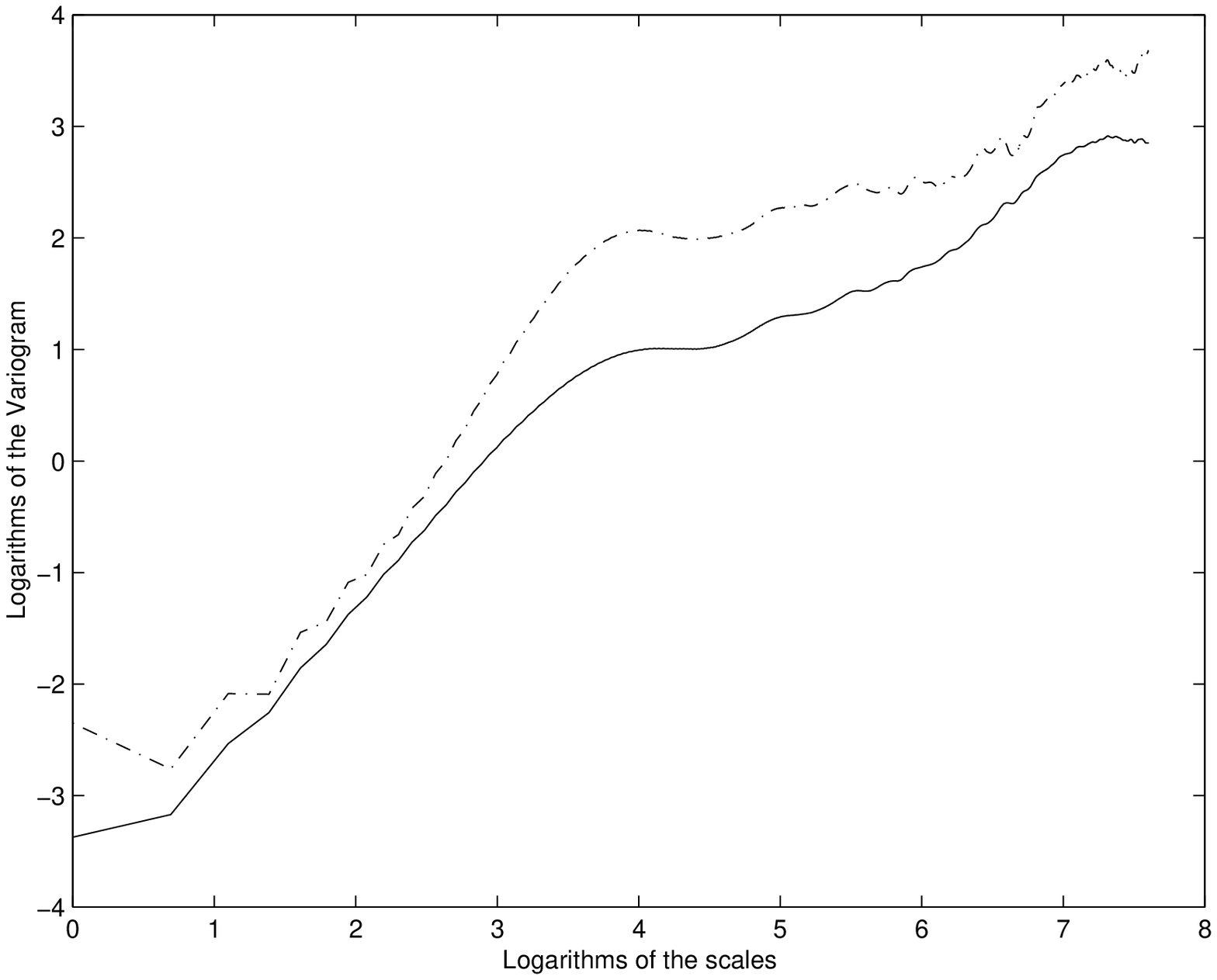}
\]
{\bf Figure 2 :} An example of the log-log graph of the variogram
for the previous trajectories X (-.) and Y (-).\\
~\\
They found $H_0 > 0.5$, $H_1 <0.5$ and a critical time lag
$\delta_c\simeq 1~s$. These results were interpreted as
corresponding to two different kinds of regulation  of the human
stance : in the {\em long term} $H_1 <0.5$ and the process is
anti-persistent, in the {\em short term} $H_0 > 0.5$ and the
process is persistent. This method was employed several times in
biomechanics under the various experimental conditions  (opened
eyes versus closed eyes, different feet angles,...). But, a lack
of mathematical models and of statistical studies has made
impossible to obtain confidence intervals on the two slopes
$2H_0$, $2H_1$ and the critical time lag $\delta_c$.
\section{The multiscale fractional Brownian motion and its statistical
study based on wavelet analysis}
\subsection{Description of the model}
A fractional Brownian motion  $B_H=\left\{ B_H(t),\, t \in \R
\right \}$ of parameters $(H,\,\sigma)$ is a real centered
Gaussian process with stationary increments and $\; \E\left|B_H(s)
- B_H(t) \right|^2 \,=\,\sigma^2\, |t-s|^{2H},\,$  for all $(s,t)
\in \R^2$ where $H \in]0, 1[$ and $\sigma >0$. The fractional
Brownian motion (F.B.M.) has been proposed by Kolmogorov (1940)
who defined it by the harmonizable representation~: \be
\label{def:BF}
B_H(t) =%%\frac{1}{C_H}
\int_{\R} \frac{\left(e^{it\xi} - 1\right)}{|\xi|^{H+1/2}} \,
\overline{\widehat{W}}(d\xi),~~~~\mbox{for all}~~t \in \R, \ee
where $W(dx)$ is a Brownian measure and $\widehat{W}(d\xi)$ its
Fourier transform (namely for any function $f\in L^{2}(\R)$ one
has almost surely, $\int_{\R} f(x) W(dx) \,=\, \int_{\R}
\widehat{f}(\xi) \,\widehat{W}(d\xi) $, with the convention that
$\widehat{f}(\xi)=\int_{\R} e^{-i\xi \, x} \,f(x) \,dx$ when $f\in
L^{1}(\R)\bigcap L^{2}(\R)$). We refer to Samorodnitsky and Taqqu
(1994) for the question of the equivalence of the different
representations of the F.B.M. From the harmonizable
representation, a natural generalization is the multiscale
fractional Brownian motion with a Hurst index depending on the
frequency. More precisely, we define~:
\begin{Def}
For $K\in \N$, a $(M_K)$-multiscale fractional Brownian motion
$X=\{X(t),t\in \R\}$ (simplify by $(M_K)$-F.B.M.) is a process
such as
\be \label{def:msm.b.f.} X(t) = 2 \sum _{j=0}^K
\int_{\omega_j}^{\omega_{j+1}} \sigma_j \frac {(e^{it\xi} -1\
)}{|\xi|^{H_j+1/2}}\, \overline{\widehat{W}}(d\xi)~~~~\mbox{for
all}~~t\in \R \ee with $\omega_0 =
0<\omega_1<\dots<\omega_K<\omega_{K+1} =\infty$ by convention,
$\sigma_i>$ and $H_i \in ]0,1[$ for $i \in \{0,1,\cdots,K\}$.
\end{Def}
The $(M_K)$-F.B.M. was notably introduced in order to relax the
self-similarity property of F.B.M. Indeed, the self-similarity is
a form of invariance with respect to changes of time scale
\cite{MvN:68} and it links the behavior to the high frequencies
with the behavior to the low frequencies. In Bardet and Bertrand
(2003), the main properties of these processes are provided : $X$
is a Gaussian centered process with stationary increments, its
trajectories are a.s. of H\"older regularity $\alpha$, for every
$0\le\alpha<H_K$ and its increments form a long-memory process
(except if the different parameters satisfy a particular
relationship, {\em i.e.}, if its spectral density is a continuous
function with $0<H_i<1/2$ for $i=0,1,\cdots,K$).
\subsection{The question of the choice of the estimator}
In the remainder of this paper, we suggest a statistical study of
such a model based on wavelet analysis. In this subsection, we
explain
the reason of this choice.  \\
~\\
To begin with, we will describe the statistical framework
precisely. Let $X=\{X(t),t\in \R_+\}$ be a $(M_K)$-F.B.M. defined
by (\ref{def:msm.b.f.}). We observe one path of the process $X$ on
the interval $[0, T_N]$  at the discrete times $t_i = i  \cdot
\Delta_N$ for $i=1,\dots, N$ with $T_N = N \cdot \Delta_N$. Therefore,
$$
(X({\Delta_N}),X({2\Delta_N}), \ldots, X({N\Delta_N}))
~~~~\mbox{is known,}
$$
and we consider the asymptotic  $N \to \infty$, $\Delta_N \to 0$
and $T_N \to \infty$. We want to estimate the parameters of the
$(M_K)$-F.B.M. that are $(H_0, H_1, \dots, H_K)$,
$(\sigma_0,\sigma_1, \dots, \sigma_K)$ and $( \omega_1, \dots,
\omega_K)$. \\
~\\
Even if the model is defined as a parametric one, we prefer to
use a semi-parametric statistics based on the wavelet analysis.
This choice is justified by the following reasons. First, the spectral density of $X$ is %%in the general case
not continuous in the general case. Thus, one cannot use the
classical results   on the consistency of  the maximum likelihood
or Whittle maximum likelihood estimators for long memory processes
(see Fox and Taqqu, 1986, Dahlhaus, 1989 or Giraitis and
Surgailis, 1990). Moreover, this is not a classical time series
parametric estimation~: indeed, we consider
$(X(\Delta_N),X(2\Delta_N),\ldots,X(N\Delta_N))$ instead of
$(X(1),X(2),\ldots,X(N))$ and therefore this is also an estimation
problem of the parameters of a continuous stochastic process.
Secondly, the following semi-parametric statistics are more robust
than a parametric one if the model is misspecified. Consider the
example where the function $H(\xi)$ is a not exactly a piece-wise
constant function, but instead a constant function on several
intervals and some unknown function on the other intervals. In
this case,  a parametric estimator could not work while the
semi-parametric method based on the wavelet analysis will remain
efficient. \\
~\\
Another semi-parametric method was developed from the seminal
paper of Istas and Lang (1997). This method of estimation is
derived from  the variogram and provides good results in the case
of F.B.M. (see Bardet, 2000) or of  multifractional F.B.M. (see
Benassi {\it et al.}, 1998). However, one faces difficulties in
identifying the model $(M_K)$-F.B.M. with this kind of method.
Indeed, one can easily  satisfy that for $\delta>0$~: \ba
\label{dec:var} \mathcal{V}(\delta) = \E\left ( X(t+\delta)
-X(t)\right )^2 &=& 4 \sum_{j=0}^{K} \delta^{2 H_j}\, \sigma_j^2\;
\int_{\delta \omega_j}^{\delta \omega_{j+1}} \frac {(1-\cos
v)}{v^{2H_j +1}} \, dv. \ea The principle of the variogram's
method ensues from the writing of $\log \Big (\mathcal{V}(\delta)
\Big )$ as an affine function of $\log \delta$. For a
$(M_K)$-F.B.M., with $\displaystyle{C(H_i)= \int _0 ^\infty \frac
{(1-\cos v)}{v^{2H_i +1}} \, dv }$ for $i=0,1,\ldots,K$, two cases
could provide such a relation~:
\begin{enumerate}
\item for $\delta \to \infty$,~~$\displaystyle{\log \Big (\mathcal{V}
(\delta) \Big )= 2H_0 \cdot \log \delta + \log \big ( 4 \cdot \sigma_0^2 \cdot
C(H_0) \big ) + O(\delta^{-2H_0})}$;
\item for $\delta \to 0$,~~$\displaystyle{\log \Big (\mathcal{V}(\delta)
\Big )= 2H_K \cdot  \log \delta + \log \big (4 \cdot \sigma_K^2 \cdot C(H_K)
\big ) + O(\delta^{2-2H_K})}$
\end{enumerate}
(the proof of such expansions is in the proof of Lemma
\ref{maj:cov}). In those cases, if one can show that there is a
convergent estimator $V_N(\delta)$ of $\mathcal{V}(\delta)$, then
a log-log regression of $\log \Big ( V_N(\delta)\Big )$ onto $\log
\delta$ could provide an estimation of the different parameters.
Nevertheless, such a method would have a lot of drawbacks. On one
hand, the estimation of "intermediate" parameters $(H_j)_{1\leq j
\leq K-1}$ and $(\sigma^2_j)_{1\leq j \leq K-1}$ requires very
specific asymptotic properties between all the frequency changes
$(\omega_j)_{1\leq j \leq K-1}$. This implies  a lack of generality
of the methods based on the variogram. Moreover, concretely, the
frequency changes are fixed and one obtains rough approximation
instead of asymptotic properties. For instance, numerical
simulations show that in some cases the log-log plot of the
variogram does not exhibit any intermediate linear part. On the
other hand, when the model is misspecified the variogram model
could lead to inadequate results. For example the following
picture gives the case of a $(M_2)$-F.B.M. where the variogram
method would detect only one frequency change and could not
precisely estimate its value. Finally, the variogram's method
could perhaps be applied in the two first previous situations 1.
and 2., {\em i.e.} for the estimation of $(H_0,\sigma_0^2)$ or
$(H_K,\sigma_K^2)$ with $\delta$ will have to be a function of $N$
(number of data). But this choice of function will depend on the
unknown parameters $H_0$ or $H_K$ for obtaining central limit
theorems for $\log \Big ( V_N(\delta)\Big )$... (see
the same kind of problem in Abry {\em et al.}, 2002). \\
\[
\epsfxsize 15cm \epsfysize 7.5cm \epsfbox{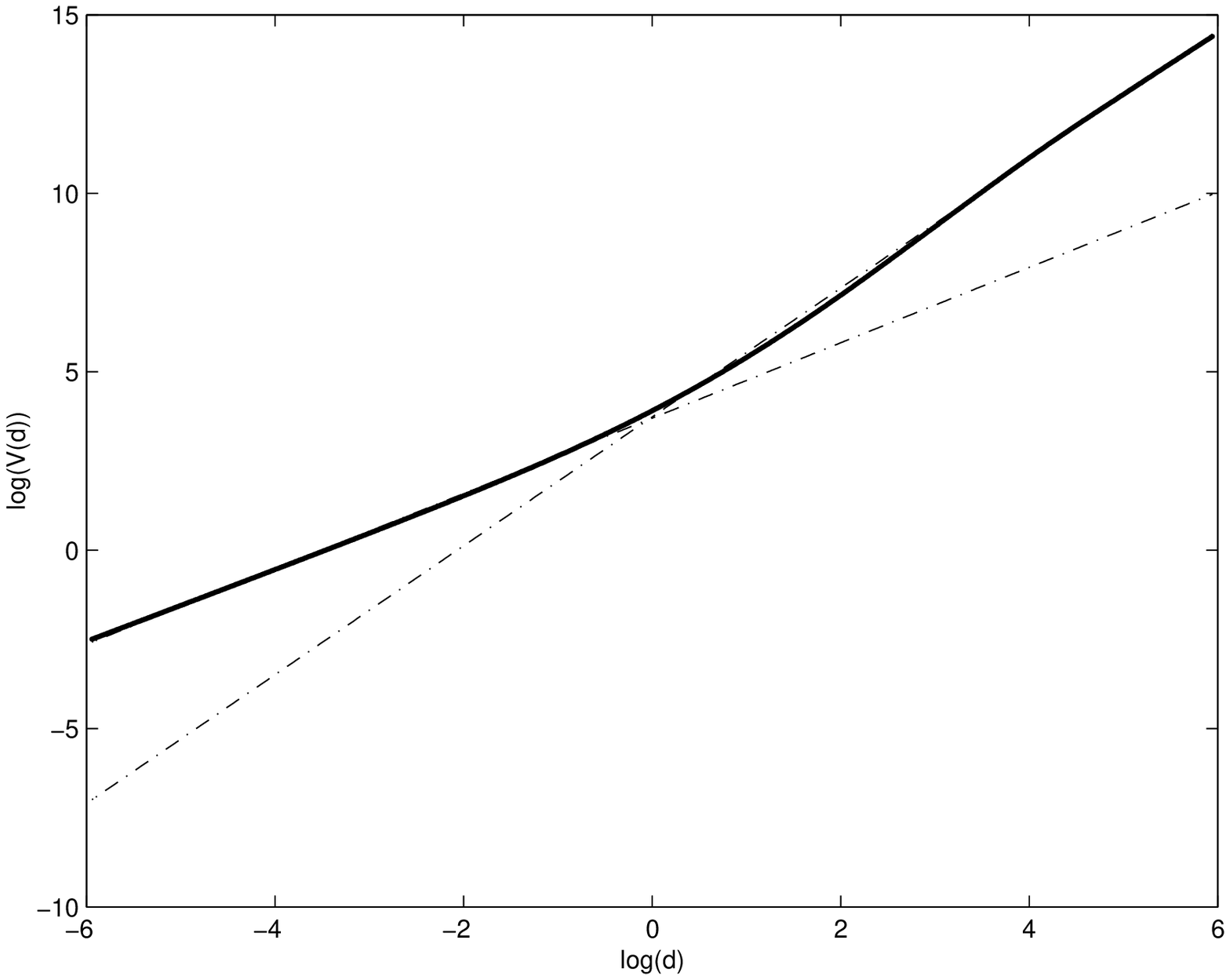}
\vspace{-0.5cm}
\]
{\bf Figure 3:} An example of a theoretical variogram for a
$(M_2)$-f.B.m, with $H_0=0.9$, $H_1=0.2$, $H_2=0.5$, and
$\sigma_0=\sigma_1=\sigma_2=5$ and $\omega_1=0.05$, $\omega_2=0.5$
(in solid, the theoretical variogram, in dot-dashed, its
theoretical asymptotes for $\delta \to 0$ and $\delta \to \infty$).  \\
~\\
We deduce from the definition of the model and the previous
discussion that a wavelet analysis could be an interesting
semi-parametric method for estimating the parameters of a
$(M_K)$-F.B.M. Indeed, such a method is based on the change of
scales (or frequencies). Therefore, as it is developed below, a
wavelet analysis is able to detect the different spectral domain
of self-similarity and then estimate the different parameters of
the model.
\subsection{A statistical study  based on wavelet analysis}
This method has been introduced by Flandrin (1992) and was
developed by Abry {\em et al.} (2002) and Bardet {\em et al.}
(2000). We also use in the following similar results on wavelet
analysis for $(M_K)$-F.B.M. obtained in Bardet and Bertrand
(2003).
Let $\psi$ be a %%``mother''
wavelet satisfying the following assumption~:
%%\\ \\ {\bf je sugg\`ere de supprimer "mother" dans la ligne ci-dessus. En effet, ceci est contradictoire avec la
 %%remarque que $\psi$ n'est pas obig\'e d'etre une ondelette m\`ere.}
\\ \\
{\bf  Assumption (A1):} {\em $\psi:~\R \mapsto \R$ is a ${\cal
C}^\infty$ function satisfying~:
\begin{itemize}
\item for all $m\in \R$, $\displaystyle{\int_{\R} \left |t^m\psi(t)\right |dt
<\infty}$;
\item its Fourier transform $\widehat{\psi}(\xi)$ is an even function
compactly supported on  $[-\beta,-\alpha]\cup[\alpha,\beta]$ with
$0 <\alpha<\beta$.
\end{itemize}}
\noindent We stress these conditions are sufficiently mild and are
satisfied in particular by the Lemari\'e-Meyer {\em "mother"}
wavelet. The admissibility property, i.e. $\displaystyle{\int_{\R}
\psi(t)dt = 0}$, is a consequence of the second one and more
generally, for all $m \in \N$,
\begin{eqnarray}\label{moments}
\int_{\R}t^m \psi(t)dt = 0.
\end{eqnarray}
%%Note that it is not mandatory to choose $\psi$ to be a ``mother'' wavelet associated to a multiresolution analysis
%%of $\LL^2(\R)$ and the whole theory can be developed without resorting to this assumption~:
Note that it is not necessary  to choose $\psi$ to be a {\em
"mother"} wavelet associated to a multiresolution analysis of
$\LL^2(\R)$. The whole theory can be developed without resorting
to this assumption. The choice of $\psi$ is then very large. ~\\
~\\
Let $(a,b)\in \R_+^*\times \R$ and denote $\lambda=(a,b)$. Then define
the family of functions $\psi_\lambda$ by
$\displaystyle{\psi_{\lambda}(t)=\frac 1{\sqrt a}\,
\psi\left(\frac{t}{a}-b \right)}$. Parameters $a$ and $b$ are
so-called the scale and the shift of the wavelet transform. Let us
underline that we consider a continuous wavelet transform. Let
$d_X(a,b)$ be the wavelet coefficient of the process $X$ for the
scale $a$ and the shift $b$, with
$$
d_X(a,b)
=\frac 1{\sqrt a} \int_{\R}  \psi(\frac{t}{a}-b)X(t)dt
=<\psi_{\lambda},X>_{L^2(\R)}.
$$
If $\psi$ satisfies Assumption (A1) and $X$ is a $(M_K)$-F.B.M.,
the family of wavelet coefficients verifies the following
properties (see Bardet and Bertrand, 2003)~:
\begin{enumerate}
\item for $a>0$, $(d_{X}(a,b))_{b\in \R}$ is a stationary
centered Gaussian process such as~:
\begin{eqnarray} \label{stoch}
\E\left(d^2_{X}(a,.) \right)=\mathcal{I}_{1}(a)=a \int_{\R}
|\widehat{\psi}(au)|^2\ \cdot  \rho^{-2}(u)\, du.
\end{eqnarray}
\item for all $i=0,1,\cdots,K$, if the scale $a$ is such as
$\displaystyle{[\frac \alpha a, \frac \beta a] \subset [\omega_i,\omega_{i+1}]}$, then
\begin{eqnarray}
\label{d(a)} \E\left(d^2_{X}(a,.) \right)=a^{2H_i+1}\cdot
\sigma_i^2 \cdot  K_{H_i}(\psi),~~\mbox{with}~
K_{H}(\psi)=\int_{\R} \frac { \left|\widehat{\psi}(u)\right|^2}{
|u|^{2H+1}}du.
\end{eqnarray}
\end{enumerate}
Property (\ref{d(a)}) means that the logarithm of the variance of
the wavelet coefficient is an affine function of the logarithm of
the scale with slope $2H_i+1$ and intercept $\log \sigma_i^2 +
\log K_{H_i}(\psi)$. This property is the key tool for  estimating
the parameters of $X$. Indeed, if we consider a convergent
estimator of $\log \left (\E\left(d^2_{X}(a,.) \right) \right)$,
it provides a linear model in $\log a$ and $\log \sigma_i^2$.
Before specifying such an estimator, let us stress that  one  only
observes a  discretized path
$(X(0),X(\Delta_N),\ldots,X(N\Delta_N))$ instead of a
continuous-time path. %%$(X(t))_{t\in\R}$.

As a consequence, for $a>0$ and $N
\in \N^*$, a natural estimator is the logarithm of
the empirical variance of the wavelet coefficient, that is $\log
I_N(a)$ where~: \be \label{def:IN(a)} I_N(a) = \frac{1}{|D_N(a)|}
\sum_{k\in D_N(a)} d_{X}^2(a,k\Delta_N),\ee with~:
\begin{itemize}
\item $r \in ]0,1/3[$;
\item $m_N=\left[ r(N/a) \right ]$ and $M_N=\left[ (1-r)(N/a) \right ]$
where $[x]$ is the integer part of $x \in \R$;
\item $D_N(a)=\{m_N,m_N+1,\ldots,M_N \}$ and $|D_N(a)|$ is the cardinal of the set
$D_N(a)$.
\end{itemize}
For $0<a_{min}<a_{max}$, a functional central limit theorem for $(
\log I_N(a))_{a_{min} \leq a \leq a_{max}}$ can be established
~(see a similar proof in Bardet and Bertrand, 2003) :
\begin{Prop} \label{propTLC}
Let $X$ be a $(M_K)$-F.B.M., $0<a_{min}<a_{max}$ and  $\psi$
satisfy Assumption (A1). Then~:
\begin{eqnarray} \label{TLC}
\sqrt {N\Delta_N} \left ( \log I_N(a)-\log  \mathcal{I}_{1}(a)
\right )_{a_{min} \leq a \leq a_{max}} \limiteloi (Z(a))_{a_{min}
\leq a \leq a_{max}}
\end{eqnarray}
with $(Z(a))$ a centered Gaussian process such as for $(a_1,a_2)
\in [a_{min},a_{max}]^2$, \be \label{def:skl} \cov
(Z(a_1),Z(a_2))= \frac {2a_1 \,a_2}{(1-2r)\,
\mathcal{I}_{1}(a_1)\, \mathcal{I}_{1}(a_2)} \int _{\R} \left(
\int _{\R} \frac {\overline{\widehat{\psi}}(a_1
\xi)\widehat{\psi}(a_2 \xi)}{|\rho(\xi)|^2}e^{-iu\xi}d \xi \right
)^2du. \ee
\end{Prop}
Then, if we specify   the locations of the change points in terms
of scales, {\em i.e.} frequencies, we obtain the  following:
\begin{cor} \label{corTLC}
Let $i \in \{0,1,\cdots,K\}$ and assume that $\displaystyle{\frac
\beta \alpha \leq \frac {\omega_{i+1}}{\omega_i}}$. Then,
\begin{eqnarray} \label{TLCprecis}
\nonumber \hspace{-1.5cm} \sqrt {N\Delta_N} \left ( \log I_N(
1/f)\hspace{-1mm} +\hspace{-1mm}(2H_i\hspace{-1mm}
+\hspace{-1mm}1)\log f\hspace{-1mm} -\hspace{-1mm}\log
\sigma_i^2\hspace{-1mm} -\hspace{-1mm}\log K_{H_i}(\psi)
 \right )_{\omega_i/\alpha \leq f \leq \omega_{i+1}/\beta}\hspace{-2mm}&& \\
 \limiteloi \hspace{-2mm} (Z(1/f))_{\omega_i /\alpha \leq f
 \leq \omega_{i+1}/\beta} \hspace{-1.4cm} & &
\end{eqnarray}
with the centered Gaussian process $(Z(.))$ such as for
$\displaystyle{(f_1,f_2) \in [\frac {\omega_i} \alpha, \frac
{\omega_{i+1} } \beta]^2}$, \be \label{def:sklprecis} \cov
(Z(1/f_1),Z(1/f_2))=\frac {2\left(f_1\, f_2 \right)^{2H_i}}{(1-2r)
\,  K^2_{H_i}(\psi)} \int _{\R} \left( \int _{\R} \frac
{\overline{\widehat{\psi}}(  \xi/f_1)
\widehat{\psi}(  \xi /f_2)}{|\xi|^{2H_i+1}}e^{-iu\xi}d \xi \right )^2du.\\
\ee
\end{cor}
For $\Delta_N$ small enough, this result shows that all parameters
$H_i$ and $\sigma_i^2$ could be estimated by using a linear
regression of $\log I_N( 1/f_j)$ versus $\log f_j$, when the
frequencies $\omega_i$ are known. Moreover, this central limit
theorem shows that a graph of ($\log f,\,\log I_N( 1/f)$) for
$f>0$ exhibits different areas of asymptotic linearity~: it
suggests the procedure of the following section to estimate and
test the frequency changes (see for instance figures 4 or 6).
\subsection{The discretization problem}
In the applications, we only observe a finite  time series $( X
(0), X (\Delta_N),\cdots, X ((N-1)\times\Delta_N))$ and we must derived the
empirical wavelet coefficients from this time series. Since the
process $X$ has almost a continuous path but with a regularity
$\alpha_X <1 $ almost surely, we should use the Riemann sum. Thus,
for $(a,b)\in \R_+^*\times \R$ we define the empirical wavelet
coefficient by \be \label{def:eab} e_{X}(a,b)=\frac
{\Delta_N}{\sqrt a} \sum _{p=0}^{N-1}\psi ( \frac {p \Delta_N} a
-b)\times X (p\Delta_N)\ee and the discretized estimator by
 \be \label{def:JN(a)} J_N(a) =
\frac{1}{|D_N(a)|} \sum_{k\in D_N(a)} e^2_{X}(a,k\Delta_N). \ee We also
define for every $k \in D_N(a)$ the   error
\begin{eqnarray}\label{def:epsilon}
\varepsilon_N(a,k) = e_{X}(a,k\Delta_N) - d_{X}(a,k\Delta_N).
\end{eqnarray}
Now, it is possible to provide the functional central limit
theorem for $( \log J_N(a))_{a_{min} \leq a \leq a_{max}}$
computed from $( X (0), X (\Delta_N),\cdots, X (N\Delta_N))$:
\begin{Theo} \label{Theo:TCL:discret}
Under assumptions of Proposition \ref{propTLC} and with $\Delta_N$
such as $N\Delta_N\to\infty$ and $N (\Delta_N)^2\to 0$ when $N\to
\infty$. Then, with the same process $Z$ than in (\ref{TLC}),
\begin{eqnarray} \label{TLC2}
\sqrt {N\Delta_N} \left ( \log J_N(a)-\log \mathcal{I}_1(a) \right
)_{a_{min} \leq a \leq a_{max}} \limiteloi (Z(a))_{a_{min} \leq a
\leq a_{max}}.
\end{eqnarray}
As a particular case, for $i \in \{0,1,\cdots,K\}$ and if
$\displaystyle{\frac \beta \alpha \leq \frac
{\omega_{i+1}}{\omega_i}}$, then
\begin{eqnarray}  \label{TLCdis}
\nonumber \hspace{-1.5cm} \sqrt {N\Delta_N} \left ( \log J_N(
1/f)\hspace{-1mm} +\hspace{-1mm}(2H_i\hspace{-1mm}
+\hspace{-1mm}1)\log f\hspace{-1mm} -\hspace{-1mm}\log
\sigma_i^2\hspace{-1mm} -\hspace{-1mm}\log K_{H_i}(\psi)
 \right )_{\omega_i/\alpha \leq f \leq \omega_{i+1}/\beta}\hspace{-2mm}&& \\
\limiteloi \hspace{-2mm} (Z(1/f))_{\omega_i /\alpha \leq f \leq
\omega_{i+1}/\beta}. \hspace{-1.4cm} & &
\end{eqnarray}
\end{Theo}
The convergence rate of the central limit theorem (\ref{TLC2}) is
$\sqrt{N\Delta_N}$. Thus, the discretization problem implies that
the maximum convergence rate is $o(N^{1/4})$ from the previous
conditions on $\Delta_N$.
\section{Identification of the parameters}
First, let us describe the method on a heuristic level. From
Proposition \ref{propTLC}, Formula (\ref{TLCdis}), we have
\begin{eqnarray} \label{TLC3}
\log J_N( 1/f)=-(2H_i+1)\times\log(f)+\log\left( \sigma_i^2\right)
+\log \left(K_{H_i}(\psi)\right) + \varepsilon_f^{(N)},
\end{eqnarray}
for the frequencies $f$ which satisfy the condition \be
\label{cond2} \log\left(\omega_{i}\right) -\log(\alpha) \le
\log\left(f\right) \le \log\left(\omega_{i+1}\right) -\log(\beta).
\ee Moreover we have $ \displaystyle{
\left({N\Delta_N}\right)^{1/2} \left(
\varepsilon_{f_j}^{(N)}\right )_{1\leq j \leq m}\limiteloi
(Z(1/f_j))_{1\leq j \leq m}.}$ Formula (\ref{TLC3}) and condition
(\ref{cond2}) mean that for $\,\displaystyle{\log(f) \in
[\log\left(\omega_{i}\right) -\log(\alpha),
\log\left(\omega_{i+1}\right)}$ $-\log(\beta)]$, we have a linear
regression of $\;\log J_N(1/f)\;$ onto $\;\log(f)\;$ with slope
$\;-(2H_i+1)\;$ and intercept\\
 $\;\log \sigma_i^2 +\log
K_{H_i}(\psi)\;$ and for $\,\log(f) \in
\left[\log\left(\omega_{i+1}\right) -\log(\alpha),
\log\left(\omega_{i+2}\right) -\log(\beta)\right]\;$ a linear
regression with slope $\;-(2H_{i+1}+1)\;$ and intercept $\;\log
\sigma_{i+1}^2 +\log K_{H_{i+1}}(\psi)$. This is  a problem of
detection of abrupt change on the parameters of a linear
regression, but with a transition zone for $\;\log(f) \in
\left]\log\left(\omega_{i+1}\right) -\log(\beta),\;
\log\left(\omega_{i+1}\right) -\log(\alpha) \right[$.
\begin{rem}
Condition (\ref{cond2}) implies that $\displaystyle{\omega_{i+1} >
\frac{\beta}{\alpha}\times \omega_{i}}$. Therefore we could only detect
the frequency changes sufficiently spaced. For instance, if we
choose the Lemari\'e-Meyer wavelet, we get $\beta/\alpha = 4$ which
leads to the condition $\omega_{i+1} > 4\times \omega_{i}$.
\end{rem}
In this section, we describe the estimation of the parameters and
a goodness of fit test. Both of them are based   on the following
assumption :
\\
\\ {\bf Assumption $(B_K)$ : } {\sl The process  $X$ is a
$(M_K)$-multiscale fractional Brownian motion. This process is
characterized by the parameters   $\Omega^*$, $H^*$ and $\sigma^*$
where $\Omega^* = \left(\omega_1^*, \cdots, \omega_K^* \right)$ with $H^*
= \left(H_0^*, H_1^*,\dots,H_K^*\right)$ and $\sigma^*=\left(
\sigma_0^*, \sigma_1^*, \dots,\sigma_K^* \right)$. Moreover the
following conditions are fulfilled
\begin{itemize}
\item  $\displaystyle{\omega_{i+1}^* > \frac{\beta}{\alpha}\times \omega_{i}^*}$
for $i=1,\cdots,K-1$;
\item $\displaystyle{\min _{0\leq i  \leq (K-1)}
\left \{ \Big(H_{i+1}^*-H_i^*\Big)^2+\Big(\sigma_{i+1}^*-
\sigma_i^*\Big)^2 \right \}>0}$ and
\item there exists a compact set ${\cal K} \subset ]0,1[
\times ]0,\infty[$ such as $(H_i^*,\sigma_i^*) \in {\cal K}$ for all
$i=0,1,\cdots,K$.
\end{itemize}
}
\subsection{Estimation of the parameters}
Let $X$ be a $(M_K)$-F.B.M. satisfying the assumption $(B_K)$ with
$K$ a known integer number. We observe one path of the process at
$N$ discrete times, that $( X (0), X (\Delta_N),\cdots, X
(N\Delta_N))$. Let $\left[f_{min},\,f_{max}\right]$, with
$0<f_{min}<f_{max}$, be the chosen frequency band (see section 5,
for an example). We discretize a (slightly modified) frequency
band and compute the wavelet  coefficients at the frequencies
$(f_k)_{0\leq k \leq a_N}$ where \ba \nonumber f_k=\frac {f_{min}}
\beta \, (q_N)^k&&\mbox{for}~k =0,\cdots,a_N,\quad q_N=  \left ( \frac
{f_{max}}{f_{min}} \frac \beta \alpha \right )^{1/a_N}\quad
\mathrm{and}\quad a_N=N\Delta_N  . \ea For notational convenience
, we assume here that $N \Delta_N$ is an integer number.  By
definition, we have $f_0=f_{min}/\beta$ and $f_{a_N}
=f_{max}/\alpha$, then, using the wavelet coefficients at the
frequencies $(f_k)_{0\leq k \leq a_N}$, we could detect   all
frequency changes $(\omega_i^*)$  included in the band
$]f_{min},f_{max}[$. To simplify the notations, we use the
following assumption :
\par
\medskip
\noindent {\bf Assumption (C) : } {\sl $\omega_i^* \in
]f_{min},f_{max}[$ for all $i=1,\dots, K$. }
\par
\medskip
\noindent In this framework, the estimation of the different
parameters of $X$ becomes a problem of linear regression with a
known number of changes; thus, we follow the same method as in Bai
(1994), Bai and Perron (1998), Lavielle (1999) or Lavielle and
Moulines (2000) and define the estimated parameters
$(\widehat{T}^{(N)},\widehat{\Lambda}^{(N)})$ as the couple of
vectors which minimize the quadratic criterion~:
$$Q^{(N)}(T,\Lambda)= \sum_{j=0}^{K+1}\,\sum_{i=1+t_j}^{t_{j+1}-\tau_N}
\left|Y_i-X_i \lambda_j \right|^2,\mbox{and thus }$$
$$(\widehat{T}^{(N)},\widehat{\Lambda}^{(N)})=\mbox{Argmin}
\left \{ Q^{(N)}(T,\Lambda);~ T \in {\cal A}_K^{(N)},\Lambda \in
{\cal B}_K \right \} $$ with
\begin{itemize}
\item $Y_i= \log \left (J_N(1/f_i) \right ) $, $X_i=( \log f_i,1)$ for $i=0,\cdots,a_N$;
\item $\displaystyle{\tau_N= \left [  \frac { \log ( \beta /\alpha) }{\log q_N}  \right ] }$,
where $[x]$ is the integer part of $x$.
\item $T=(t_0,t_1, \cdots,t_{K+1})\in {\cal A}_K^{(N)}$ where
$$\hspace{-1cm}{\cal A}_K^{(N)}\hspace{-1mm}=\hspace{-1mm}
\left \{ (t_0,\cdots,t_{K+1}) \in \N^{K+2};\hspace{-1mm}~
t_0=0,\hspace{-1mm} ~t_{K+1}=a_N+\tau_N,\hspace{-1mm}~t_{j+1}-t_j
> \tau_N~\mbox{for}~j=0,\cdots,K \right \};$$
\item $\Lambda=(\lambda_0,\cdots,\lambda_K) \in {\cal B}_K$ where
$\displaystyle{\lambda_j=\left ( \begin{array}{c}
-(2H_j+1) \\
\log \sigma^2_j + \log K_{H_j}(\psi)
\end{array} \right ) }$ and then
$$\hspace{-1cm}{\cal B}_K \hspace{-1mm}=\hspace{-1mm}
\Big \{(\lambda_0,\cdots,\lambda_K) ~~\mbox{with}~~(H_j,\sigma^2_j)\in
{\cal K}~~\mbox{for all}~~j\in \{0,1,\cdots,K\} \Big \}.$$
\end{itemize}
The integer $\tau_N$ corresponds to  the number of frequencies in
the transition zones and  $\log f_{i+\tau_N}=\log f_i + \log
(\beta/\alpha)$. Obviously, for $j=0,\cdots,K$, the vector
$\widehat{\lambda}_j^{(N)}$ provides the estimators
$\widehat{H}_j^{(N)}$ of $H_j^*$ and $\widehat{\sigma}_j^{(N)}$ of
$\sigma_j^*$  by the relation
$\displaystyle{\widehat{\lambda}_j^{(N)}=\left ( \begin{array}{c}
- (2 \widehat{H}_j^{(N)} +1) \\
\log \Big ( (\widehat{\sigma}_j^{(N)})^2\Big )  + \log
K_{\widehat{H}_j^{(N)}}(\psi)
\end{array} \right ) } $.
For a given $T \in {\cal A}_K^{(N)}$, each
$\widehat{\lambda}_j^{(N)}$ is obtained from a linear regression
of $(Y_i)$ onto $(X_i)$ for $i=t_j+1,\cdots,t_{j+1}-\tau_N$. Thus,
with $\widehat{T}=(\widehat{t_j})_{0\leq j \leq K+1}$ obtained
from the minimization in $T$ of $Q^{(N)}(T,\widehat{\Lambda})$, we
define the different estimators of the change frequencies as \be
\label{def:omega:chapeau} \widehat{\omega}_j^{(N)}=\alpha \,
f_{\widehat{t}_j^{(N)}}=\alpha \cdot  \frac {f_{min}} \beta \left
(
\frac {f_{max}}{f_{min}} \frac \beta \alpha \right )  ^{\frac
{\widehat{t}_j^{(N)}}{a_N}}~~\mbox{for}~j=1,\cdots,K. \ee We have
the following convergence :
\begin{Prop} \label{omegas}
Let $X$ satisfy Assumptions (C) and ($B_K$) with a known $K$,
$(X_{\Delta_N}, \cdots, X_{N\Delta_N})$ be a discretized path, and
$\psi$ satisfy Assumption (A1). Let $\Delta_N$ be such as
$N\Delta_N \to \infty$ and $N(\Delta_N)^2 \to 0$ when $N\to
\infty$. Assume that
$(\widehat{H}_i^{(N)},\widehat{\sigma}_i^{(N)}) \in {\cal K}$ for
all $i=0,\cdots,K$. Then for all $\varepsilon>0$, there exists
$0<C<\infty $ such as for all large $N$, \ba \Pr
\label{conv_omega}
\left ( (N\Delta_N )^{1/4}\left | \widehat{\omega}_j^{(N)}
-\omega_j^*\right |\geq C\right ) \leq
\varepsilon~~\mbox{for}~j=1,\cdots,K. \ea
\end{Prop}
\begin{rem}
The proof of this proposition shows a more general result, {\em
i.e.} for $(p,q) \in [3/4,1] \times [0,1]$, for $\varepsilon >0$,
there exists $C>0$ such as
$$\Pr \left ( a_N^{1-p}\left | \widehat{\omega}_j^{(N)} -\omega_j^*\right |
\geq C\right )  \leq \varepsilon~~\mbox{for}~j=1,\cdots,K$$ with
$a_N=(N\Delta_N)^q$. For numerical considerations and convergence
rate of the following estimators of the parameters, we are going
to fix now on $p=3/4$ and $q=1$ and then $a_N=N \Delta_N.$
\end{rem}
For $j=0,\cdots,K$, the natural estimates of $H_j^*$ and
$\sigma_j^{2*}$ are given by  the regression of $(Y_i)$ onto
$(\log f_i)$ for $i \in \{ \widehat{t}_j^{(N)},\cdots,
\widehat{t}_{j+1}^{(N)}-\tau_N \}$. But the probability that
$[\widehat{t}_j^{(N)},\widehat{t}_{j+1}^{(N)}-\tau_N] \subset
[t_j^*, t_{j+1}^* -\tau_N]$ does not increase fast enough to $1$
as $N \to \infty$, in order to obtain a sufficiently fast
convergence rate for these estimators. We address this difficulty
as follows. We fix an integer number $m\geq 3$ and for
$j=0,\cdots, K$, we consider $[\tilde{U}_j^{(N)},\tilde{V}_j^{(N)}
]$ an interval strictly included in
$[\widehat{t}_j^{(N)},\widehat{t}_{j+1}^{(N)}-\tau_N]$,  such as
\ba \label{UV} \tilde{U}_j^{(N)}=\widehat{t}_j^{(N)}+ \left [\frac
{\widehat{t}_{j+1}^{(N)}-\widehat{t}_{j}^{(N)}-\tau_N }{m+1}
\right ]~~\mbox{and}~~\tilde{V}_j^{(N)}=\widehat{t}_j^{(N)}+ m
\left [\frac {\widehat{t}_{j+1}^{(N)}-\widehat{t}_{j}^{(N)}-\tau_N
}{m+1} \right ]. \ea Then we estimate the parameters from a
regression onto $m$ points uniformly distributed in
$[\tilde{U}_j^{(N)},\tilde{V}_j^{(N)} ]$; it provides the
following estimator $\tilde{\lambda}_j^{(N)}$ from a regression of
$(Y_i)$ onto $(X_i)$ for \\
$\displaystyle{i \in \{\tilde{U}_j^{(N)},\cdots,\tilde{V}_j^{(N)}
\}=
\left \{\tilde{U}_j^{(N)}+(k-1) \left [\frac
{\widehat{t}_{j+1}^{(N)}-\widehat{t}_{j}^{(N)}-\tau_N }{m+1}
\right ] \right  \}_{1\leq k \leq m}}$. By this way, define \ban
\tilde{\lambda}_j^{(N)}&= &\Big(- (2 \tilde{H}_j^{(N)} +1),
 \log \tilde{\sigma^2}_j^{(N)} + \log K_{\tilde{H}_j ^{(N)}}(\psi)\Big )'\\
&= &\left (  (\tilde{X}^{(N)}_j)'  \tilde{X}^{(N)}_j  \right
)^{-1}( \tilde{X}^{(N)}_j )'\tilde{Y}^{(N)}_j~~\mbox{with}~~\left
\{
\begin{array}{l} \tilde{X}^{(N)}_j=\displaystyle{\left ( \log f_i~,~1
\right )_{i \in
\{\tilde{U}_j^{(N)},\cdots,\tilde{V}_j^{(N)} \}}} \\
\tilde{Y}^{(N)}_j=(Y_i)_{i \in
\{\tilde{U}_j^{(N)},\cdots,\tilde{V}_j^{(N)} \}} \end{array} \right .
, \ean and for all $k=1,\cdots,m$, define $\displaystyle{ g_0^*(k)
=
\frac {f_{min}} {\beta}\left ( \frac { \omega_{1}^*} {f_{min}}
\right ) ^{k/(m+1)} }$, $\displaystyle{ g_K^*(k) = \frac
{\omega_K^*} {\alpha}\left ( \frac {f_{max}} {f_{min}} \right )
^{k/(m+1)} }$ and \\
$\displaystyle{ g_j^*(k) = \frac {\omega_j^*} {\alpha}\left (
\frac {  \alpha \omega_{j+1}^*}
{\beta\omega_j^*} \right ) ^{k/(m+1)} }$ for all $j\in\{1,\cdots,K-1\}$, . \\ \\
We get the following central limit theorems for the corresponding
estimators $(\tilde{H}_j^{(N)},\tilde{\sigma^2}_j^{(N)})$~:
\begin{Prop}\label{HetK}
Under the same assumptions as in Proposition \ref{omegas}, for all
$j=0,\cdots,K$, \ba \label{conv_lambda} (N\Delta_N )^{1/2}\left
(\tilde{\lambda}_j^{(N)}-\lambda_j^*\right ) & \limiteloi& {\cal
N}(0,\Gamma_1^{\lambda_j^*}) \ea where
$\displaystyle{\Gamma_1^{\lambda_j^*}=\left (  X^{*'}_j  X^{*}_j
\right )^{-1}X^{*}_j\Sigma_j^* X^{*'}_j\left (  X^{*'}_j X^{*}_j
\right )^{-1}}$, with $\displaystyle{X^{*}_j=\left ( \log
g_j^{*}(k)~,~1 \right )_{1\leq k \leq m}}$ and
$\Sigma_j^*=(s^{*j}_{kl})_{1 \leq k,l \leq m}$ the following
matrix~:
%and $\Gamma^\sigma=(\gamma^\sigma_{ij})_{1\leq i,j \leq K+1}$ two
%a diagonal matrix such as~:
\ba \label{Sigmaj} s^{*j}_{kl}= 2 \cdot \Big( g_j^{*}(k)
g_j^{*}(l)\Big )^{2H^*_j}\cdot  \frac {\displaystyle{  \int _{\R}
\left( \int _{\R} \overline{\widehat{\psi}} \left (  \frac \xi {
g_j^{*}(k)} \right ) \widehat{\psi} \left (  \frac \xi {
g_j^{*}(l)} \right ) |\xi|^{-(2H_j^*+1)}e^{-iu\xi}d \xi \right
)^2du}} {\displaystyle{\left ( \int_{\R}
\left|\widehat{\psi}(u)\right|^2 |u|^{-(2H_j^*+1)}du \right )^2 }}
. \ea
\end{Prop}
\begin{rem}
Another possible choice would be to consider the regression for
all the available frequencies in the interval
$[\tilde{U}_j^{(N)},\tilde{V}_j^{(N)} ]$. The number of considered
frequencies increases then with the rate $a_N=N \Delta_N$.
However, it does not  improve significantly   the convergence
since the remainders of the regression are very strongly
dependent.
\end{rem}
\subsection{Goodness of fit test}
It is also possible to estimate parameters $H_j^*$ and
$\sigma_j^*$ from an feasible (or estimated) generalized least
squares estimation (for more details, see Amemiya, chap. 6.3,
1985). Indeed, we can identify the asymptotic covariance matrix
$\Sigma_j^*$ for $j=0,\cdots,K$~: this matrix has the form
$\Sigma_j^*=\Sigma(H_j^*,\omega_j^*,\omega_{j+1}^*)$ and, from the
previous limit theorems,
$\widehat{\Sigma}_j^{(N)}=\Sigma(\tilde{H_j}^{(N)},
\widehat{\omega}_j^{(N)},\widehat{\omega}_{j+1}^{(N)})$ converges in
probability to $\Sigma_j^*$. Thus, it is possible to construct an
estimator $\underline{\lambda}_j^{(N)}$ of $\lambda_j^*$ with a
feasible generalized least squares (F.G.L.S.) regression {\em
i.e.} by minimizing
$$
\parallel  \tilde{Y}^{(N)}_j- \tilde{X}^{(N)}_j
\lambda \parallel ^2_{\widehat{\Sigma}_j^{(N)}}
=(\tilde{Y}^{(N)}_j- \tilde{X}^{(N)}_j\lambda)' \left (
\widehat{\Sigma}_j^{(N)}\right ) ^{-1}(\tilde{Y}_j^{(N)}-
\tilde{X}^{(N)}_j \lambda).
$$
First, we give asymptotic behavior of
$\displaystyle{\underline{\lambda}_j^{(N)}= \left \{
\begin{array}{c}
\left (- (2 \underline{H}_j^{(N)} +1), \log
\underline{\sigma^2}_j^{(N)}
+ \log K_{\underline{H}_j ^{(N)}}(\psi) \right )'\\
\left (  (\tilde{X}^{(N)}_j)' \left (
\widehat{\Sigma}_j^{(N)}\right ) ^{-1} \tilde{X}^{(N)}_j  \right
)^{-1}( \tilde{X}^{(N)}_j)'  \left (
\widehat{\Sigma}_j^{(N)}\right ) ^{-1} \tilde{Y}^{(N)}_j
\end{array} \right .} $.
%%  and for the corresponding estimators
%% $(\underline{H}_j^{(N)},\underline{\sigma^2}_j^{(N)})$.
\begin{Prop}\label{HetK2}
Under the same assumptions as in Proposition \ref{HetK}, for all
$j=0,\cdots,K$, \ba (N\Delta_N )^{1/2}\left
(\underline{\lambda}_j^{(N)}-\lambda_j^*\right ) & \limiteloi&
{\cal N}(0,\Gamma_2^{\lambda_j^*}) \ea with
$\displaystyle{\Gamma_2^{\lambda_j^*}=\left (  X^{*'}_j \left (
\Sigma_j^{*}\right ) ^{-1} X^{*}_j \right )^{-1}}$.
\end{Prop}
For $j=0,\cdots,K$, the vectors $ \tilde{Y}_j ^{(N)}$ and
$\tilde{X}^{(N)}_j \underline {\lambda}_j^{(N)}$ are two different
estimators of the vector\\
$\displaystyle{\left ( -(2H_j^* +1)\log f_i+\log \sigma_j^{2*}
+\log K_{H_i^*}(\psi)\right ) _{i \in
\{\tilde{U}_j^{(N)},\cdots,\tilde{V}_j^{(N)} \}} } $. It suggests to
define the following goodness of fit test. The test statistic
$T_K^{(N)}$ is defined as the sum of the squared distances between
these two estimators for all $K+1$ frequency ranges:
$$
T_K^{(N)}=(N \Delta_N )\cdot \left ( \sum _{j=0}^K \parallel
\tilde{Y}^{(N)}_j -\tilde{X}^{(N)}_j \underline{\lambda}_j^{(N)}
\parallel ^2_{\widehat{\Sigma}_j^{(N)}}\right ) .
$$
This distance is the F.G.L.S. distance between points $(\log
f_i,Y_i)_{i \in \{\tilde{U}_j^{(N)},\cdots,\tilde{V}_j^{(N)} \}}$
for $j=0,\cdots,K$ and the $(K+1)$ F.G.L.S. regression lines. As a
consequence, we get
\begin{Prop} \label {test}
Under assumptions of Proposition \ref{omegas}, we have
\begin{eqnarray} \label{Test}
T_K^{(N)} \limiteloi {\chi ^2 ((K+1)(m-2))}.
\end{eqnarray}
\end{Prop}
\begin{rem}
Proposition \ref{test} may be explained with heuristic arguments.
Remainders are turned white, thus it is only natural for the sum
of the second regression remainder squares to asymptotically form
a $\chi^2$ process. The number of degrees of freedom   is
$(K+1)(m-2)$ because one loses two degrees of freedom   after the
twice estimation of the $(K+1)$ vectors $\lambda_j^*$  (we also
show that these vectors are asymptotically independent).
\end{rem}
\subsection{Estimation of the number of frequency changes}
Throughout the previous study, the number of frequency change,
$K$, is assumed to be known. But the previous test provides a way
for estimating $K$. In fact, it can be recursively done by
beginning with $K=0$ and continuing till the assumption ``$X$ is a
$(M_K)$-F.B.M.'' is accepted. The following applications in
biomechanics provide different examples of the power of
discrimination of such a procedure. However, this estimation of
the number of frequency changes must be carefully applied~: from
numerical and heuristic arguments, it does not seem reasonable to
work with $K>2$.
\subsection{Estimation procedure and on the choice  of parameters}
Thus, for identifying a $(M_K)$-multiscale fractional (with $K$
unknown) from a time series
$(X_0,X_{\Delta_N},\cdots,X_{N\Delta_N})$ we suggest the following
procedure:
\begin{enumerate}
\item Begin with $K=0$.
\item Choose a mother wavelet $\psi$ (and thus $\alpha$ and $\beta$),
a frequency band $[f_{min},f_{max}]$ and $m$ (see below for these
different choices).
\item Compute the different frequencies $(f_i)_{0\leq i \leq a_N}$.
\item Compute the vector $(Y_i)_{0\leq i \leq a_N}=(\log J_N(1/f_i))_{0 \leq i \leq a_N}$.
\item Minimize $Q^{(N)}(T,\Lambda)$ and thus compute the different
values of $\widehat{\omega}_j^{(N)}$ for $j=1,\cdots,K$.
\item Compute the different regression moments $\{\tilde{U}_j^{(N)},\cdots,\tilde{V}_j^{(N)}\}$
 and then the estimators $\tilde{\lambda}_j^{(N)}$ (for $j=0,\cdots,K$).
\item Compute the different matrices $\widehat{\Sigma}_j^{(N)}$ and then
$\underline{\lambda}_j^{(N)}$ (for $j=0,\cdots,K$).
\item Compute $T_K^{(N)}$ and compare its value to the $95 \%$-quantile of a
$\chi ^2 ((K+1)(m-2))$. If the test is rejected then go back to
step 2. with $K=K+1$.
\end{enumerate}
How to chose the function  $\psi$ and the parameters $f_{min}$,
$f_{max}$ and $m$~?
\begin{enumerate}
\item  {\bf Choice of $\psi$ :} The mother wavelet $\psi$ has to satisfy
Assumptions (A1) but as we say previously it is not mandatory to
associate this function to orthogonality properties. However, the
Lemari\'e-Meyer wavelet is a natural choice with good numerical
properties of asymptotic decreasing but a too large ratio
$\beta/\alpha$ which implies a too large transition zone of
frequencies. The function $\psi$ can also be deduced from an
arbitrary construction of its Fourier transform $\widehat \psi$;
for instance, we propose $\displaystyle{{ \widehat \psi}_1
(\lambda)= \exp \left ( \frac {-1}
{(|\lambda|-\alpha)(\beta-|\lambda|)} \right )\1 _{\alpha \leq
|\lambda| \leq \beta} }$ and the function $\psi_2$ built from a
translation of the Fourier transform of the Lemari\'e-Meyer
function to $[-2\pi,-\pi] \cup [\pi,2\pi]$ (thus the ratio is now
$\beta/\alpha=2$). The results obtained from those functions
$\psi_1$ and $\psi_2$ are essentially the same than with the
Lemari\'e-Meyer mother function, they appear more precise for the
detection of frequency changes $\omega_j^*$ (because $\log
\beta/\alpha$ and thus the transition band, could be as small as
wanted) and less precise for the estimation of parameters $H_j^*$
(because $\psi_1$ and $\psi_2$ are not concentrated as well around
$0$).
\item  {\bf Choice of $f_{min}$ and $f_{max}$ :}
(we assume here that the frequencies are given in the inverse of
$(X_1,X_2\cdots)$ time unity). The choice of  $f_{min}$ and
$f_{max}$ is first driven by the selection of a frequency band
inside which the process has to be studied; the inspected
frequency band is then $\displaystyle{[\frac { f_{min}
}\beta,\frac { f_{min} }\alpha]}$. Secondly,
$\displaystyle{N\times \frac { f_{min}}\beta}$ should be large
enough for computing $\displaystyle{I_N(\frac
\beta {f_{min}})}$ in (\ref{def:IN(a)}). Formally one only needs
to have $\displaystyle{N\times \frac { f_{min}}\beta \geq 1}$ but
numerically $\displaystyle{N\times\frac {f_{min}} \beta \geq 10}$
seems to be necessary to use correctly the central limit theorem.
Finally, the discretization problem implies that $f_{max}$ cannot
be too large for providing a good estimation of
$\displaystyle{d_{X}(\frac \alpha {f_{max}},k\Delta_N)}$ by
$\displaystyle{e_{X}(\frac \alpha {f_{max}},k\Delta_N)}$. In
practice $\displaystyle{\frac {f_{max}} \alpha \leq \frac 1
{\Delta_N}}$ appears as a minimal condition.
\item {\bf Choice of $m$ :} Formally, $m$ could be chosen such as
$3 \leq m < \min _j (t_{j+1}^* -\tau_N -t_{j}^*)$. Theoretically,
the larger the $m$, the closer to $1$ the power of the test. But
numerical considerations imply that if $m$ is too large then the
different matrix $\widehat{\Sigma}_j^{(N)}$ are extremely
correlated and the quality of the test is very dependent to the
quality of the different estimations of $\widehat{\lambda}_j^*$.
As a consequence, we chose $5\leq m \leq 10$.
\end{enumerate}
\section{Numerical simulation and applications in Biomechanics}
\subsection{Simulations}
{\bf 1.} Initially, we apply the estimators and tests to several
simulated trajectories of classical F.B.M. (generated according to
a Choleski decomposition) with different values of
$H=0.2,~0.4,~0.6$ and $0.8$). The selected of values of different
parameters are~: $N=6000$, $\Delta_N=0.03$, $m=5$, $f_{min}=0.05$
and $f_{max}=20$. There are $30$ independent replications of each
time series. The results are presented in the following table~: \\
\begin{center}
\begin{tabular}{|l|c|c|c|c|}
\hline  Theoretical values of $H$ & 0.2 & 0.4 & 0.6 & 0.8 \\
\hline Empirical mean of $\widehat{H}$ & 0.148 & 0.384 & 0.599 &
0.821 \\ \hline Standard deviation of $\widehat{H}$ & 0.034 &
0.031 & 0.041 & 0.048 \\ \hline
\end{tabular}
\end{center}
~\\
The Figure 4 presents the log-log representation for one
trajectory~: the linearity is seeming. Moreover, Figure 5 exhibits
a histogram of the distribution of the test statistic $T_0^{(N)}$
(in this case $K=0$ and $30 \times 4=120$ independent realizations)
compared to a $\chi^2$-distribution with $3$ degrees of freedom.
The goodness-of-fit Kolmogorov-Smirnov test for $T_0^{(N)}$ to the
$\chi^2(3)$ distribution  is also accepted
(with $D\simeq 0.091$ and $p-value\simeq 0.272$).  \\ \\
\[
\epsfxsize 6 cm \epsfysize 4 cm \epsfbox{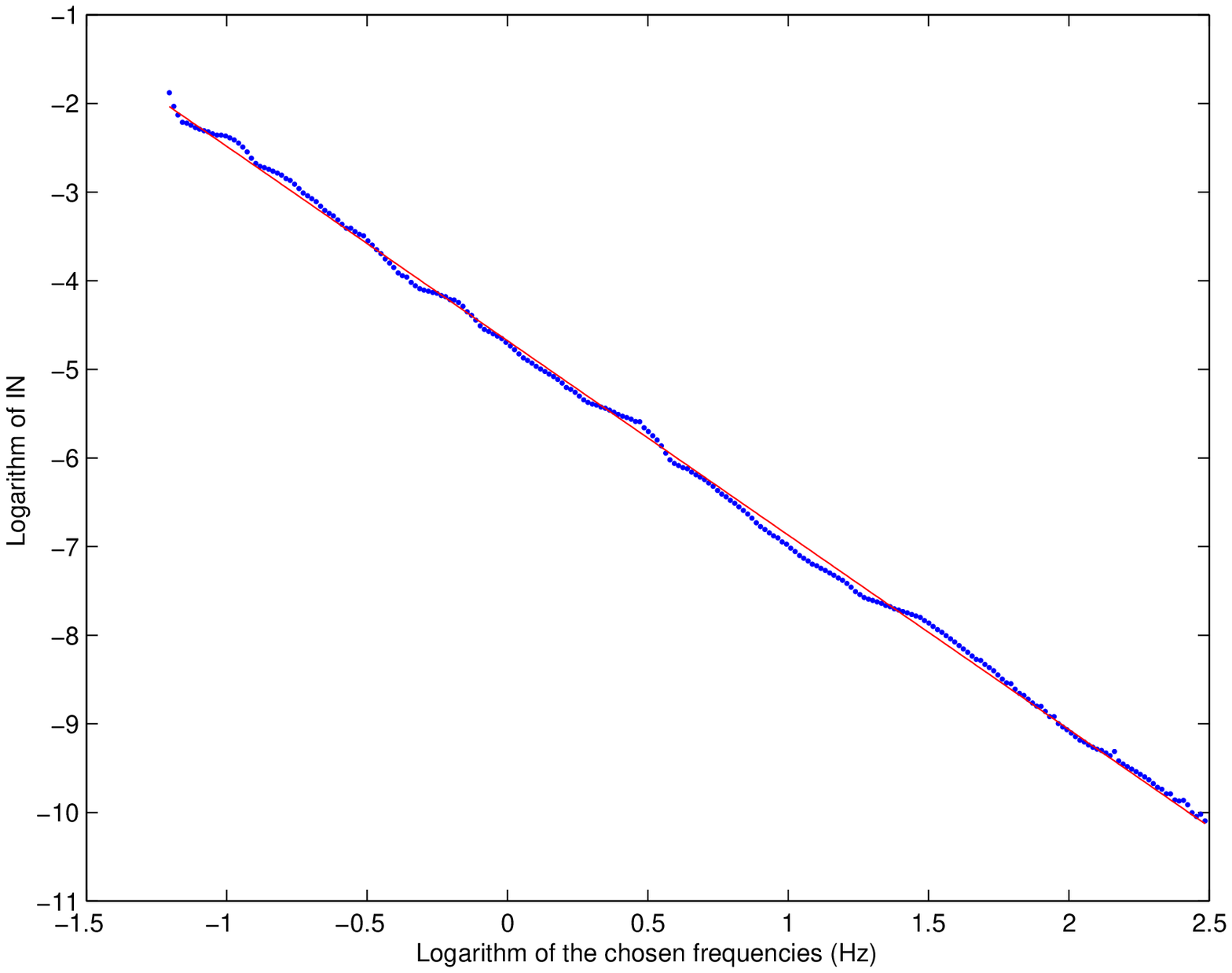}
\hspace{2cm} \epsfxsize 6 cm \epsfysize 4 cm
\epsfbox{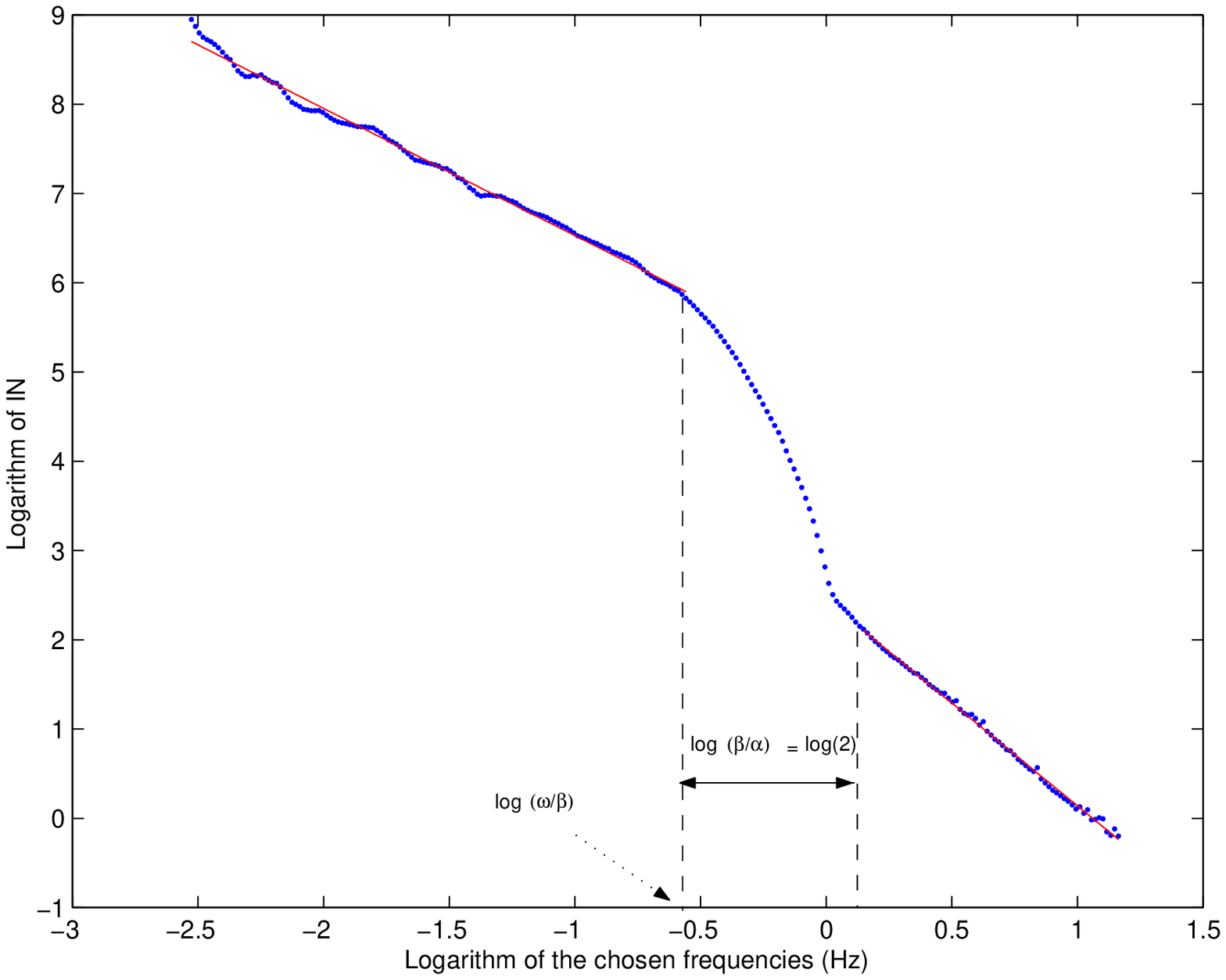}
\]
{\bf Figure 4 :} The log-log representation for a trajectory of a
$(M_0)$-FBM (left, with $H=0.6$) and $(M_1)$-FBM (right) \\
~\\
{\bf 2.} Then, we apply to $30$ independent replications
trajectories of $(M_1)$-FBM (generated according to a Choleski
decomposition with numerical approximations of the covariances)
with $H_0=0.2$ and $\sigma_0^2=10$, $H_1=0.7$, and $\sigma_1^2=5$,
and $\omega_1=5$. The results (with parameters~: $N=6000$,
$\Delta_N=0.03$, $m=5$, $f_{min}=0.8$ and $f_{max}=16$) are the
following~:
\\
\begin{center}
\begin{tabular}{|l|c|c|c|}
\hline Theoretical value & $H_0=0.2$ & $H_1=0.7$ & $\omega_1=5$  \\
\hline Empirical mean  &  0.197 &0.693 & 5.18  \\
\hline Standard deviation & 0.110 & 0.068 & 0.491
\\ \hline
\end{tabular}
\end{center}
~\\
Figure 4 presents the log-log representation for one trajectory,
with the $2$ regression lines. The hypothesis of the modelling
with a simple FBM (therefore with $K=0$) is always rejected (in
such a case, the model is misspecified and the statistic
$T_0^{(N)}$ is then between $39.3$ and $126.8$, very different
from the realizations of $\chi^2$-distribution with $3$ degrees of
freedom). On the contrary, the hypothesis of the modelling with a
$(M_1)$-FBM is always accepted and a histogram of the realizations
of the test statistic $T_1^{(N)}$ is presented in Figure 5
(compared to a $\chi^2$-distribution with $6$ degrees of freedom).
The goodness-of-fit Kolmogorov-Smirnov test for $T_1^{(N)}$ to the
$\chi^2(6)$ distribution is also accepted (with $D\simeq 0.187$
and $p-value\simeq 0.059$).
\[
\epsfxsize 6 cm \epsfysize 4 cm \epsfbox{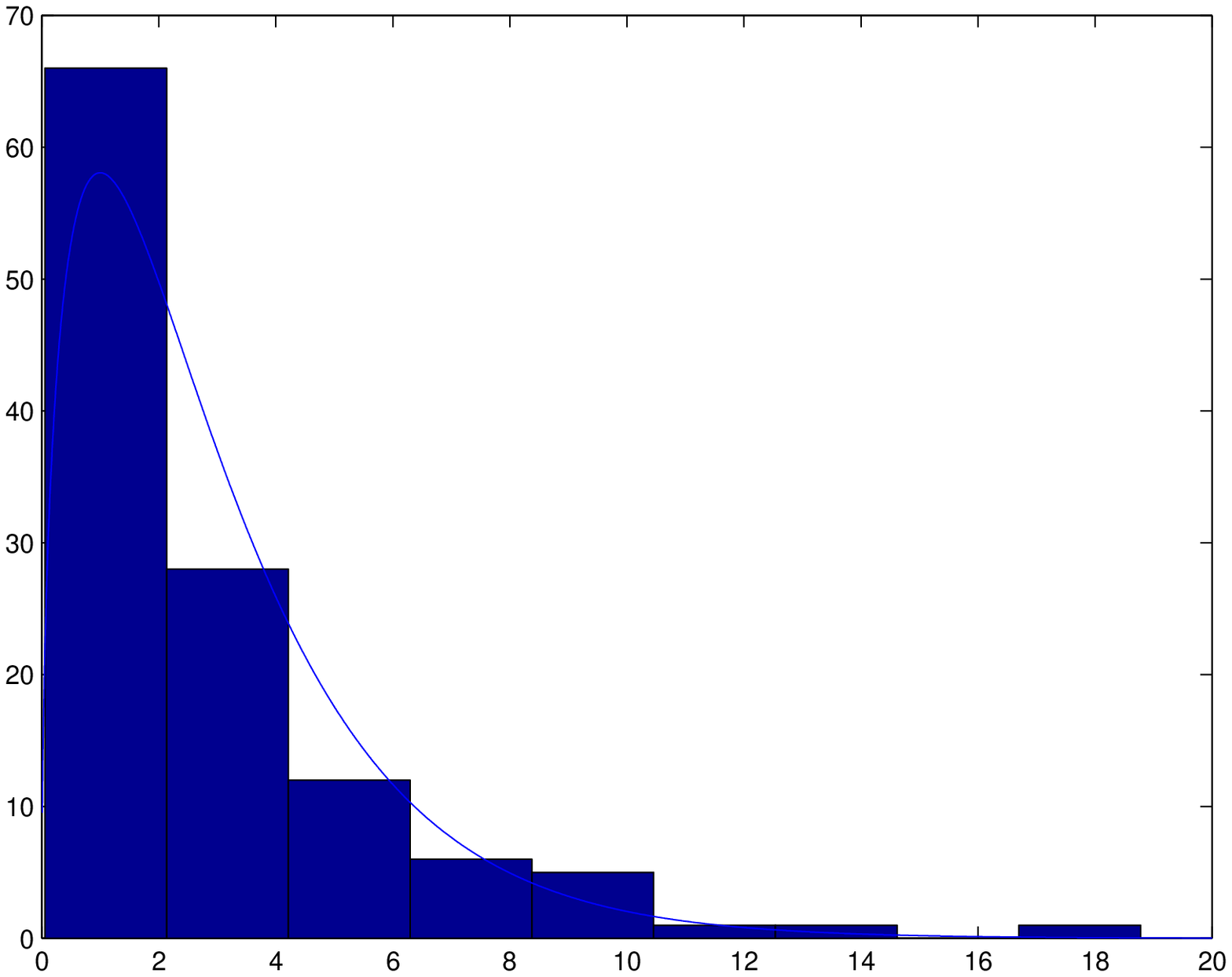} \hspace{2cm}
\epsfxsize 6 cm \epsfysize 4 cm \epsfbox{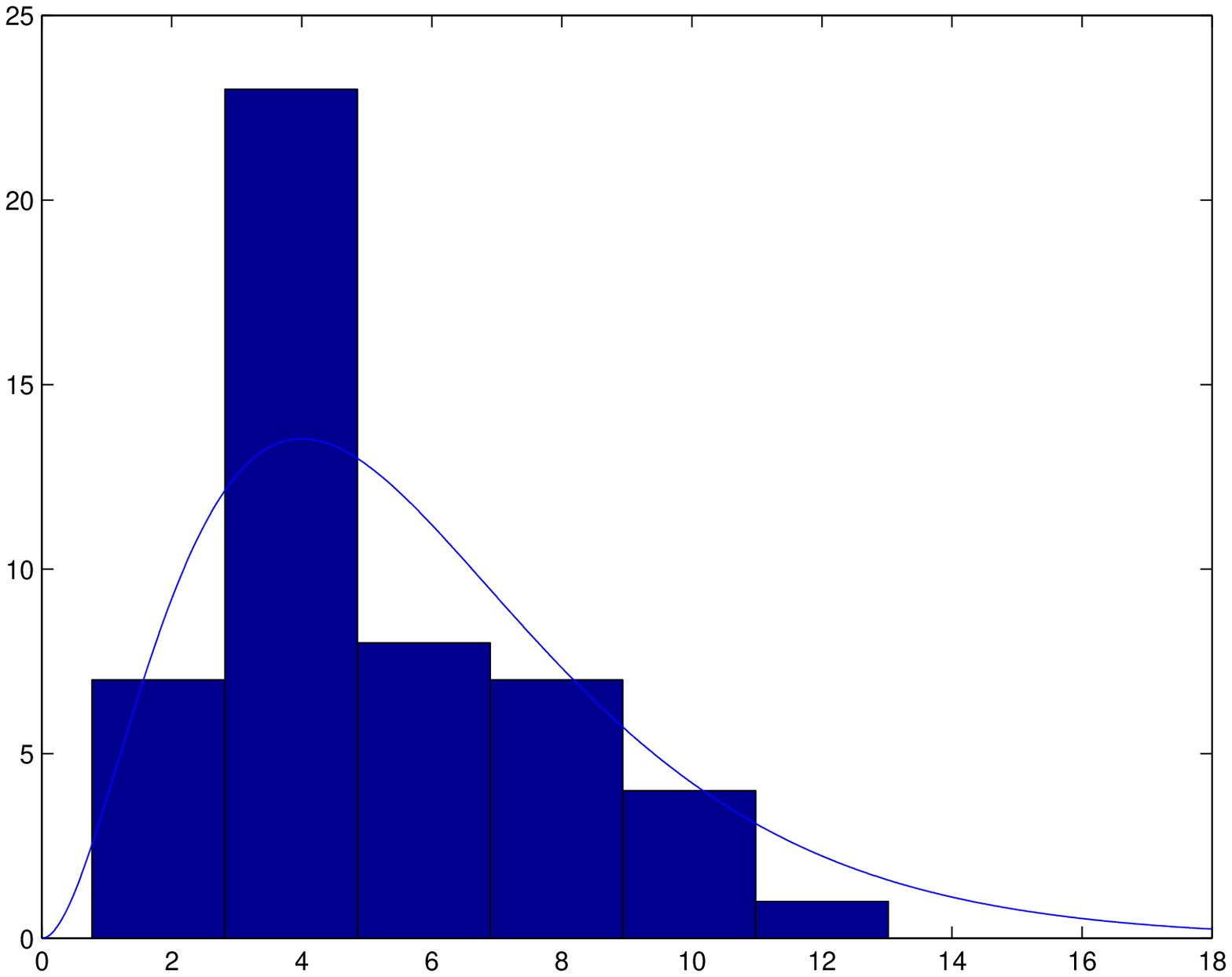}
\]
{\bf Figure 5 :} The empirical distribution of $T_0^{(N)}$ and
$T_1^{(N)}$ (respectively) compared to the corresponding $\chi^2$
distribution in the cases of simulated trajectories of $(M_0)$-FBM
(left) and $(M_1)$-FBM (right) \\
~\\
{\bf Conclusion of these simulations~:} the results are
surprisingly good compared with the complexity of the method. The
asymptotic distribution of the test statistics can be used for
real data. However, the computation time is important (especially
for the computation of the test statistic)~: $3$ hours are
necessary for the treatment of each $(M_1)$-FBM replication.
\subsection{Applications in Biomechanics}
We apply our statistics to different trajectories (see the
description in the Introduction) with the following parameters~:
\begin{itemize}
\item $N=6000$ and $\Delta_N=0.03$;
\item The mother wavelet is $\psi_1$ (with $\alpha=5$ and $\beta=10$).
\item The choice of the frequency band is $f_{min}=0.15$ and $f_{max}=15$
which corresponds to a detection frequency band $[0.52~,~38.32]$
Hz (with mother wavelet $\psi_1$);
\item $m=5$.
\end{itemize}
First, we study the $X$-trajectories of one subject (fore-aft
direction) for different feet position ($0;~2;~10;~20 \, cm$
clearance and $0;~15;~30;~45^o$ angle). In all the cases, the test
(with a type I error of $5 \%$) rejects the hypothesis of a
modelling with a simple $(M_0)$-FBM. But the modelling with a
$(M_1)$-FBM is accepted by the test $12$ times out of $16$, with
an empirical mean of $\widehat{\omega}_1\simeq 3.5$ and a standard
deviation  of $\widehat{\omega}_1\simeq 1$ (the different values
of $\tilde{H}_0$ and $\tilde{H}_1$ are in $[0.9,1]$ in the different cases). ~\\
~\\
For the different $Y$-trajectories (medio-lateral direction) of
the same patient, the test rejects the hypothesis of a modelling
with a simple $(M_0)$-FBM in all the case. The modelling with a
$(M_1)$-FBM is accepted by the test $13$ times out of $16$, with
an empirical mean of $\widehat{\omega}_1\simeq 3.1$ and a standard
deviation  of $\widehat{\omega}_1\simeq 1$ (the different values
of $\tilde{H}_0$ and $\tilde{H}_1$ are in $[0.8,1]$ in the different cases). ~\\
~\\
Figure 6 presents log-log plots of $J_N(f_k)$ versus $ f_k$ (i.e.
$\log J_N(f_k)$ vs. $\log f_k$) of all the experiments, for
$X$-trajectories (left)
and $Y$-trajectories (right).  ~\\
~\\
\[
\epsfxsize 6 cm \epsfysize 4 cm \epsfbox{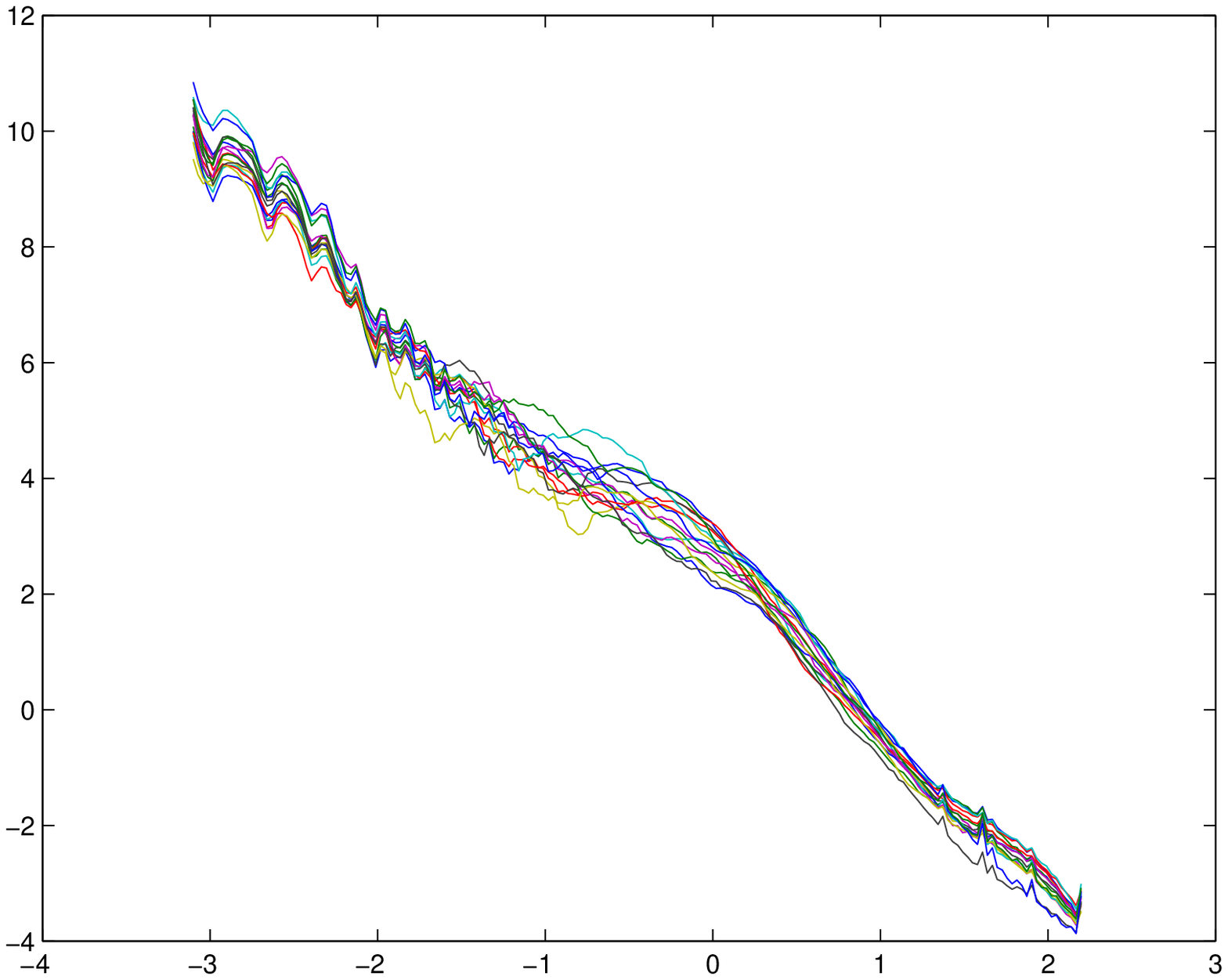} \hspace{2cm}
\epsfxsize 6 cm \epsfysize 4 cm \epsfbox{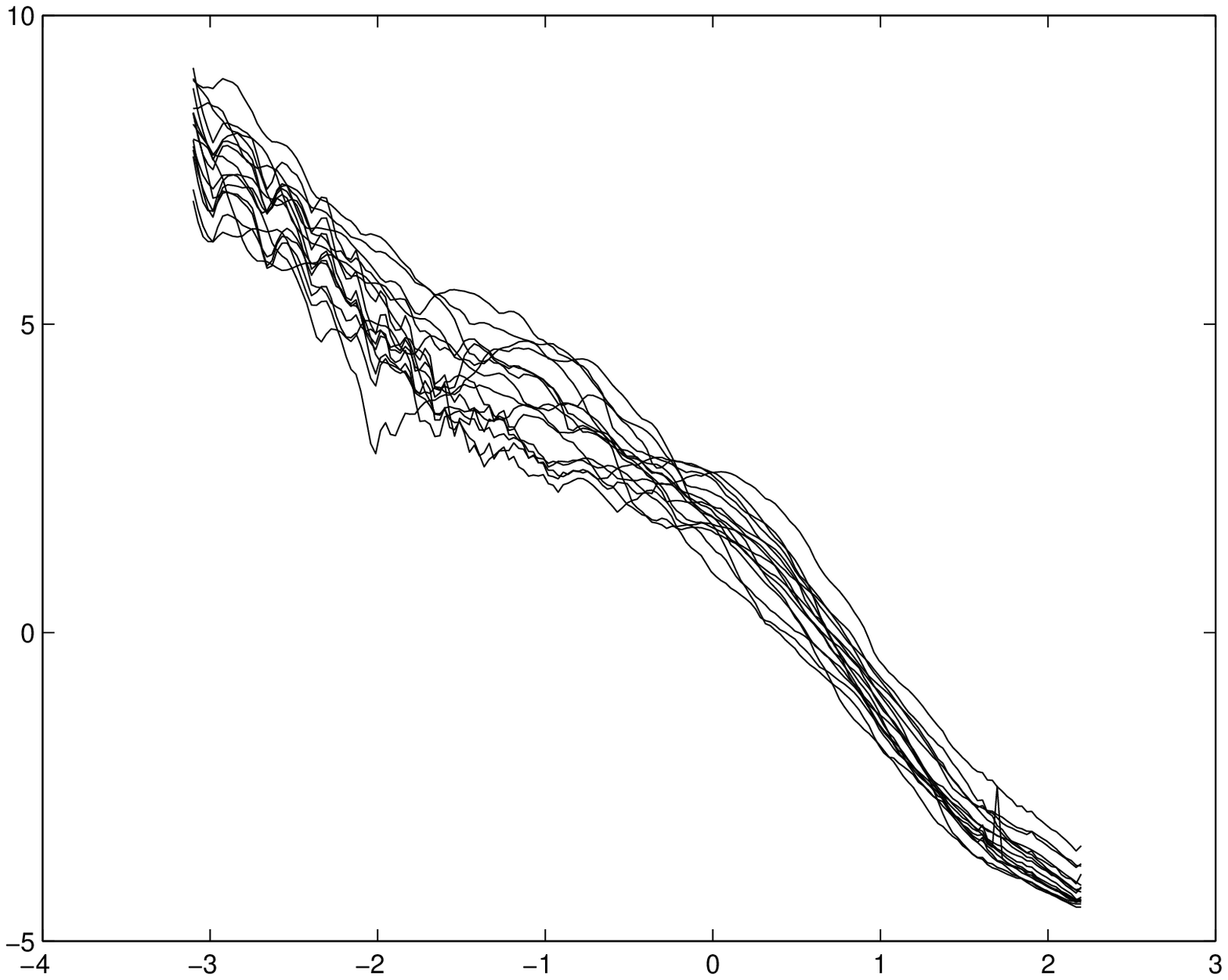}
\]
{\bf Figure 6 :} The log-log representation of the $16$ different
$X$-trajectories (left) and $Y$-trajectories (right) \\
~\\
{\bf Conclusion of these applications to biomechanics data~:} all
these results allow us  to give new interpretations on the upright
position. The behavior of $X$-trajectories and $Y$ trajectories
are very similar, for all the positions of the feet (the studied
statistics do not seem to depend on the angle and clearance of the
feet). The  $(M_1)$-FBM models of these trajectories fit well,
which suggests  two different type of behavior for low and high
frequencies. The frequency change is around $3$ Hz, which
corresponds to a physiological change~: this could be interpreted
for instance as the passage of a cerebral control of the stability
by the inner ear to a muscular auto-stabilization. We return to
\cite{B2DMV:IEEE} for a more detailed discussion of the
biomechanical interpretations. Such an estimation of this
frequency change would be very interesting for a better detection
of certain pathologies and to help in their cure.
\appendix
\section{Proofs}
\subsection{Proof of Theorem \ref{Theo:TCL:discret}}
First, we prove the following technical Lemma~:
\begin{lem}\label{maj:cov}
Let $X$ be a $(M_K)$-MBM. For $(t,t')\in \R^2$, and $(u,u')\in
\R_+^2$, define~: \ba \label{maj:cov:def} S(t,u,t',u')=\E \Big [
(X(t+u)-X(t))\cdot (X(t'+u')-X(t'))\Big]. \ea
\begin{enumerate}
\item For all $(u,u',t,t')\in \R_+^2\times \R^2$, there exists a constant $C>0$  depending only on the parameters $(\omega_j)_j$,
$(\sigma_j)_j$ and $(H_j)_j$ such that ~:
\ba \label{maj:cov:maj2}
\big|S(t,u,t',u')\big| \leq C \cdot \left ( u^{H_K} \cdot  \1 _{u\leq 1} +
u^{H_0} \cdot  \1_{u>1} \right )  \times \left ( u^{,H_K} \cdot
\1 _{u'\leq 1} + u^{,H_0} \cdot  \1_{u'>1} \right ); \ea
\item  More precisely, if $\big (
\max(u,u') \cdot \omega_K \big ) <1~~$ and $\displaystyle{~~\max(u,u')
\leq \frac 1 4 \cdot |t'-t|}~~$, there exists a constant $C>0$  depending only on
the parameters $(\omega_j)_j$, $(\sigma_j)_j$ and $(H_j)_j$ such
that ~:
\ba \label{maj:cov:maj}
\big|S(t,2u,t',2u')\big| \leq C \cdot \Big (u \cdot
u'+\max(u,u')^4 \Big ) \left ( \frac 1 {|t-t'+u'-u|} +
\max_{i=0,1,\cdots,K} \left \{\frac 1 {|t-t'+u'-u|^{2-2H_i}} \right \}
\right ). \ea
\end{enumerate}
\end{lem}
\begin{dem}
1/ First, the Cauchy-Schwarz inequality implies that
$$
\big|S(t,u,t',u')\big| \leq  \sqrt{\E \big [ (X(t+a)-X(t))^2\big ]} \times
\sqrt{\E \big [ (X(t'+a')-X(t')^2)\big ]}.
$$
But, $\displaystyle{~~\E \big [ (X(t+a)-X(t))\big ]^2 = 4
\sum_{j=0}^{K} \sigma_j^2 \cdot a^{2 H_j} \int_{a \omega_j}^{a
\omega_{j+1}} \frac {(1-\cos v)}{v^{2H_j +1}} \, dv.~~}$ Then, the
following expansions~:
$$
\int_{0}^{x} \frac {(1-\cos v)}{v^{2H +1}} \, dv = \left \{
\begin{array}{ll}
\displaystyle{\frac 1 {2 (2-2H)}\, x^{2-2H} + O(x^{4-2H})}
& \mbox{for}~~~x \to 0; \\
\displaystyle{ C(H)  - \frac 1 {2H}\, \frac 1 {x^{2H}} + O \big (
\frac 1 {x^{2H+1}} \big ) } & \mbox{for}~~~x \to \infty.
\end{array} \right .
$$
with $\displaystyle{C(H)=\int _0 ^\infty \frac {(1-\cos v)}{v^{2H
+1}} \, dv,~~~}$ imply that~:
\begin{eqnarray} \label{vitessevario}
\E \big [ (X(t+u)-X(t))\big ]^2 =  \left \{
\begin{array}{ll}
 4 \cdot  \sigma_K^2 \cdot  C(H_K) \cdot u^{2H_K} +O(u^2) &
\mbox{when}~~u
\to 0; \\
& \\
 4  \cdot \sigma_0^2 \cdot  C(H_0) \cdot u^{2H_0} +O(1) &
\mbox{when}~~u \to \infty;
\end{array} \right .
\end{eqnarray}
that achieves the proof of the majoration (\ref{maj:cov:maj2}). \\
~\\
 2/ We turn now to the proof of the upper bound
(\ref{maj:cov:maj}). To begin with, we remark that for  all
$(t,t',u,u')\in
\R^4$,
the  following equalities   are true~:
\begin{eqnarray}
\nonumber S(t,2u,t',2u')&=&\int _{\R} \frac {(e^{-i(t+2u)\xi}
-e^{-it\xi})(e^{i(t'+2u')\xi}
-e^{it'\xi})}{\rho^2(\xi)} \, d\xi  \\
\nonumber&=& \int _{\R}  \frac {(e^{-iu\xi} -e^{iu\xi})(e^{iu'\xi}
-e^{-iu'\xi})}{\rho^2(\xi)} \,e^{i\xi(t'-t)+i\xi(u'-u)} \,
d\xi \\
\nonumber&=& 8 \int _0 ^\infty \frac {\sin (u\xi) \cdot  \sin(u'\xi) \cdot
\cos \big ( \xi(t'-t+u'-u) \big )} {\rho^2(\xi)} \,
d\xi \\
\label{S} &=& 8 \sum _{i=0}^K \sigma_i^2   \int
_{\omega_i}^{\omega_{i+1}}\frac {\sin (u\xi) \cdot \sin(u'\xi)
\cdot \cos
\big ( \xi(t'-t+u'-u) \big )} {\xi^{2H_i+1}} \, d\xi.
\end{eqnarray}
Then, we bound the different integrals. \\$\bullet$
\mbox{First,} we threat the case
$i=K$ that is when the upper limit of the integral is $\infty$. In
this case, we can rewrite the integral between $\omega_K$ and
${\infty}$ as the difference of the  integral between $0$
 and ${\infty}$ and the one between $0$ and $\omega_K$, that is
 $\displaystyle{\int _{\omega_K}^{\infty}\gamma(\xi)\,d\xi = \int_{0}^{\infty}\gamma(\xi)\,d\xi\, - \int _0^{\omega_K}\gamma(\xi)\,d\xi}$
where $\gamma(\xi)=\sin (u\xi) \cdot \sin(u'\xi) \cdot \cos
\big ( \xi(t'-t+a'-a) \big )\,\xi^{-(2H_K+1)}$.
%%\\ ({\bf (2eme version au choix)} we have
%%\ban \int _{\omega_K}^{\infty} \frac {\sin (u\xi)\cdot \sin(u'\xi) \cdot \cos
%%\big ( \xi(t'-t+a'-a) \big )} {\xi^{2H_i+1}} \, d\xi  &=& \int_{0}^{\infty} \frac {\sin (u\xi) \cdot \sin(u'\xi) \cdot \cos
%%\big ( \xi(t'-t+a'-a) \big )} {\xi^{2H_i+1}} \, d\xi
%%\\ &&- \int _0^{\omega_K}\frac {\sin (u\xi) \cdot \sin(u'\xi) \cdot \cos
%%\big ( \xi(t'-t+a'-a) \big )} {\xi^{2H_i+1}} \, d\xi\ean
The second integral of the right hand side can be bounded by the
same argument than the terms of order $i=0$ in (\ref{S}). The
first one corresponds to the expression of the covariance of the
increments of a F.B.M. $B_{H_K}$ with Hurst parameter $H_K$ %%\in(0,1)$
and variance $1$. Thus, for all $(t,t')\in
\R^2$, $(u,u')\in
\R_+^2$ such that $4 \cdot \max (u,u') \leq |t'-t|$, we get~:
\begin{eqnarray}
\nonumber&&  \hspace{-1.5cm}
\left | 8  \int _0 ^\infty \frac {\sin (u\xi) \cdot
\sin(u'\xi) \cdot \cos \big (
\xi(t'-t+u'-u) \big )} {\xi^{2H_K+1}} \, d\xi \right |
\\
\nonumber&= &
\left |\E \Big ((B_{H_K}(t+2u)-B_{H_K}(t))\cdot (B_{H_K}(t'+2u')-B_{H_K}(t')) \Big )  \right |
  \\
\nonumber & =& \frac 1 {2C^2(H_K)}\left | \Big (|t-t'+2u|^{2H_K} -
|t-t'+2u-2u'|^{2H_K}-|t-t'|^{2H_K}+|t'-t-2u'|^{2H_K} \Big ) \right |\\
\label{inegal1} & \leq & D(H_K)  \cdot \frac {u \cdot u'}{|t-t'+u'-u|^{2-2H_K}},
~~~~\mbox{with}~~D(H_K)>0.
\end{eqnarray}
$\bullet$
\mbox{Next,} we consider the integrals with a finite upper limit and a
non-zero lower limit. This corresponds to  $i=1,\dots,K-1$. In
these cases, for $b>0$, an integration by parts provides us
\begin{multline}
\label{IPP}
\int_{\omega_i}^{\omega_{i+1}}\frac {\sin (u\xi) \cdot \sin(u'\xi) \cdot
\cos (b \xi)} {\xi^{2H_i+1}} \, d\xi = \frac 1 b  \left ( -
\int_{\omega_i}^{\omega_{i+1}} \frac {u \cdot \cos (u\xi) \cdot
\sin(u'\xi) + u' \cdot \cos (u'\xi) \cdot \sin(u\xi)} {\xi^{2H_i+1}}\cdot
\sin (b \xi) \, d\xi
  \right .\\
\left .  + (2H_i+1)\,\int_{\omega_i}^{\omega_{i+1}} \frac {
\sin(u\xi) \cdot \sin(u'\xi)} {\xi^{2H_i+2}}\cdot \sin (b \xi) \, d\xi +
 \left
[\frac {\sin (u\xi) \cdot \sin(u'\xi) \cdot \sin (b \xi)}
{\xi^{2H_i+1}}
\right ]_{\omega_i}^{\omega_{i+1}} \right ).
\end{multline}
By using  the   majoration $|\sin(ux)| \leq ux$, $|\sin(u'x)|
\leq u'x$, $|\cos(ux)|\leq 1$, $|\cos(u'x)|\leq 1$ and $|\sin(bx)|
\leq 1$ for $x \geq 0$, we deduce that for all $(u,u',b) \in
\R_+^3$,
\begin{eqnarray}\label{casii}
\Big | \int_{\omega_i}^{\omega_{i+1}}\frac {\sin (u\xi) \cdot
\sin(u'\xi) \cdot \cos (b \xi)} {\xi^{2H_i+1}} \, d\xi \Big | \leq C_i
\cdot \frac {u \cdot u'} b,
\end{eqnarray}
where $C_i>0$ is a constant depending only on $H_i, \omega_i$ and
$\omega_{i+1}$.
\\
$\bullet$ \mbox{Finally, } it remains to bound the two integrals
with lower limit $0$. We will show only how to bound \\
$\displaystyle{\int_0^{\omega_1}\frac {\sin (u\xi) \cdot \sin(u'\xi) \cdot
\cos (b \xi)} {\xi^{2H_i+1}} \, d\xi}$, since the other integral
can be treated similarly. The integration by part formula
(\ref{IPP}) remains valid even  when the lower limit is $0$.
Indeed, the integrand can be bounded by $\displaystyle{C \times
\xi^{1-2H_i}}$ and $\displaystyle{\int_0^1 \xi^{1-2H_i}\,d\xi
<\infty}$ as soon as $H_i<1$.
 After this remark, we bound the three terms of
the right hand side of (\ref{IPP}).
\\
i) From   $|\sin(ux)| \leq ux$, $|\sin(u'x)|
\leq u'x$  and $|\sin(bx)|
\leq 1$ for $x \geq 0$, we deduce that for all $(u,u',b) \in
\R_+^3$,
\ba
\label{maj:crochet}
\Big |\left [\frac {\sin (u\xi)
\cdot
\sin(u'\xi) \cdot \sin (b \xi)} {\xi^{2H_0+1}}
\right ]_{0}^{\omega_1} \Big | &\leq& \big
(\omega_{1}^{1-2H_0} \big ) \cdot u \cdot u'.
\ea
 ii) For all $(\xi,\xi') \in
[0,\omega_{1}]$, the  power series expansion of $x \mapsto
\sin(x)$ implies that
$$
\sin(u\xi) \cdot \sin(u\xi) = u \cdot u'\cdot  \xi ^2 +
\sum_{k=1}^\infty \Big ( \sum_{j=0}^k \frac {u^{2j+1} \cdot
(u')^{2(k-j)+1}}{(2j+1)! \cdot (2(k-j)+1)!} \Big ) \cdot (-1)^k
\xi^{2k+2}.
$$ One can remark that
\begin{eqnarray*}
\sum_{j=0}^k \frac {u^{2j+1} \cdot (u')^{2(k-j)+1}}{(2j+1)! \cdot
(2(k-j)+1)!} &\leq & \max(u,u')^{2k+2}
\sum_{j=0}^k \frac {1}{(2j+1)! \cdot (2(k-j)+1)!}  \\
&\leq &   \max(u,u')^{2k+2},
\end{eqnarray*}
because $\sum_{j\geq 0} \frac {1}{(2j+1)!} \leq 2$. As a
consequence, when $(\max(u,u') \cdot \omega_{1}) <1$ and $b>0$,
integration and summation can be interchanged and
\begin{eqnarray*}
\left | \int_{0}^{\omega_{1}}\frac {\sin (u\xi) \cdot \sin(u'\xi) \cdot
\sin (b \xi)} {\xi^{2H_0+2}} \, d\xi -u \cdot u'
\int_{0}^{\omega_{1}}\frac { \sin ( b\xi)} {\xi^{2H_0}}
\, d\xi \right | &\leq &  \sum_{k=1}^\infty  \left (
\max(u,u')^{2k+2}  \cdot
\int_{0}^{\omega_{1}}  \xi ^{2k-2H_0} \, d\xi\right )
\\
&   \leq &\frac { \max(u,u')^4 \cdot
\omega_{1}^{3-2H_0}} { 1 -(\max(u,u')\cdot \omega_{1})^2 }.
\end{eqnarray*}
%%%%%%%%%%%%%%%%%%%%%%%%%%%%%%%
%% MISE EN PAGE DE JMB
%%\begin{eqnarray*}
%%\left | \int_{0}^{\omega_{1}}\frac {\sin (u\xi) \cdot \sin(u'\xi) \cdot
%%\sin (b \xi)} {\xi^{2H_0+2}} \, d\xi -u \cdot u'
%%\int_{0}^{\omega_{1}}\frac { \sin ( b\xi)} {\xi^{2H_0}}
%%\, d\xi \right | &\leq &  \\
%%&& \hspace{-5cm} \leq  \sum_{k=1}^\infty  \left (
%%\max(u,u')^{2k+2}  \cdot
%%\int_{0}^{\omega_{1}}  \xi ^{2k-2H_0} \, d\xi\right ) \\
%%&& \hspace{-5cm} \leq \frac { \max(u,u')^4 \cdot
%%\omega_{1}^{3-2H_0}} { 1 -(\max(u,u')\cdot \omega_{1})^2 }.
%%\end{eqnarray*}
But $\displaystyle{\int_{0}^{\omega_{1}}\frac { \sin ( b\xi)}
{\xi^{2H_0}} \, d\xi =b^{2H_0-1}\int_{0}^{b \cdot \omega_{1}}\frac
{ \sin (\xi)} {\xi^{2H_0}} \, d\xi}$. Denote $\displaystyle{M(H)=
\sup_{x \in \R_+} \Big |\int _0 ^x \frac { \sin (\xi)} {\xi^{2H}}
\, d\xi \Big |}$ for $0<H<1$. Thus
\begin{eqnarray}\label{sin0}
\Big | \int_{0}^{\omega_{1}}\frac {\sin (u\xi) \cdot \sin(u'\xi) \cdot \sin
(b \xi)} {\xi^{2H_0+2}} \, d\xi \Big | \leq M(H_0)  \cdot u \cdot
u'\cdot b^{2H_0-1}+ \frac {\omega_{1}^{3-2H_0}\cdot \max(u,u')^4}
{ 1 -(\max(u,u')\cdot \omega_{1})^2}.
\end{eqnarray}
iii)  Similarly for  $(\xi,\xi')
\in [0,\omega_{1}]$, we have
$$
\cos(u\xi) \cdot \sin(u'\xi) = u'\cdot \xi + \sum_{k=1}^\infty \Big (
\sum_{j=0}^k \frac {u^{2j} \cdot (u')^{2(k-j)+1}}{(2j)! \cdot (2(k-j)+1)!}
\Big ) \cdot (-1)^k \xi^{2k+1},
$$
But for $k \geq 1$, $\displaystyle{\Big ( \sum_{j=0}^k \frac
{u^{2j} \cdot (u')^{2(k-j)+1}}{(2j)! \cdot (2(k-j)+1)!} \Big ) \leq
\max(u,u')^{2k+1}}$. As a consequence, when $(\max(u,u') \cdot
\omega_{1}) <1$ and $b>0$, integration and summation can be
interchanged and we get
\begin{eqnarray}\label{cos0}
\left | \int_{0}^{\omega_{1}}\frac {u \cdot \cos (u\xi) \cdot \sin(u'\xi) \cdot
\sin (b \xi)} {\xi^{2H_0+1}} \, d\xi -u \cdot u'
\int_{0}^{\omega_{1}}\frac { \sin ( b\xi)} {\xi^{2H_0}} \, d\xi
\right |  \leq \frac { \max(u,u')^4 \cdot \omega_{1}^{3-2H_0}} { 1
-(\max(u,u')\cdot \omega_{1})^2 }.
\end{eqnarray}
Therefore from (\ref{maj:crochet}), (\ref{sin0}), (\ref{cos0}) and
(\ref{inegal1}), we deduce~for $(u,u')$ such that $\max(u,u')\cdot
\omega_{K}<1/2$~:
%\begin{eqnarray}
%\Big | \int_{0}^{\omega_{1}}\frac {a \cdot \cos (u\xi) \cdot
%\sin(u'\xi) \cdot \sin (b \xi)} {\xi^{2H_0+1}} \, d\xi \Big | \leq
%M(H_0)  \cdot u \cdot u'\cdot b^{2H_0-1}+ \frac {\omega_{1}^{3-2H_0}\cdot
%\max(u,u')^4} { 1 -(\max(u,u')\cdot \omega_{1})^2}.
%\end{eqnarray}
%Finally,
%$$
%\Big | \int_{0}^{\omega_{1}}\frac {\sin (u\xi) \cdot \sin(u'\xi) \cdot \cos
%(b \xi)} {\xi^{2H_0+1}} \, d\xi \Big | \leq  3 \cdot M(H_0)  \cdot
%\frac{ u \cdot u'}{b^{2-2H_0}}+ 4 \cdot
%\frac{\max(u,u')^4}{b}\cdot \omega_{1}^{3-2H_0} + \frac {u \cdot u'}b \cdot \omega_{1}^{1-2H_0},
%$$
%and by using (\ref{inegal1}), we get
$$
\Big | \int_{\omega_K}^{\infty}\frac {\sin (u\xi) \cdot \sin(u'\xi) \cdot
\cos (b \xi)} {\xi^{2H_K+1}} \, d\xi \Big | \leq  (D(H_K)+3
M(H_K)) \cdot \frac{ u \cdot u'}{b^{2-2H_K}}+  4 \cdot
\frac{\max(u,u')^4}{b}\cdot \omega_{K}^{3-2H_K}+ \frac {u \cdot u'}b \cdot \omega_{K}^{1-2H_K}.
$$
By combining the two previous bounds with (\ref{S}) and
(\ref{casii}), we deduce (\ref{maj:cov:maj}) and this finishes the
proof.
\end{dem}
~\\
The proof of Theorem \ref{Theo:TCL:discret} uses  the two
following lemmas:
\begin{lem}\label{maj:sum:epsilon} Under the same notations and assumptions
as in  Theorem \ref{Theo:TCL:discret}, there exists two constants
$C_1>0$ and $C_2>0$ depending only on $r$, $a_{min}$ and $a_{max}$
such that for all $N$
\ba
\label{maj:unif:Eeps2}
i)\hspace{1.5cm} \sup_{a\in[a_{min},\,a_{max}]}\max_{k\in
D_N(a)}\E\varepsilon^2_N(a,k ) &\le& C_1\times
\varphi(N)\hspace{7cm}
\ea
\ba
\label{maj:unif:diff:Eeps2}
ii)\hspace{0.3cm}\sup_{a_1,\,a_2\in[a_{min},\,a_{max}]}\,\max_{k\in
D_N(a_1)\cap D_N(a_2)} \E\left|\sqrt{a_1}\,\varepsilon_N(a_1,k
)-\sqrt{a_2}\,\varepsilon_N(a_2,k )\right|^2 &\le& C_2\times
\varphi(N)\times \left|a_2-a_1\right|^2\hspace{1cm}
\ea
where $\displaystyle{\varphi(N) =  \Delta_N^2+ \Delta_N^{1+2H_K} +
\Delta_N^2
\log N +
\Delta_N^{1 + 2
\overline{H}} N ^{-1+2\overline{H}}
 +
(N\Delta_N)^{-2} }$ with $\overline{H}=\max
\{H_i~, ~i=0,\cdots,K\}$. \\
iii) Moreover $\displaystyle{(N \, \Delta_N )\, \varphi(N) \to 0}$
when $N\to\infty$.
\end{lem}
\begin{dem}
The error $ \varepsilon_N(a,k)$ contains three different terms,
the first  one corresponds to the replacement of the integral onto
the interval $[0, T_N ]$ by its Riemann sum, the second and the
third ones correspond to the replacement of the integral onto $\R$
by the integral onto the interval $ [0, T_N ]$ where $T_N = N
\Delta_N$. More precisely,  we have
\ba
\label{dec:epsi}
 \varepsilon_N(a,k) &= &
\frac 1{\sqrt a} \times
\left(\varepsilon_{1,N}(a,k)+
\varepsilon_{2,N}(a,k) + \varepsilon_{3,N}(a,k)\right)\ea with
\ban \varepsilon_{1,N}(a,k) &=&  \int_{ 0}^{T_N }
\psi(\frac{t}{a}-k\Delta_N)\,X(t)\,dt \,-\, \Delta_N\, \sum
_{p=0}^{N-1}\psi ( \frac {p \Delta_N} a -k\Delta_N)\,X
(p\Delta_N),
\\
\varepsilon_{2,N}(a,k) &=&\int_{T_N  }^{\infty}
\psi(\frac{t}{a}-k\Delta_N)\,X(t)\,dt,
\\
\varepsilon_{3,N}(a,k) &=& \int_{-\infty }^{ 0}
\psi(\frac{t}{a}-k\Delta_N)\,X(t)\,dt. \ean By using $(x+y+z)^2
\le 3\,(x^2+y^2+z^2)$ for all real numbers $x$, $y$, $z$, we
deduce \ba\label{maj:E:epsilon}   \E \varepsilon_N^2(a,k) &\le&
\left( \frac{3}{a}\right)\times \left(\sum_{i=1}^3\E
\varepsilon_{i,N}^2(a,k) \right).\ea
We now bound the different
terms $\E
\varepsilon_{i,N}^2(a,k)$ for $i=1,2, 3$~: \\
~\\
{\bf (1) Bound of $\displaystyle{\E
\varepsilon_{1,N}^2(a,k)}$.} \\
~\\
We have the decomposition
$\displaystyle{\,\varepsilon_{1,N}(a,k)\, =
\,I_{1,N}(a,k)\,+\,I_{2,N}(a,k)\,}$, where
\begin{eqnarray*}
I_{1,N}(a,k) &= &\sum_{p=0}^{N-1}
\int_{p\Delta_N}^{(p+1)\Delta_N } \psi \Big (\frac{t}{a}-k\Delta_N
\Big  )\, \Big (X(t)-X(p\Delta_N) \Big )\,dt \\
\mathrm{and}\qquad\qquad I_{2,N}(a,k)&=&   \sum_{p=0}^{N-1} \int_{p\Delta_N}^{(p+1)\Delta_N
}
\left (
\psi \Big (\frac{t}{a}-k\Delta_N \Big  ) - \psi  \Big ( \frac {p
\Delta_N} a -k\Delta_N \Big ) \right )\,X (p\Delta_N)\,dt.
\end{eqnarray*}
Then, the  inequality $(x+y)^2 \leq 2(x^2+y^2)$ which is valid for
all $(x,y)\in
\R^2$, implies
$$
\E \varepsilon_{1,N}^2(a,k) \leq
2\, \E \big ( I_{1,N}^2(a,k)\big ) + 2\, \E \big
(I_{2,N}^2(a,k)\big ).
$$
On one hand, we have
\begin{eqnarray*}
\E \big ( I_{1,N}^2(a,k)\big )&\hspace{-0.3cm} =\hspace{-0.3cm}&
\sum_{p=0}^{N-1} \sum_{p'=0}^{N-1}\int_{p\Delta_N}^{(p+1)\Delta_N
}\hspace{-0.3cm} \int_{p'\Delta_N}^{(p'+1)\Delta_N
}\hspace{-1cm}dtdt' \psi \Big (\frac{t}{a}-k\Delta_N \Big  ) \psi
\Big (\frac{t'}{a}-k\Delta_N \Big  ) \E
\Big ((X(t)-X(p\Delta_N))(X(t')-X(p'\Delta_N)) \Big ) \\
&\hspace{-0.3cm} =\hspace{-0.3cm}& \sum_{p=0}^{N-1}
\sum_{p'=0}^{N-1} \int_{0}^{\Delta_N }\hspace{-0.3cm}
\int_{0}^{\Delta_N }\hspace{-0.5cm}dudu'\psi \Big
(\frac{u+p\Delta_N}{a}-k\Delta_N \Big  ) \psi \Big
(\frac{u'+p'\Delta_N}{a}-k\Delta_N \Big  )
S(p\Delta_N,u,p'\Delta_N,u').
\end{eqnarray*}
Afterwards, we use Lemma \ref{maj:cov} to bound the terms
$S(p\Delta_N,u,p'\Delta_N,u')$, where   we use  different types of
bounds depending wether   $(p,p')$ is in the vicinity of the
diagonal or not. Namely, when $|p-p'| \leq 3$ we use the upper
bound (\ref{maj:cov:maj2}), and otherwise we  use the one in
(\ref{maj:cov:maj}). Observe that the assumptions of Lemma
\ref{maj:cov} are satisfied for large enough $N$ since $N\Delta_N
\to \infty$, as $N\to \infty$. Thus, when $N$ is large enough, we
get
\begin{eqnarray*}
\E \big ( I_{1,N}^2(a,k)\big )
 &\hspace{-0.3cm}  \leq \hspace{-0.3cm}&\hspace{-0.3cm} \sum_{\scriptsize
\begin{array}{c} p=0,p'=0 \\\ |p-p'| \leq 3  \end{array}
}^{N-1}   \int_{0}^{\Delta_N }\hspace{-0.3cm} \int_{0}^{\Delta_N
}\hspace{-0.3cm}du du'\left |\psi \Big
(\frac{u+p\Delta_N}{a}-k\Delta_N \Big  )\right | \left |\psi \Big
(\frac{u'+p'\Delta_N}{a}-k\Delta_N \Big  )\right |C \cdot (u \cdot
u')^{H_K}+  \\
 && \hspace{0cm} +
\sum_{\scriptsize
\begin{array}{c} p=0,p'=0 \\ |p-p'| \geq 4  \end{array}
}^{N-1}\hspace{-0.5cm} \int_{0}^{\Delta_N }\hspace{-0.3cm}
\int_{0}^{\Delta_N }\hspace{-0.5cm}du du' \left | \psi \Big
(\frac{u+p\Delta_N}{a}-k\Delta_N \Big  ) \psi \Big
(\frac{u'+p'\Delta_N}{a}-k\Delta_N \Big  )  \right | \times \cdots
\\
&& \hspace{3cm} \cdots \times C \cdot \Delta_N^2 \left ( \frac {
1} {(|p'-p|-1)\Delta_N} + \max_{i=0,\cdots,K} \left
\{
\frac { 1} {((|p'-p|-1)\Delta_N)^{2-2H_i} } \right
\} \right )
\\
&\le&C\cdot \Delta_N^{2+2H_K}\hspace{-0.3cm} \sum_{\scriptsize
\begin{array}{c} p=0,p'=0 \\\ |p-p'| \leq 3  \end{array}
}^{N-1} \sup_{\theta\in(0,1)}\left |\psi \left
(\Big(\frac{\theta+p}{a}-k \Big)\Delta_N
\right)\right | \times \sup_{\theta'\in(0,1)}\left |\psi \left
(\Big(\frac{\theta'+p'}{a}-k \Big)\Delta_N
\right)\right |
\\ && \hspace{0cm} + C\cdot \Delta_N^{4}
\sum_{\scriptsize
\begin{array}{c} p=0,p'=0 \\ |p-p'| \geq 4  \end{array}
}^{N-1}   \sup_{\theta\in(0,1)}\left |\psi \left
(\Big(\frac{\theta+p}{a}-k \Big)\Delta_N
\right)\right | \times \sup_{\theta'\in(0,1)}\left |\psi \left
(\Big(\frac{\theta'+p'}{a}-k \Big)\Delta_N
\right)\right | \times \cdots
\\&& \hspace{3cm} \cdots\times \left ( \frac { 1} {(|p'-p|-1)\Delta_N} +
\max_{i=0,\cdots,K} \left
\{
\frac { 1} {((|p'-p|-1)\Delta_N)^{2-2H_i} } \right
\} \right ).
\end{eqnarray*}
However, according to  Assumption (A1), for every integer $m \in
\N^*$ there exists a constant   $C>0$
such that  for all $x
\in \R$, $|\psi(x)| \leq C \cdot (1+|x|)^{-m}$. In particular, $\displaystyle{\sup_{\theta'\in(0,1)}\left |\psi \left
(\Big(\frac{\theta'+p}{a}-k \Big)\Delta_N
\right)\right |}$ is bounded. Therefore,
\begin{eqnarray*}
\E \big ( I_{1,N}^2(a,k)\big )
&\le&C\cdot \Delta_N^{2+2H_K}  \sum_{ p=0}^{N-1}
\sup_{\theta\in(0,1)}\Big(1+
\left|\theta+p-ak \right|\Delta_N/a
\Big)^{-m}
\\ &+&   C\cdot \Delta_N^{4}
\sum_{ p=0}^{N-1}   \sup_{\theta\in(0,1)}\Big(1+
\left|\theta+p-ak \right|\Delta_N/a
\Big)^{-m}
\,
\sum_{ q=3}^{N} \left ( \frac { 1} {q\Delta_N} +
\max_{i=0,\cdots,K} \left
\{
\frac { 1} {(q\Delta_N)^{2-2H_i} } \right
\} \right ).
\end{eqnarray*}
But
\ban
\Delta_N \sum_{p=0}^{N-1}  \sup_{\theta\in(0,1)}\left(1+
\left|\theta+p-ak \right|\frac{\Delta_N}{a}\right)^{-m}&\le&C\,\Delta_N \sum_{\ell=-\infty}^{\infty}
\left(1+
\frac{\left|\ell \Delta_N\right|}{a}\right)^{-m}\,\le\, C\,  \int_{-\infty}^\infty \left(1+ \frac{|x|}{a}\right)^{-m} \, dx
\\ &=& C\,|a|\,   \int_{-\infty}^\infty (1+ |y|)^{-m} \, dx\quad\le\,C\,|a_{max}|.
\ean
 Let us  denote $\overline{H}=\max\{H_i~, ~i=0,\cdots,K\}$ and
$\underline{H}=\min \{H_i~, ~i=0,\cdots,K\}$. We deduce
\begin{eqnarray}
\nonumber
\E \big ( I_{1,N}^2(a,k)\big )
&\le&C\cdot \Delta_N^{1+2H_K}  +   C\cdot \Delta_N^{3}
\,
\sum_{ q=3}^{N} \left ( \frac { 1} {q\Delta_N} +
\max_{i=0,\cdots,K} \left
\{
\frac { 1} {(q\Delta_N)^{2-2H_i} } \right
\} \right )
\\
\nonumber
&\le&
 C
 \left (\Delta_N^{1+2H_K} + \Delta_N^2\,\log N+
\Delta_N^3  \Big (  \frac 1 {\Delta_N
^{2-2\underline{H}}} \sum_{q=3}^{\Delta_N^{-1}} \frac 1 {q
^{2-2\underline{H}}} + \frac 1 {\Delta_N ^{2-2\overline{H}}}
\sum_{q=\Delta_N^{-1}}^N  \frac 1 {q ^{2-2\overline{H}}} \Big )
 \right )
\\
\label{maj:I1} &\hspace{-0.3cm}  \leq \hspace{-0.3cm}& C
\left (
\Delta_N^{1+2H_K} + \Delta_N^2\,\log N\,+
\Delta_N^2
\,+\,\Delta_N^{1+2\overline{H}} N^{-1+2\overline{H}}\right).
\end{eqnarray}
On the other hand, by using  Lemma \ref{maj:cov}, formula
(\ref{maj:cov:def}), we get
\begin{eqnarray*}
\E \big ( I_{2,N}^2(a,k)\big )&\hspace{-0.3cm} =\hspace{-0.3cm}&
\sum_{p=0}^{N-1} \sum_{p'=0}^{N-1}\int_{p\Delta_N}^{(p+1)\Delta_N
}\hspace{-0.3cm} \int_{p'\Delta_N}^{(p'+1)\Delta_N
}\hspace{-1cm}dudu' \left ( \psi \Big (\frac{u}{a}-k\Delta_N \Big
)-\psi \Big
(\frac{p\Delta_N}{a}-k\Delta_N \Big  ) \right ) \\
&&\hspace{2.3cm} \times  \left ( \psi \Big (\frac{u'}{a}-k\Delta_N
\Big )-\psi \Big
(\frac{p'\Delta_N}{a}-k\Delta_N \Big  ) \right )S(0,p\Delta_N,0,p'\Delta_N)
%%\\
%%&\hspace{-0.3cm} \le \hspace{-0.3cm}& C\,\sum_{p=0}^{N-1}
%%\sum_{p'=0}^{N-1}\int_{p\Delta_N}^{(p+1)\Delta_N }\hspace{-0.3cm}
%%\int_{p'\Delta_N}^{(p'+1)\Delta_N }\hspace{-1cm}dudu' \left ( \psi
%%\Big (\frac{u}{a}-k\Delta_N \Big )-\psi \Big
%%(\frac{p\Delta_N}{a}-k\Delta_N \Big  ) \right ) \\
%%&&  \times  \left ( \psi \Big (\frac{u'}{a}-k\Delta_N
%%\Big )-\psi \Big (\frac{p'\Delta_N}{a}-k\Delta_N \Big  ) \right)\cdot
%%\big ((p\Delta_N)^{H_0} +
%%(p\Delta_N)^{H_K} \big  )\cdot
%%\big ((p'\Delta_N)^{H_0} +
%%(p'\Delta_N)^{H_K} \big  )
\\
&\hspace{-0.3cm}  \leq
\hspace{-0.3cm}& C \left ( \sum_{p=0}^{N-1}
 \int_{0}^{\Delta_N } \hspace{-0.5cm}du \left |\psi \Big
(\frac{u+p\Delta_N}{a}-k\Delta_N \Big )-\psi \Big
(\frac{p\Delta_N}{a}-k\Delta_N \Big  ) \right | \times\big|
(p\Delta_N)^{H_0} +
(p\Delta_N)^{H_K} \big|  \right) ^2 \\
&\hspace{-0.3cm}  \leq \hspace{-0.3cm}& C \left ( \sum_{p=0}^{N-1}
 \int_{0}^{\Delta_N } \hspace{-0.5cm}du  \cdot \frac u a \cdot \sup_{t\in[0,\Delta_N]} \left |\psi' \Big
(\frac{t+p\Delta_N}{a}-k\Delta_N \Big ) \right |
\times\big|(p\Delta_N)^{H_0} +
(p\Delta_N)^{H_K} \big| \right) ^2
 \\
&\hspace{-0.3cm}  \leq \hspace{-0.3cm}& C \Delta_N^4 \left (
\sum_{p=0}^{N-1}
\sup_{t\in[0,\Delta_N]} \left |\psi' \Big
(\frac{t+p\Delta_N}{a}-k\Delta_N \Big ) \right |
\times \big |(p\Delta_N)^{H_0} +
(p\Delta_N)^{H_K} \big| \right) ^2.
\end{eqnarray*}
But  Assumption (A1) implies that for   $m=4$, there exists a
constant $C>0$ such that  for all $x \in \R$, $|\psi'(x)|
\leq C
\cdot (1+|x|)^{-m}$. We deduce
\begin{multline*}
\Delta_N \sum_{p=0}^{N-1} \sup_{t\in[0,\Delta_N]} \left |\psi'
\Big (\frac{t+p\Delta_N}{a}-k\Delta_N \Big ) \right |
\times \left|(p\Delta_N)^{H_0} +
(p\Delta_N)^{H_K} \right|  \\ \leq  C \Delta_N
\sum_{p=-\infty}^{\infty} \frac 1 {\Big ( 1 + \left | p\Delta_N
\right |  \Big )^{m}} \cdot \left|(p\Delta_N)^{H_0} +
(p\Delta_N)^{H_K} \right| \leq C \cdot \int_{-\infty}^\infty \frac
{|x|^{H_0}+ |x|^{H_K}}{(1+ |x|)^m} \, dx\,<\,\infty.
\end{multline*}
Therefore,
\begin{eqnarray}
\label{maj:I2}  \E \big (
I_{2,N}^2(a,k)\big ) &\le&
 C \Delta_N^2.
\end{eqnarray}
{\bf (2) Bound of $\displaystyle{\E
\varepsilon_{2,N}^2(a,k)}.$} \\
By using  Lemma
\ref{maj:cov} and  Cauchy-Schwartz inequality, we deduce that for $N$ large enough,
\ban \E \varepsilon_{2,N}^2(a,k)) &=&
\int_{T_N}^{\infty}\int_{T_N}^{\infty}
\psi\left(\frac{u}{a}-k\Delta_N  \right)
\,\psi\left(\frac{u'}{a}-k\Delta_N  \right) \,S(0,u,0,u')
\,du\,du' \\
&\le& C  \left ( \int _{T_N}^{\infty}
\left(1+\left|\frac{u}{a}-k\Delta_N\right|\right)^{-m}\,u^{2H_0}
du \right)^2. \ean On one hand,   $k \in D_N(a)$  implies that $k
\leq [(1-r)N/a]$. On the other hand,    $u \geq T_N=N \cdot \Delta_N$. Therefore, we
have$\displaystyle{\left(1+\left|\frac{u}{a}-k\Delta_N\right|\right)\ge
\left( u-(1-r) N\Delta_N\right)/a}$.
This implies that for $m\ge 4$ and $N$ large enough,
\ban \int _{T_N}^{\infty}
\left(1+\left|\frac{u}{a}-k\Delta_N\right|\right)^{-m}\,u^{2H_0}
du & \leq &  \int _{T_N}^{\infty} \left(\frac{u}{a}-(1-r) \frac {
N\Delta_N}a \right)^{-m}\,u^{2H_0} du  \\
&=& \frac{a^m}{(N \Delta_N) ^{m-2H_0-1}}\,\int_{r}^{\infty}
\frac{ (v+1-r)^{2H_0}}{v^m}\, dv,
 \ean
by making the change of variable $u= (N \Delta_N)\,(v+1-r)$.
Consequently,
\begin{eqnarray}
\label{maj:epsilon2} \E\varepsilon_{2,N}^2(a,k) \leq C \cdot \frac 1 {(N \Delta_N)^2}.
\end{eqnarray}
~\\
{\bf (3) Bound of $\displaystyle{\E\varepsilon_{3,N}^2(a,k)}.$\;} \\
By using the same kind of argument than in (2), one obtains that
for $N$ large enough
 \ban \E \varepsilon_{3,N}^2(a,k)) &=&
\int_{\infty}^0\int_{\infty}^0 \psi\left(\frac{u}{a}-k\Delta_N
\right) \,\psi\left(\frac{u'}{a}-k\Delta_N  \right) \,S(0,u,0,u')
\,du\,du' \\
&\le& C  \left ( \int_{\infty}^0
\left(1+\left|\frac{u}{a}-k\Delta_N\right|\right)^{-m}\,(|u|^{2H_K}+|u|^{2H_0})
du \right)^2 \\
&\le& C  \left ( \int_{r/2 \cdot N \Delta_N}^\infty \frac
{1+(v-r/2 \cdot N \Delta_N)^2 }{v^m} \, dv \right)^2
 \ean
As a consequence, for $m\geq 4$ and $N$ large enough,
\begin{eqnarray}
\label{maj:epsilon3} \E
\varepsilon_{3,N}^2(a,k) \leq C \cdot \frac 1 {(N \Delta_N)^2}.
\end{eqnarray}
Finally, from (\ref{maj:I1}), (\ref{maj:I2}),
%%(\ref{maj:epsilon1}),
(\ref{maj:epsilon2}) and (\ref{maj:epsilon3}), we deduce that
(\ref{maj:unif:Eeps2}) holds. This finishes the proof of the point
i).  Since $N \Delta_N \to \infty$ and $N \Delta_N ^2 \to 0$,
$(N\Delta_N) \varphi(N)$  converges to $0$ when $N\to \infty$.
This proves the point iii). To complete the proof of Lemma
\ref{maj:sum:epsilon}, it remains to proves the point ii). We
deduce from the decomposition (\ref{dec:epsi}) that \ba
\E\left|\sqrt{a_1}\varepsilon_N(a_1,k
)-\sqrt{a_2}\varepsilon_N(a_2,k )\right|^2 &\le& 3\,
 \sum_{i=1}^3 \E\left|\varepsilon_{i,N}(a_2,k )-\varepsilon_{i,N}(a_1,k )\right|^2
\ea The same calculations than the ones used to prove the point i)
provide the upper bound on the terms \\
$\displaystyle{\E\left|\varepsilon_{i,N}(a_2,k
)-\varepsilon_{i,N}(a_1,k )\right|^2}$. Indeed, consider for
instance the terms with $\varepsilon_{2,N}$, then by using Taylor
formula, for every pair $(a_1,a_2)$ with $a_{min}\le a_1<a_2\le
a_{max}$ there exists a real number $\theta \in (a_1,\,a_2)$ such
that \ban \varepsilon_{2,N}(a_2,k )-\varepsilon_{2,N}(a_1,k ) &=&
\left(a_2-a_1\right) \times
\int_{T_N}^{\infty}\left(\frac{-t}{\theta^2}\right)\psi'\left(\frac{t}{\theta}-k\Delta_N\right)\,X(t)\,dt
\ean Next, by using the same kind of arguments than for the bound
of $\displaystyle{\E \varepsilon_{2,N}^2(a,k)}$
 in point {\em i)}, we get that for every integer $m>4$, every
$a_1, a_2$ in $[a_{min},\,a_{max}]$ and $k\in D_N(a_1)\cap
D_N(a_2)$
\ban
\E \left|\varepsilon_{2,N}(a_2,k )-\varepsilon_{2,N}(a_1,k
)\right|^2 &\le& \frac{C}{a^4_{min}}\,\left|a_2-a_1\right|^2
\left(\int _{T_N}^{\infty}
\left(1+\left|\frac{u}{\theta}-k\Delta_N\right|\right)^{-m}\,u^{1+2H_0}
du\right)^2
\\
&\le& C \,\left|a_2-a_1\right|^2 \,(N \Delta_N)^{-2}.
\ean
We deduce similarly that
\ban
\E \left|\varepsilon_{3,N}(a_2,k )-\varepsilon_{3,N}(a_1,k
)\right|^2 &\le& C \,\left|a_2-a_1\right|^2 \,(N \Delta_N)^{-2}.
\ean
At this point, it remains to show
\ba
\label{maj:Delta:eps1}
\E \left|\varepsilon_{1,N}(a_2,k )-\varepsilon_{1,N}(a_1,k
)\right|^2 &\le& C \,\left|a_2-a_1\right|^2 \,\varphi(N)
\ea
to finish the proof of item ii). But, we have the decomposition
\ban
\E \left|\varepsilon_{1,N}(a_2,k )-\varepsilon_{1,N}(a_1,k
)\right|^2 &\le& 2\,
\E \left|I_{1,N}(a_2,k )-I_{1,N}(a_1,k
)\right|^2+ 2\,\E \left|I_{2,N}(a_2,k )-I_{2,N}(a_1,k )\right|^2
\ean
However,   Taylor Formula implies the existence of two real
numbers $\theta_1, \,\theta_2 \in (a_1,\,a_2)$ such that
$$I_{i,N}(a_2,k )-I_{1,N}(a_1,k ) = (a_2-a_1)\cdot
\widetilde{I}_{i,N}(\theta_i,k )\qquad \mathrm{for}\;i=1\; \mathrm{or} \;2$$  where $\widetilde{I}_{i,N}(a,k
)$ is obtained by replacing into the expression of $I_{i,N}(a,k )$
the map $\displaystyle{\psi \Big (\frac{t}{a}-k\Delta_N \Big  )}$
by the map $\displaystyle{\left(\frac{-t}{a^2}\right)\times\psi'
\Big (\frac{t}{a}-k\Delta_N \Big  )}$ and
$\displaystyle{\psi \Big (\frac{p\Delta_N}{a}-k\Delta_N \Big  )}$
by $\displaystyle{\left(\frac{-p\Delta_N}{a^2}\right)\times\psi'
\Big (\frac{p\Delta_N}{a}-k\Delta_N \Big  )}$. So,
\ban
\E \left|\varepsilon_{1,N}(a_2,k )-\varepsilon_{1,N}(a_1,k
)\right|^2 &\le& C
\,\left|a_2-a_1\right|^2\times\left\{\E\widetilde{I}^2_{1,N}(\theta_1,k
) \,+\, \E\widetilde{I}^2_{2,N}(\theta_2,k )
\right\}
\ean
Since the map $\displaystyle{t\mapsto
\left(\frac{t}{a}\right)\times\psi'
\Big (\frac{t}{a}-k\Delta_N \Big  )}$ is still continuously differentiable and fast decreasing,
 one can lead same calculations that in
the bound of $\E I^2_{2,N}(a,k )$. We finally get
$\;\displaystyle{\E\widetilde{I}^2_{2,N}(\theta_2,k )+
\E\widetilde{I}^2_{2,N}(\theta_2,k ) \le\, C\, \varphi(N)
}$. This implies (\ref{maj:Delta:eps1}) and completes the proof of
Lemma \ref{maj:sum:epsilon}.
\end{dem}
\begin{lem}\label{maj:EI-J}
Under the same assumptions as in Theorem \ref{Theo:TCL:discret},
there exists a positive constant $C>0$ such that for every real
number $a>0$ and $N\in\N^*$, we have $\displaystyle{\E
\left|I_N(a) - J_N(a)
\right| \le C
\times\varphi(N)^{1/2}}$.
\end{lem}
\begin{dem}
Since the variables $d  = d (a,k \Delta_N)$ and $e= e(a,k
\Delta_N)$ are Gaussian, the variables $d^2-e^2\, $ have finite
second order moment and Jensen's inequality implies
 \ban
\E \left|I_N(a) - J_N(a) \right| &=& \E \left|\frac{1}{|D_N(a)|}\,
\sum_{k \in D_N(a)}  \left(d^2(a,k \Delta_N)-e^2(a,k
\Delta_N)\right)\right|
\\
&\le& \frac{1}{|D_N(a)|}\,\sum_{k \in D_N(a)}  \sqrt{\E
\left(d^2(a,k
\Delta_N)-e^2(a,k \Delta_N)\right)^2}
 \ean
Then we derive an upper bound for the expectations $\E
\left(d^2(a,k \Delta_N)-e^2(a,k \Delta_N)\right)^2$. Indeed, $d$
and $e$ are jointly Gaussian variables with zero means. One has
$$\E(d^2-e^2)^2 = \E(d-e)^2 (d+e)^2 =: \E \varepsilon^2 Z^2,
$$ where $\varepsilon = d-e$ and $Z = d+e$ are also jointly Gaussian and
have mean zero. By using that $Z  = \sigma_2\, \sigma_1^{-1}
\rho\, \varepsilon + \xi$, where $\sigma^2_1 = \E \varepsilon^2$,
$\sigma^2_2 = \E Z^2$,   $\rho= corr(\varepsilon,\,Z)$ and where
$\xi$ is independent of $\varepsilon$ and Gaussian, one can show
that \ban
 \E \varepsilon^2 Z^2&=& (\E \varepsilon^4)
\frac{\sigma_2^2\rho^2}{\sigma_1^2}\,+\,\E \varepsilon^2\, \E
\xi^2\,=\, 3\sigma_1^2 \sigma_2^2 \rho^2 + \sigma_1^2 \sigma_2^2
(1-\rho^2) \,\le\,3 \sigma_1^2 \sigma_2^2
\\
&=& 3\,\E \varepsilon^2 \times \E(d+e)^2.\ean But $(d+e)^2 =
(2d-\varepsilon)^2 \le 8 \,d^2 + 2\varepsilon^2$, therefore
\ban
\E \left|I_N(a) - J_N(a) \right| &\le & \frac{\sqrt{6}}{|D_N(a)|}\,\sum_{k
\in D_N(a)} \sqrt{\E \varepsilon^2_N(a,k )} \times \sqrt{\E\,\left [
4\,d^2(a,k \Delta_N) + \varepsilon_N^2(a,k)\right ] }
\\
&\le& \frac{\sqrt{6}}{|D_N(a)|}\, \left\{ \sum_{k \in D_N(a)}
\E
\varepsilon_N^2(a,k) \right\}^{1/2} \times\left\{\sum_{k \in D_N(a)}
\E\left[4\,d^2(a,k \Delta_N) + \varepsilon_N^2(a,k)\right]
\right\}^{1/2}
\\
&\leq & \sqrt{6} \,\left\{\frac{1}{|D_N(a)|}\sum_{k \in D_N(a)}
\E
\varepsilon_N^2(a,k) \right\}^{1/2}\times  \left\{4
\,\mathcal{I}_{1}(a)+\frac{1}{|D_N(a)|}\sum_{k \in D_N(a)} \E
\varepsilon_N^2(a,k) \right\}^{1/2} , \ean where the two last
inequalities follow from Cauchy-Schwartz inequality and $\E
d^2(a,k \Delta_N) = \mathcal{I}_{1}(a)$ for every integer $k$.
Since $\widehat{\psi}$ is compactly supported, then
$\displaystyle{\sup_{a\in[a_{min},\, a_{max}]}
|\mathcal{I}_{1}(a)| <
\infty}$. By combining this remark with Lemma \ref{maj:sum:epsilon}  i), this provides
$\displaystyle{\E\left|I_N(a) - J_N(a)
\right| \le C\times\varphi(N)^{1/2}}$
 and finishes
the proof of the lemma.
\end{dem}
~\\ \noindent Now, the following proof of Theorem
\ref{Theo:TCL:discret} can be established~: \\
\begin{dem}[Theorem \ref{Theo:TCL:discret}]
From Lemma \ref{maj:EI-J} combined with Lemma
\ref{maj:sum:epsilon}   iii), we deduce \be\label{lim:NdeltaN:EI-J}
\lim_{N\to \infty} (N \, \Delta_N )^{1/2}  \E \left|I_N(a) -
J_N(a) \right|= 0. \ee Combined with (\ref{TLC}), this implies the
convergence of the finite-dimensional distribution in
(\ref{TLC2}). Indeed, it suffices to show that \be \label{ref:3.8}
\sqrt{N \, \Delta_N}\, \left(\log J_N(a) - \log I_N(a)\right)
\limiteproba 0. \ee Let $\varepsilon
>0$. By using the inequality $\displaystyle{\left|\log(x) -
\log(y)\right|\le 2\, |x/y-1|}$, valid for all $|x/y-1| \le 1/2$,
$x,y >0$ one can show that
\ba \nonumber \Pr\left( \sqrt{N \,
\Delta_N}\, \left|\log J_N(a) - \log I_N(a)\right| \ge \varepsilon
\right) && \\
\nonumber \hspace{-3cm}&&\hspace{-3cm}\leq \Pr\left(2 \, \left|
J_N(a)/ I_N(a) - 1\right|\ge \frac{\varepsilon}{\sqrt{N \,
\Delta_N}}\right) +\Pr\left(  \left|  J_N(a)/I_N(a) - 1\right|>
\frac{1}{2} \right)
\\
\label{ref:3.9} \hspace{-3cm}&&\hspace{-3cm}\leq 2\,\Pr\left(
\left| J_N(a) - I_N(a)\right|\ge \frac{\varepsilon
\,I_N(a)}{2\,\sqrt{N \, \Delta_N}} \right)
\\
\nonumber \hspace{-3cm}&&\hspace{-3cm}\leq 2\,\Pr\left(\left|
J_N(a) - I_N(a)\right|\le \frac{\varepsilon
\,\mathcal{I}_1(a)}{4\,\sqrt{N \, \Delta_N}}
\;\mathrm{and}\;\left| J_N(a) -   I_N(a)\right|\ge
\frac{\varepsilon \,I_N(a)}{2\,\sqrt{N \, \Delta_N}} \right)
\\
\nonumber &&\hspace{2cm}+2\,\Pr\left( \left| J_N(a) -
I_N(a)\right|\ge \frac{\varepsilon \,\mathcal{I}_1(a)}{4\,\sqrt{N
\, \Delta_N}} \right)
\\
\label{ref:3.10} \hspace{-3cm}&&\hspace{-3cm}\leq 2\,\Pr\left( \left|
J_N(a) -   I_N(a)\right|\ge \frac{\varepsilon
\,\mathcal{I}_1(a)}{4\,\sqrt{N \, \Delta_N}} \right) +
2\,\Pr\left(I_N(a) \le \frac{\mathcal{I}_1(a)}{2} \right).
\ea
The second inequality in (\ref{ref:3.9}) is valid for all $N$ such
that $ \varepsilon /\sqrt{N \Delta_N}\le 1/2$, that is, for all
sufficiently large $N$. The second term in the right-hand side of
(\ref{ref:3.10}) vanishes, as $N\to \infty$, because
$\displaystyle{I_N(a)\limiteproba \mathcal{I}_1(a)}$. By using the
Markov inequality, one can bound above the first term in the
right-hand side of (\ref{ref:3.10}) by
$$
\frac{8\,\sqrt{N \, \Delta_N}}{\varepsilon \,\mathcal{I}_1(a)} \,
\E\left| J_N(a) - I_N(a)\right|.
$$
Thus, from (\ref{lim:NdeltaN:EI-J}), one obtains Relation
(\ref{ref:3.8}), which completes the proof of the convergence of
the finite distributions. To finish with the proof of Theorem
\ref{Theo:TCL:discret}, we have to show the tightness of the
sequence $\left(L_N(a)\right)_{a_{min}\le a \le a_{max}}$ where
$\displaystyle{L_N(a) = \sqrt{N \, \Delta_N} \, \Big(J_N(a)
-\mathcal{I}_1(a)\Big)}$. Observe one has the decomposition
$L_N(a) =  L_{1,N}(a) + L_{2,N}(a)$ with $\displaystyle{\left \{
\begin{array}{ccl} L_{1,N}(a) &= &
\sqrt{N \, \Delta_N}\, \Big(I_N(a) -\mathcal{I}_1(a)\Big) \\
L_{2,N}(a) &=& \sqrt{N \, \Delta_N} \, \Big(J_N(a) -I_N(a)\Big)
\end{array} \right . }$. In \cite{BaBe}, one have proved the tightness and the weak convergence of
$\left(L_{1,N}(a)\right)_{a_{min}\le a \le a_{max}}$ in Skorokhod
topology on the space of c\`ad-l\`ag functions on $[a_{min},a_{max}]$.
From the other hand, (\ref{lim:NdeltaN:EI-J}) implies that
$L_{2,N}(a) \limiteloi 0$ for all $a\in [a_{min},a_{max}]$. Note
that the limit process is null, thus it is obviously continuous.
Then, provided one have shown the tightness of
$\left(L_{2,N}(a)\right)_{a_{min}\le a \le a_{max}}$, one can
deduce the tightness of $\left(L_{N}(a)\right)_{a_{min}\le a \le
a_{max}}$, see for instance Jacod and Shyriaev, Cor 3.33, p. 317.
Next, one deduce the weak convergence of $ L_{N}(a)$ to $Z(a)$ in
the Skorokod topology on the space of c\`ad-l\`ag functions on
$[a_{min},a_{max}]$. The last step is the proof of the tightness
of $\left(L_{2,N}(a)\right)_{a_{min}\le a \le a_{max}}$. Following
Ikeda and Watanabe, Th.4.3 p. 18, it suffices to show the
existence of a positive constant $M_2$ such that for all
$a_1,\,a_2 \in [a_{min},\, a_{max}]$ \be \label{CS:tension}
\E\left(L_{2,N}(a_2)-L_{2,N}(a_1)\right)\le M_2\,
\left|a_2-a_1\right|^2. \ee Nowever, from (\ref{def:IN(a)}) and
(\ref{def:JN(a)}), we get \ban J_N(a) - I_N(a) &=&
\left|D_N(a)\right|^{-1} \sum_{k\in D_N(a)}
\left(e^2(a,k)-d^2(a,k)\right). \ean Therefore, for $a_1<a_2$, we
have \ban L_{2,N}(a_2)-L_{2,N}(a_1) &=& \Big( N \,
\Delta_N\Big)^{1/2}\times \left\{\sum_{k=[rN/a_2]}^{[rN/a_1]}
\left|D_N(a_2)\right|^{-1}\left(e^2(a_2,k)-d^2(a_2,k)\right)
\right.
\\
&+&  \sum_{k=[rN/a_1]}^{[(1-r)N/a_2]}
\left[\left|D_N(a_2)\right|^{-1}\left(e^2(a_2,k)-d^2(a_2,k)\right)
-\left|D_N(a_1)\right|^{-1}\left(e^2(a_1,k)-d^2(a_1,k)\right))\right]
\\
&+& \left. \sum_{k=[(1-r)N/a_1]}^{[(1-r)N/a_2]}
\left|D_N(a_1)\right|^{-1}\left(e^2(a_1,k)-d^2(a_1,k)\right)
\right\} \ean Then, one remarks that for any finite family $I$ of
random variables $\left(X_i\right)_{i\in I}$ with finite variance
we have \\$\displaystyle{\E \left(\sum_{i\in I}X_i\right)^2\,=\,
\sum_{(i,j)\in I^2}\E\left(X_i\,X_j\right)\,\le\, \sum_{(i,j)\in
I^2}\sqrt{\E X_i^2}\times\sqrt{\E X_j^2}}\,=\, \left(\sum_{i\in
I}\sqrt{\E X_i^2}\right)^2\;$ which combined with
\\$\displaystyle{(x+y+z)^2 \,\le\, 3\, (x^2+y^2+z^2)}\;$ implies
\ban \E \left|L_{2,N}(a_2)-L_{2,N}(a_1)\right|^2 &\le& C\,\left( N
\, \Delta_N\right)\times\left(S_1^2+S_2^2+S_3^2\right) \ean where \ban
S_1&=&\left|D_N(a_2)\right|^{-1} \sum_{k=[rN/a_2]}^{[rN/a_1]}
\sqrt{\E\left(e^2(a_2,k)-d^2(a_2,k)\right)^2},
\\
S_2&=&\left|D_N(a_1)\right|^{-1}\,\sum_{k=[(1-r)N/a_1]}^{[(1-r)N/a_2]}
\sqrt{\E\left(e^2(a_1,k)-d^2(a_1,k)\right)^2},
\\
S_3&=& \sum_{k=[rN/a_1]}^{[(1-r)N/a_2]}
\sqrt{\E\left[\left|D_N(a_2)\right|^{-1}\left(e^2(a_2,k)-d^2(a_2,k)\right)
-\left|D_N(a_1)\right|^{-1}\left(e^2(a_1,k)-d^2(a_1,k)\right)\right]^2}.
\ean
From (\ref{def:epsilon}), we get $\displaystyle{e^2(a,k)-d^2(a,k)
= \varepsilon^2(a,k) + 2\, \varepsilon(a,k)\,d(a,k)}$. Moreover,
the random  variables $X=\varepsilon(a,k)$ or $X=d(a,k)$ are
centred Gaussian random variables, thus we have
$\displaystyle{\sqrt{\E\left(X^4\right)} =
\sqrt{3}\,
\E\left(X^2\right)}$. Then, by combining this remark with
Cauchy-Schwarz inequality and Lemma \ref{maj:sum:epsilon}, we
deduce
\ban
\E\left(e^2(a,k)-d^2(a,k)\right)^2 &\le& C\,\left\{\E \varepsilon^4(a,k) + \E  \varepsilon^2(a,k)\,d^2(a,k) \right\}
\\
&\le&
 C\,\left\{\E\varepsilon^4(a,k) + \sqrt{\E  \varepsilon^4(a,k)}\times\sqrt{\E d^4(a,k)} \right\}
\\
&\le& C\,\left\{\Big(\E\varepsilon^2(a,k) \Big)^2+
\E\varepsilon^2(a,k)\times\E d^2(a,k)\right\}
\\
&\le& C\times\varphi(N)\times
\Big\{\mathcal{I}_1(a)+\varphi(N)\Big\}.
\ean
Afterwards
\ba
\nonumber
\left( N \,
\Delta_N\right)\,S_1^2 &\le& C\,\left( N \,
\Delta_N\right)\times \varphi(N)\times
\Big\{\mathcal{I}_1(a)+\varphi(N)\Big\}\times\left|D_N(a_2)\right|^{-2}\times
\Big\{[rN/a_1]-[rN/a_2]\Big\}^2
\\
\label{maj:S1}
&\le&
 C\,\left( N \,
\Delta_N\right)\times \varphi(N)\times\left|a_2-a_1\right|^2.
\ea
since $\left|D_N(a_2)\right|\sim (1-2r)N/a$ as $N$ goes to
$\infty$ and $\mathcal{I}_1(a)$ is bounded. The same calculations
provide
\ba
\label{maj:S2}\left( N \,
\Delta_N\right)\,S_2^2
&\le&
 C\,\left( N \,
\Delta_N\right)\times \varphi(N)\times\left|a_2-a_1\right|^2.
\ea
Next, we derive the upper bound for $S_3^2$. Let us stress that
the functions $\displaystyle{a\mapsto N\,|D_N(a)|^{-1}}$ converges
uniformly to $a (1-2r)^{-1}$ when $N$ goes to $\infty$. Thus one
can replace  $S_3^2$ by $\widetilde{S}_3^2$ where
\ban
\widetilde{S}_3&=&\frac{1} {(1-2r)N}
 \sum_{k=[rN/a_1]}^{[(1-r)N/a_2]}
\sqrt{\E\left[ a_2\,\left(e^2(a_2,k)-d^2(a_2,k)\right)
- a_1\,\left(e^2(a_1,k)-d^2(a_1,k)\right)\right]^2}.
\ean
Then,  by using (\ref{def:epsilon}), we get the following
expansion of the term $f_k$ define below
\ban
f_k&:=&a_2\,\left(e^2(a_2,k)-d^2(a_2,k)\right) -
a_1\,\left(e^2(a_1,k)-d^2(a_1,k) \right)
\\&=&
\Big(\sqrt{a_2}\,\varepsilon(a_2,k)-\sqrt{a_1}\,\varepsilon(a_1,k)\Big)
\times\Big(\sqrt{a_2}\,\varepsilon(a_2,k)+\sqrt{a_1}\,\varepsilon(a_1,k)\Big)
\\&&+
2\,\sqrt{a_2}\,
d(a_2,k)\times\Big(\sqrt{a_2}\,\varepsilon(a_2,k)-\sqrt{a_1}\,\varepsilon(a_1,k)\Big)
+\sqrt{a_1}\,\varepsilon(a_1,k)\times\Big(\sqrt{a_2}\,d(a_2,k)-\sqrt{a_1}\,d(a_1,k)\Big)
\ean
We lay the emphasize on the fact that all the random variables in
the above formula  are Gaussian centred variables. But for two
Gaussian centred random variables, say $X$ and $Y$, we get $\E
\left(X^2Y^2\right) \le  \sqrt{\E X^4}\times \sqrt{\E Y^4} = 3 \left(\E
X^2\right)\times \left(\E Y^2\right)$. By combining this remark
with
 Lemma
\ref{maj:sum:epsilon}, one obtains
\ban
\E f_k^2 &\le& C\,\left\{
\E\Big(\sqrt{a_2}\,\varepsilon(a_2,k)-\sqrt{a_1}\,\varepsilon(a_1,k)\Big)^2
\times\E\Big(\sqrt{a_2}\,\varepsilon(a_2,k)+\sqrt{a_1}\,\varepsilon(a_1,k)\Big)^2\right.
\\&&\hspace{-0.9cm}+
\left.a_2\,\E
d^2(a_2,k)\times\E\Big(\sqrt{a_2}\,\varepsilon(a_2,k)-\sqrt{a_1}\,\varepsilon(a_1,k)\Big)^2
+a_1\,\E\varepsilon^2(a_1,k)\times\E\Big(\sqrt{a_2}\,d(a_2,k)-\sqrt{a_1}\,d(a_1,k)\Big)^2\right\}
\\
&\le& C\,
\varphi(N)\,\left|a_2-a_1\right|^2\,\left\{a_1\,\mathcal{I}_1(a_1)+(a_1+a_2)\,\varphi(N) \right\}
+ C\,
\varphi(N)\,\E\Big(\sqrt{a_2}\,d(a_2,k)-\sqrt{a_1}\,d(a_1,k)\Big)^2.
\ean
But Taylor Formula implies the existence of a real numbers
$\theta_t\in (a_1,\,a_2)$ such that
\ban
\sqrt{a_2}\,d(a_2,k)-\sqrt{a_1}\,d(a_1,k) &=&
 \big(a_2-a_1 \big)
 \, \int_{\R} \left(\frac{-t}{\theta_t^2}\right)\,
\psi'\left(\frac{t}{\theta_t}-k\Delta_N\right)\,X(t)\,dt
\ean
and after
$\displaystyle{\E\Big(\sqrt{a_2}\,d(a_2,k)-\sqrt{a_1}\,d(a_1,k)\Big)^2\,
\le\, C\, \big|a_2-a_1 \big|^2 }$. Indeed,  one observe that since $\theta_t \in (a_{min},a_{max})$,
one haves $1/\theta_t^2\le 1/a_{min}^2$. This implies
\ban
\E\Big(\sqrt{a_2}d(a_2,k)-\sqrt{a_1}d(a_1,k)\Big)^2
&=& |a_2-a_1 |^2  \int_{\R} \int_{\R} \frac{uv} {\theta_u^2
\theta_v^2}\,\cdot
\psi'\left(\frac{u}{\theta}-k\Delta_N\right)\,\psi'\left(\frac{v}{\theta}-k\Delta_N\right)\,S(0,u,0,v)\,du\,dv
\\
&\le& \hspace{-0.3cm} \frac{|a_2-a_1
|^2}{a_{min}^4}\left(\hspace{-0.1cm}\int_{\R}
\sup_{\theta\in(a_{min}, a_{max})}\hspace{-0.1cm}\left|\psi'\left(\frac{u}{\theta}-k\Delta_N\right)\right|\,\left ( |u|^{1+H_K}    \1_{u\leq 1} +
|u|^{1+H_0}\right) \,du\hspace{-0.1cm}\right)^2.
\ean
On the other hand, the fast decreasing of the function $\psi'$
insures    $$\int_{\R}
\sup_{\theta\in(a_{min},\,a_{max})}\left|\psi'\left(\frac{u}{\theta}-k\Delta_N\right)\right|\,\left ( |u|^{1+H_K}    \1 _{u\leq 1} +
|u|^{1+H_0}\right) \,du<\infty.$$ Therefore, since $a_1$, $a_2$,
$\mathcal{I}_1(a)$ are bounded
 and $\varphi(N)\to 0$ as $N$ goes to $\infty$,  we have
$\;\displaystyle{\E f_k^2
\,\le\,C\,\varphi(N)\,\left|a_2-a_1\right|^2}\;$.
This leads to
\ban
\left( N \,
\Delta_N\right)\cdot\widetilde{S}^2_3&\le& C\,N^{-1}\cdot\left( N \,
\Delta_N\right)\times
\varphi(N)\,\left|a_2-a_1\right|^2\Big([(1-r)N/a_2]-[rN/a_1]\Big)
\\
&\le& C\cdot\left( N \, \Delta_N\right)\times
\varphi(N)\,\left|a_2-a_1\right|^2. \ean Eventually, combined with
(\ref{maj:S1}, \ref{maj:S2}), one obtains \ban \E
\left|L_{2,N}(a_2)-L_{2,N}(a_1)\right|^2 &\le& C\cdot\left( N \,
\Delta_N\right)\times \varphi(N)\,\left|a_2-a_1\right|^2. \ean But,
Lemma \ref{maj:sum:epsilon} iii) implies that $\left( N \,
\Delta_N\right)\times \varphi(N)$ converges to $0$ when $N$ converges
to $\infty$, therefore we deduce (\ref{CS:tension}). This finishes
the proof of the tightness of the sequence
$\left(L_N(a)\right)_{a_{min}\le a \le a_{max}}$. Now, the
functional Delta method (see for instance Van der Vaart, chapter
20, p. 297), provide a central limit theorem for
$\log(I_N(.))-\log({\cal I}_1(.))$, because the function $\log(.)$
is a Hadamard-differentiable function on the space of c\`ad-l\`ag
function on $[a_{min},a_{max}]$; this completes the proof of
Theorem \ref{Theo:TCL:discret}.
\end{dem}
\subsection{Proofs of section 4}
\begin{dem}[{\bf Proposition \ref{omegas}}] %%\\~\\
\quad{\bf We lay the emphasize on the fact that, in this proof, we generalize the choice of the frequencies by considering $a_N=(N \Delta_N)^q,~\mbox{with}~q>0.$} \\
~\\
For a given $N$, denote
$T^*=(t^*_0=0,t_1^*,\cdots,t_K^*,t^*_{K+1}=a_N)$ such as~:
$$f_{t_j^*} <\frac { \omega_j^*}{\alpha} \leq f_{t_j^*+1} ,\quad \mathrm{for\; all}
\quad j=1,\dots,K$$ and for $T=(0,t_1,\cdots,t_K,a_N) \in {\cal A}_K
^{(N)}$, we denote $\displaystyle{Z_i ^{(N)}=\sqrt{N
\Delta_N}(Y_i-\log {\cal I}_1(1/f_i))}$, \\
$Y_{]t_j,t_{j+1}]}=(Y_{t_j+1},\cdots,Y_{t_{j+1}-\tau_N})'$,
$X_{]t_j,t_{j+1}]}=(\log f_{t_j+i},1)_{1\leq i \leq (t_{j+1}-t_j)}
$, $Z ^{(N)}_{]t_j,t_{j+1}]}=(Z ^{(N)}_{t_j+1},\cdots,Z
^{(N)}_{t_{j+1}-\tau_N})'$.\\ \\
{\bf First step :} We would like to prove~:
$\widehat{\omega}_j^{(N)}
 \limiteproba \omega_j^*$ for all $j=1,\dots,K$. \\
Denote $Q^{(N)}_*=Q^{(N)}(t^*,\widehat{\Lambda}(t^*))$ where
$\widehat{\Lambda}(t^*)$  is obtained from a linear regression of
$(Y_i)$ on $(\log f_i)$ for $i=t_j^*+1,\cdots,t_{j+1}^*-\tau_N$. Let
$\varepsilon >0$ and $\displaystyle{\parallel T-T' \parallel
_\infty=\max_{j \in \{1,\cdots,K \} }|t_j-t'_j| }$ for
$T=(0,t_1,\cdots,t_K,a_N) \in {\cal A}_K ^{(N)}$ and
$T'=(0,t'_1,\cdots,t'_K,a_N) \in {\cal A}_K ^{(N)}$. Then, we get,
$$\Pr   \left (\parallel \widehat {T}-t^* \parallel _\infty \geq
\varepsilon a_N  \right ) \leq  \Pr  \left ( \min _{T \in
V_{\varepsilon a_N}} Q^{(N)}(T,\widehat{\Lambda}(T)) \leq
Q^{(N)}_* \right ), $$ where $\displaystyle{V_{\varepsilon
a_N}=\left \{ T \in {\cal A}_K ^{(N)},\parallel T -t^*
\parallel_\infty \geq \varepsilon a_N \right \}}$. We want to show
that for all $T \in V_{\varepsilon a_N}$,\\ $Q^{(N)}_*=o(Q
^{(N)}(T,\widehat{\Lambda}(T)))$. In fact,
\begin{eqnarray*}
Q^{(N)}_* &\hspace{-3mm} = &\hspace{-3mm}  \frac 1 {N \Delta_N}
\sum_{j=0}^{K+1}
 (Z ^{(N)}_{]t_j^*,t_{j+1}^*]})'\left [ Id-X_{]t_j^*,t_{j+1}^*]}
 \left ( X'_{]t_j^*,t_{j+1}^*]}X_{]t_j^*,t_{j+1}^*]} \right )^{-1} X'_{]t_j^*,t_{j+1}^*]}\right ] Z ^{(N)}_{]t_j^*,t_{j+1}^*]} \\
&\hspace{-3mm}  \leq  &\hspace{-3mm}\frac 1 {N \Delta_N}
\sum_{j=0}^{K+1} (Z ^{(N)}_{]t_j^*,t_{j+1}^*]})'
Z^{(N)}_{]t_j^*,t_{j+1}^*]} \\
&\hspace{-3mm}  \leq  &\hspace{-3mm} \frac 1 {N \Delta_N}(Z
^{(N)}_{[1,a_N]})' Z^{(N)}_{[1,a_N]}.
\end{eqnarray*}
From Proposition \ref{propTLC}, we deduce
\begin{eqnarray}\label{IZ}
\frac 1 {a_N}(Z ^{(N)}_{[1,a_N]})' Z ^{(N)}_{[1,a_N]} \limiteloi
I_Z= \int _0 ^1 Z^2 \left ( \frac \beta {f_{min}} \left ( \frac
{\alpha f_{min}}{\beta f_{max}} \right )^u \right ) du ,
\end{eqnarray}
which is a positive and $\LL ^\infty$ random variable because $Z$
is a continuous Gaussian process. Afterward, for a sequence
$(\psi_k)_k \in \R^{\N}$ and a sequence of random variables
$(\xi_k)_{k\in \N}$, we will write $\xi_N=O_P(\psi_N)$ as $N \to
\infty$, if for all $\varepsilon >0$, there exists $c>0$, such as
,
$$
P\Big (|\xi_N| \leq c \cdot \psi_N \Big ) \geq 1 -\varepsilon,
$$
for all sufficiently large $N$. Here, we obtain~: \ba \label{Q*}
Q^{(N)}_* = O_P \left ( \frac {a_N}{N \Delta_N}\right ). \ea Now,
let $T \in V_{\varepsilon a_N}$, we want a lower bound of
$Q^{(N)}(T,\widehat{\Lambda}(T))$. We use the following
decomposition \ban Q^{(N)}(T,\widehat{\Lambda}(T))& =
&\sum_{j=0}^{K+1} \sum_{i=t_j+1}^{t_{j+1}-\tau_N } \left[Y_i -\log
{\cal I}_1 (1/f_i) \right]^2 +
\left[ X_i\widehat{\lambda}_j -\log  {\cal I}_1 (1/f_i) \right]^2 + \\
&& \hspace{4cm}2  \left[ Y_i- \log  {\cal I}_1
(1/f_i)\right]\times\left[ X_i\widehat{\lambda}_j - \log  {\cal I}_1
(1/f_i)\right]\\ & = &Q_1+Q_2+Q_3. \ean Then :
\begin{enumerate}
\item Since $\displaystyle{Q_1 = \frac {1}{N \Delta _N}
\sum_{j=0}^{K+1} (Z ^{(N)}_{]t_j,t_{j+1}]})'
Z^{(N)}_{]t_j,t_{j+1}]}}$, as previously we get \ba \label{Q1} Q_1
= O_P \left (\frac {a_N}{N \Delta _N} \right ). \ea
\item Let $\displaystyle{{\underline \tau}=\left ( \log
\left ( \frac {\beta f_{max} }{ \alpha f_{min}} \right ) \right
)^{-1} \min_{j=1,\cdots,K} \left \{ \log \left ( \frac {\alpha
\omega^*_{j+1} }{\beta \omega^*_{j} } \right )\right \}}$. Then, for
all $j \in \{0,1,\cdots,K\}$, $t^*_{j+1}-\tau_N \geq t^*_j +
{\underline \tau} a_N$. Since $T \in V_{\varepsilon a_N}$, we have
$ \eta =\min \{\varepsilon,{\underline \tau},\log (\beta / \alpha)
\} >0$ and there exists an integer $j \in \{0,\cdots,K+1\}$ for which
there are no estimated abrupt change  in the interval $[t^*_j-\eta
a_N,t^*_j]$ or $[t^*_j -\tau_N ,t^*_j -\tau_N+\eta a_N]$. Thus
there exists $k \in \{0,\cdots,K+1\}$ satisfying $[t^*_j-\eta a_N ,
t^*_j]\subset [t_k,t_{k+1}-\tau_N]$ (we follow here a similar
proof than Bai and Perron in Lemma 2, p 69) and \ba \nonumber Q_2
& \geq & \sum_{i=t^*_j-\eta a_N+1}^{t^*_j}  | X_i\widehat{\lambda}
_k -\log {\cal I}_1(1/f_i) |^2 \\
\label{Q21} & \geq &   \sum_{i=t^*_j-\eta a_N+1}^{t^*_j}   \left |
A(\widehat{H}_k^{(N)},\widehat{\sigma}_k^{(N)}) + \frac i {a_N} \cdot
B(\widehat{H}_k^{(N)},\widehat{\sigma}_k^{(N)})- g\left ( \frac i
{a_N} \right ) \right |^2 , \ea with~:
\begin{itemize}
\item $\displaystyle{A(H,\sigma)=\log \Big ( \sigma^2 \cdot
K_H(\psi) \Big)-(2H+1) \cdot \log \left ( \frac {f_{min}}{\beta}
\right )}$ for all $(H,\sigma)\in {\cal K}$;
\item $\displaystyle{B(H,\sigma)=-(2H+1) \cdot\log  \left (\frac {\beta f_{max}}{\alpha f_{min}}
\right )}$ for all $(H,\sigma)\in {\cal K}$;
\item $\displaystyle{g \left ( \frac i {a_N} \right )=\log \Big ( {\cal I}_1
(1/f_i)\Big ) =\log \left ( {\cal I}_1\left ( \frac{\beta}
{f_{min}} \left ( \frac {\beta f_{max}}{\alpha f_{min}}\right )
^{-i/a_N} \right ) \right ) }$.
\end{itemize}
Since for all $(H,\sigma)\in {\cal K}$, the function
$\displaystyle{x \mapsto L_{(H,\sigma)}(x)= \Big (A(H,\sigma)+x\cdot
B(H,\sigma)-g(x)\Big )^2}$ is an infinitely differentiable
function on $\R$, we know from the theory of Riemann sums that~:
\begin{multline*}
u_N(H,\sigma)=\frac 1 {a_N} \sum_{i=t^*_j-\eta a_N+1}^{t^*_j}
\left | A(H,\sigma) + \frac i {a_N} \cdot B(H,\sigma)- g\left ( \frac
i
{a_N} \right ) \right |^2  \\
\limiteN u(H,\sigma) =\int _{s_j^*-\eta }^{s^*_j}   \left (
A(H,\sigma) + x \cdot B(H,\sigma)- g\left ( x \right ) \right )^2\,
dx,
\end{multline*}
with $\displaystyle{s_j^*= \log \left ( \frac {\omega_j^*}{f_{min}}
\right ) \left ( \log \left ( \frac {\alpha f_{max}}{\beta
f_{min}} \right ) \right ) ^{-1} = \lim _{N \to \infty}\frac
{t^*_j}{a_N}} $. Moreover, the sequence $(u_N(H,\sigma))_N$
converges uniformly to $u(H,\sigma)$ because for $N$ large enough
\ban  \sup_{(H,\sigma)\in {\cal K}} |u_N(H,\sigma)-u(H,\sigma)|
\hspace{-5mm} &\leq \hspace{-5mm}&\left ( \frac 1 {a_N^2}+ \eta
\left |s_j^*-\frac {t_j^*}{a_N} \right | \right ) \cdot \hspace{-3mm}
\sup_{(H,\sigma)\in {\cal K}} \left \{ \sup _{0 \leq x \leq
(s_K^*+1)} \left |\frac {\partial
L_{(H,\sigma)}}{\partial x}(x) \right | \right \}\\
\hspace{-5mm}& \limiteN \hspace{-5mm}& 0, \ean since ${\cal K}$ is
a compact set of $[0,1]\times ]0,\infty[$ and thus
$\displaystyle{\sup_{(H,\sigma)\in {\cal K}} \left \{ \sup _{0
\leq x \leq (s_K^*+1)} \left |\frac {\partial
L_{(H,\sigma)}}{\partial x}(x) \right | \right \}<\infty}$. As a
consequence, from (\ref{Q21}) and since we assumed that
$(\widehat{H}_i^{(N)},\widehat{\sigma}_i^{(N)}) \in {\cal K}$ for
all $i=0,\cdots,K$, for some sufficiently small, fixed $\xi >0$ and
for all sufficiently large $N$,
\begin{eqnarray}
\label{Q22} Q_2 \geq a_N \left ( \int _{s_j^*-\eta }^{s^*_j} \left
( A(\widehat{H}_k^{(N)},\widehat{\sigma}_k^{(N)}) + x \cdot
B(\widehat{H}_k^{(N)},\widehat{\sigma}_k^{(N)})- g\left ( x \right
) \right )^2\, dx-\xi \right ).
\end{eqnarray}
But it is impossible that there exists $(a,b)\in \R^2$ such as
$g(x)=a+b\cdot x $ for all $x\in [s_j^*-\eta,s^*_j]$, {\em i.e.},
$\displaystyle{{\cal I}_1 \Big (c_1\cdot e^{c_2\cdot x} \Big )=e^a \cdot e^{b
\cdot x} }$ for all $x\in [s_j^*-\eta,s^*_j]$ with $\displaystyle{c_1=
\frac{\beta} {f_{min}}}$, $\displaystyle{c_2=\log \left ( \frac
{\alpha f_{min}}{\beta f_{max}}\right ) }$, which can also be
written as~:
\begin{eqnarray}
\label{impossible} {\cal I}_1 (x )= a_1 \cdot  x^{b_1} ~~\mbox{for
all}~~x \in [\alpha/\omega_j^*,\alpha/\omega_j^*+\eta'],
\end{eqnarray}
with $\eta'>0$ and $(a_1,b_1) \in \R^2$. Indeed, assume now
(\ref{impossible}) is true. But, for all $x \in
[\alpha/\omega_j^*,\alpha/\omega_j^*+\eta']$,
$$
{\cal I}_1(x )= 2 \left ( \sigma_{j-1}^{*2}\cdot x^{2H_{j-1}^*+1} \int
_{ \alpha}^{x\cdot \omega_j^*} \frac { |\widehat{\psi}(u)|^2}
{u^{2H_{j-1}^*+1}} \,  du +  \sigma_{j}^{*2}\cdot x^{2H_{j}^*+1} \int
_{x\cdot \omega_j^*}^{\beta}\frac { |\widehat{\psi}(u)|^2}
{u^{2H_{j}^*+1}} \, du \right ).
$$
Then $\displaystyle{\frac {\partial^n{\cal I}_1}{\partial
x^n}(\alpha/\omega_j^* )= a_1 \cdot \frac {\partial^n x^{b_1}}{\partial
x^n}(\alpha/\omega_j^* )}$ for $n=0,1$, which implies that
$b_1=(2H_j^*+1)$ and $a_1=2 \sigma_{j}^{*2} K_{H_j^*}(\psi)$
(here, we use the equality $\widehat{\psi}(\alpha)=0$). Thus, for
all $x \in [\alpha/\omega_j^*,\alpha/\omega_j^*+\eta']$,
\begin{eqnarray*}
\sigma_{j-1}^{*2}\cdot x^{2H_{j-1}^*+1} \int _{ \alpha}^{x\cdot \omega_j^*}
\frac { |\widehat{\psi}(u)|^2} {u^{2H_{j-1}^*+1}} \,  du
&=&\sigma_{j}^{*2}\cdot x^{2H_{j}^*+1} \int_{ \alpha} ^{x\cdot
\omega_j^*}\frac { |\widehat{\psi}(u)|^2} {u^{2H_{j}^*+1}} \, du,
\\
\Longrightarrow&& \hspace{-0.7cm}\int _{ \alpha/x}^{\omega_j^*}
|\widehat{\psi}(x\cdot y)|^2 \left ( \frac
{\sigma_{j-1}^{*2}}{y^{2H_{j-1}^*+1}}-\frac
{\sigma_{j}^{*2}}{y^{2H_{j}^*+1}} \right )\, dy=0,
\end{eqnarray*}
and hence $\displaystyle{ \left \{ \begin{array}{l}
\sigma_{j-1}^{*2}=\sigma_{j}^{*2} \\ H_{j-1}^*=H_{j}^*
\end{array} \right . }$. But this condition is impossible from
Assumption ($B_K$) and consequently there is no $(a,b) \in \R^2$
such as $g(x)=a+b\cdot x$ for all $x\in [s_j^*-\eta,s^*_j]$. \\
\\
The function $g$ belongs to the Hilbert space $\LL^2(
[s_j^*-\eta,s^*_j];dx)$. Since ${\cal L}=\{A+B\cdot x ,~x\in
[s_j^*-\eta,s^*_j],~(A,B)\in \R^2 \}$ is a closed linear subspace
of $\LL^2( [s_j^*-\eta,s^*_j];dx)$, there exits a distance between
$g$ and $\cal L$ in $\LL^2( [s_j^*-\eta,s^*_j];dx)$, {\em i.e.}
there exists $(\tilde {A}, \tilde {B})\in \R^2$ such as
$$
\int _{s_j^*-\eta }^{s^*_j}   \left ( \tilde {A} +\tilde { B} \cdot x-
g( x) \right )^2\, dx=\inf_{(A,B)\in \R^2}\int _{s_j^*-\eta
}^{s^*_j} \left ( A + B \cdot x- g( x) \right )^2\, dx = C >0,
$$
because $g \notin {\cal L}$. Then, by choosing $\xi$ such as
$0<\xi<C/2$, the inequality (\ref{Q22}) implies~:  \ba \label{Q2}
Q_2 \geq \frac C 2  \cdot a_N \ea for all sufficiently large $N$, with
$C$ a real positive number only depending on $\eta$, $ s_j^*$,
$H_{j-1}^*$, $H_j^*$, $\sigma_{j-1}^*$, $\sigma_j^*$ and $\psi$.
\item The previous evaluations of $Q_1$ and $Q_2$ provide an upper bound of
$Q_3$.We get \ba \label{Q3} \nonumber Q_3 & \leq & 2 \left ( Q_1
\right ) ^{1/2} \left ( \sum_{k=0}^{K+1}
\sum_{i=t_k+1}^{t_{k+1}-\tau_N }
( X_i\widehat{\lambda}_k - \log  {\cal I}_1 (1/f_i))^2 \right ) ^{1/2}\\
\nonumber & \leq & 2 \left ( Q_1\right )^{1/2} \times \left ( a_N \cdot
\sup _{f_{min}\leq f \leq f_{max}}\left \{2 \sup_{\lambda \in
{\cal K}}\{(\log f,1) \cdot \lambda)^2 +2\log^2 {\cal
I}_1(1/f) \right \} \right )^{1/2}, \\
 & =& O_P \left (\frac {a_N}{\sqrt{N \Delta _N}} \right
). \ea
\end{enumerate}
We deduce from (\ref{Q1}), (\ref{Q2})  and (\ref{Q3}) that
$Q_1=o(Q_2)$ and
$Q_3=o(Q_2)$, which implies \\
$\displaystyle{\Pr \left ( \min _{T \in V_{\varepsilon a_N}}
Q^{(N)}(T,\widehat{\Lambda}(T))\geq  \frac C 4 \cdot a_N \right )
\limiteN 1}$ and thus
$$\lim_{N \to \infty} \Pr   \left (\parallel \widehat {T}-T^* \parallel _\infty \geq
\varepsilon a_N \right )
=0~~\Longrightarrow~~\widehat{\omega}_i^{(N)} \limiteproba
\omega_i^* .$$ {\bf Second step :} For $j=1,\cdots,K$, we want to prove
that if $3/4 \leq p \leq 1$ and $0\leq q \leq 1$, for all
$\varepsilon>0$, there exists $0<C<\infty $ such as for
sufficiently large $N$, $\displaystyle{\Pr \left (a_N^{1-p} \left |
\widehat{\omega}_j^{(N)} -\omega_j^*\right |\geq C\right ) \leq
\varepsilon}$. \\  \\
{\em Mutatis mutandis}, we follow the same method as in the proof
of the convergence in probability. Now, let $0<p<1$,
$\displaystyle{0< \eta=\frac 1 2 \min \{{\underline \tau},\log
(\beta / \alpha) \}}$ and consider $\min _{T \in W^\eta _{Ca_N^p}}
Q^{(N)}(T,\widehat{\Lambda}(T))$ with
$$
W^\eta_{Ca_N^p}= \left \{ T \in {\cal A}_K ^{(N)},C a_N ^p \leq
\parallel T -t^*
\parallel_\infty  \leq \eta a_N \right \}.
$$
Then, as previously, for $T \in W^\eta_{Ca_N^p}$ and $N$ large
enough, it exists $j \in \{1,\cdots,K\}$ such as
\begin{eqnarray}\label{ti}
t_j +Ca_N^p\leq t_j^*<t_{j+1} -\tau_N
\end{eqnarray}
(the following proof is valid even if one considers the
alternative $t_j^* \leq t_j -Ca_N^p$). Then \ban
Q^{(N)}(T,\widehat{\Lambda}(T))& \geq &
\sum_{i=t^*_j+1}^{t_{j+1}-\tau_N } ( Y_i -\log  {\cal I}_1 (1/f_i)
)^2 +
( X_i\widehat{\lambda}_j -\log  {\cal I}_1 (1/f_i) )^2 + \\
&& \hspace{4cm}+ 2  ( Y_i- \log  {\cal I}_1 (1/f_i))(
X_i\widehat{\lambda}_j
- \log  {\cal I}_1 (1/f_i)) \\
& \geq  &  Q'_1+Q'_2+Q'_3. \ean 1. First, we have again,
\begin{eqnarray} \label{Q'1}
Q_1'= O_P \left (\frac {a_N}{N \Delta _N} \right ).
\end{eqnarray}
2. Secondly, $\displaystyle{ Q'_2 =  \sum_{i=t_j^*+1}^{t_{j+1}
-\tau_N} ( X_i\widehat{\lambda}_j - \log {\cal I}_1(1/f_i) )^2}$.
But we know $\log {\cal I}_1(1/f_i)= X_i\lambda_j^*$ for $i \in
\{t_j^*+1,\cdots,t_{j+1}-\tau_N\}$. Moreover, for $a_i=1/f_i$, $i\in
\{t_j+1,\cdots,t_j^* \}$ and $N$ large enough, $a_i \simeq \alpha
/\omega_i^* $, and
$$
{\cal I}_1 (a_i) = {\cal I}_1 \left (\frac \alpha {\omega_j^* }
\right )+ \left (a_i-\frac \alpha {\omega_j^* } \right ) {\cal
I}'_1 \left (\frac \alpha {\omega_j^* } \right )+ O\left (a_i-\frac
\alpha {\omega_j^* } \right )^2.
$$
But $\displaystyle{ {\cal I}_1(a_i)=2 \left (
\sigma_{j-1}^{*2}a_i^{2H_{j-1}^*+1} \int _{ \alpha}^{a_i\omega_j^*}
\frac { |\widehat{\psi}(u)|^2} {u^{2H_{j-1}^*+1}}\, du +
\sigma_{j}^ {*2}a_i^{2H_{j}^*+1} \int _{a_i
\omega_j^*}^{\beta}\frac { |\widehat{\psi}(u)|^2}
{u^{2H_{j}^*+1}}\, du \right )}$ and \\
$\displaystyle{   {\cal I}'_1 \left (\frac \alpha {\omega_j^* }
\right )=2\sigma_{j}^ {*2}K_{H_j^*}(\psi) (2H_j^*+1) \left ( \frac
\alpha {\omega_j^* } \right ) ^{2H_j^*} }$; thus for $i\in
\{t_j+1,\cdots,t_j^* \}$, \ba \label{I} \log {\cal I}_1 (1/f_i) =
X_i\lambda_j^* +  \left [(2H_j^*+1)\frac {f_{min}}{\beta} \log
\left ( \frac {\beta f_{max}}{\alpha f_{min}} \right )\right ] \cdot
\left ( \frac { t_j^*-i}{a_N} \right ) + O\left (\frac {
t_j^*-i}{a_N} \right ) ^2.
\end{eqnarray}
Then, with $\widehat{\lambda}_j=(\widehat{a}_j,\widehat{b}_j)'$,
one gets for $i  \in \{t_j^*+1,\cdots,t_{j+1}-\tau_N\}$,
\begin{eqnarray}
\label{aj} \left (X_i\widehat{\lambda}_j - \log {\cal I}_1(1/f_i)
\right) =(\log f_i-\overline{log f})(\widehat{a}_j-a_j^*)+
\overline{Z},
\end{eqnarray}
$\overline{XXX}$ indicates the empirical mean of $XXX$ between
$t_j+1$ and $t_{j+1}-\tau_N$. Thus,
\begin{eqnarray}
\label{Q'21} Q'_2&\geq & \sum_{i=t_j^*+1}^{t_{j+1} -\tau_N} \left
(\left ( \log f_i -\overline{\log f}  \right )(\widehat{a}_j -
a^*_j)+\frac 1 {\sqrt{N \Delta_N}} \overline{Z} \right ) ^2.
\end{eqnarray}
We also have~:
\begin{eqnarray*}
\widehat{a}_j&= & \frac { \displaystyle{\sum _{i=t_j+1} ^
{t_{j+1}-\tau_N} \left ( \log f_i -\overline{\log f}  \right )
\left ( Y_i - \overline{Y} \right ) }} {  \displaystyle{ \sum
_{i=t_j+1}^{t_{j+1}-\tau_N} \left ( \log f_i -\overline{\log f}
  \right )^2}  } \\
&= & \frac { \displaystyle{\sum _{i=t_j+1} ^ {t_{j+1}-\tau_N}
\left ( \log f_i -\overline{\log f}  \right ) \left ( \log {\cal
I}_1(1/f_i)+ \frac 1 {\sqrt{N \Delta_N}}Z_i^{(N)} - \overline{ \log
{\cal I}_1}- \frac 1 {\sqrt{N \Delta_N}} \overline{Z} \right )  }}
{  \displaystyle{ \sum _{i=t_j+1}^{t_{j+1}-\tau_N} \left ( \log
f_i -\overline{\log f}
  \right )^2}  },
\end{eqnarray*}
and thus,
\begin{multline}
\label{borneaj} \widehat{a}_j - a^*_j= \frac { \displaystyle{\sum
_{i=t_j+1} ^ {t^*_j} \left ( \log f_i -\overline{\log f}  \right )
\left ( \log {\cal I}_1(1/f_i) -  X_i'\lambda_j^* \right ) }} {
\displaystyle{ \sum _{i=t_j+1}^{t_{j+1}-\tau_N} \left ( \log f_i
-\overline{\log f} \right )^2}  }\\
+ \frac 1 {\sqrt{N \Delta_N}} \frac { \displaystyle{\sum
_{i=t_j+1} ^ {t_{j+1}-\tau_N} \left ( \log f_i -\overline{\log f}
\right ) \left ( Z_i^{(N)} - \overline{Z} \right )  }} {
\displaystyle{ \sum _{i=t_j+1}^{t_{j+1}-\tau_N} \left ( \log f_i
-\overline{\log f}  \right )^2}  }.
\end{multline}
From the definition of $(\log f_i)$,
\begin{eqnarray}
\label{logf} \sum _{i=t_j+1}^{t_{j+1}-\tau_N} \left ( \log f_i
-\overline{\log f} \right )^2 \simeq \left [ \frac 1 {12} \log
\left ( \frac {\beta f_{max}} {\alpha f_{min}} \right )\right
](t_{j+1}-\tau_N-t_j)=O(a_N).
\end{eqnarray}
Expansions (\ref{logf}) and (\ref{I}) imply there exist two
constants $C_1>0$ and $C_2 >0$ such as for $N$ large enough~:
\begin{eqnarray}
\nonumber C_1 \left ( \frac {t^*_j-t_j}{a_N} \right )^2 \leq \left
| \frac { \displaystyle{\sum _{i=t_j+1} ^ {t^*_j} \left ( \log f_i
-\overline{\log f}  \right ) \left ( \log {\cal I}_1(1/f_i) -
X_i'\lambda_j^* \right ) }} {  \displaystyle{ \sum
_{i=t_j+1}^{t_{j+1}-\tau_N} \left ( \log f_i -\overline{\log f}
\right )^2}  } \right |  \leq C_2 \left ( \frac {t^*_j-t_j}{a_N}
\right )^2.
\end{eqnarray}
Moreover \ban \frac 1 {\sqrt{N \Delta_N}} \frac {
\displaystyle{\sum _{i=t_j+1} ^ {t_{j+1}-\tau_N} \left ( \log f_i
-\overline{\log f} \right ) \left ( Z_i^{(N)} - \overline{Z} \right
)  }} {  \displaystyle{ \sum _{i=t_j+1}^{t_{j+1}-\tau_N} \left (
\log f_i -\overline{\log f} \right )^2}  }=O_P \left ( \frac 1
{\sqrt{N \Delta_N}} \right ). \ean
%\widehat{a}_j=a_j^*+3 (2H_j^*+1)\frac {C^2 }{\eta '}\frac {f_{min}}{\alpha}\log   \left ( \frac {\alpha f_{max}}{\beta f_{min}} \right ) \left ( \frac {(N\Delta_N)^{2p}}{(N \Delta_N)^2} \right ) \left ( 1 + o(1) \right )
Thus, we deduce from (\ref{borneaj}) that~:
$$
C_1  \left ( \frac {t^*_j-t_j}{a_N} \right )^2 + O_P \left ( \frac
1 {\sqrt{N \Delta_N}} \right ) \leq \left | \widehat{a}_j-a_j^*
\right |.
$$
As a consequence, for $(p,q)$ such as $4q(1-p)\leq 1$ (for
instance, $p=3/4$ and $q=1$), then $\displaystyle{\left ( \frac
{t^*_j-t_j}{a_N} \right )^2\cdot \sqrt{N \Delta_N} \geq C^2}$, and
thus for all $\varepsilon>0$, for $N$ sufficiently large, we can
chose $C>0$ such as~:
\begin{eqnarray}
\label{encaH} \Pr  \left ( \frac {C_1^2}2 \left ( \log f_i
-\overline{\log f}  \right )^2 \left ( \frac {t^*_j-t_j}{a_N}
\right )^4 \leq \left ( \left ( \log f_i -\overline{\log f}  \right
)(\widehat{a}_j - a^*_j) \right ) ^2 \right ) \geq 1-\varepsilon.
\end{eqnarray}
Now, from (\ref{Q'21}), (\ref{encaH}) and with
$\displaystyle{\Pr(t_{j+1}-\tau_N- t_j^*\geq \frac \eta 2 \, a_N)
\limiteN 1}$, for $(p,q) \in [3/4,1]\times [0,1]$, for all
$\varepsilon>0$, for $N$ sufficiently large, we can also chose
$C>0$ such as~:
\begin{eqnarray}
\nonumber \Pr  \left ( \frac {C_1^2}4 \left ( \frac
{t^*_j-t_j}{a_N} \right )^4  \cdot \sum _{i=t_j^*+1}^{t_{j+1}-\tau_N}
\left ( \log f_i -\overline{\log f}
\right )^2  \leq Q'_2 \right ) &\geq& 1-\varepsilon, \\
\label{Q'2}\Longrightarrow ~~~ \Pr  \left ( C^4 \cdot C_2 \cdot a_N^{4p-3}
\leq Q'_2 \right ) &\geq& 1-\varepsilon,
\end{eqnarray}
with $C_2>0$ a real number not depending on $C$, $N$ and
$\varepsilon$. \\ \\
%because $\displaystyle{t_{j+1} -\tau_N -t_j > \frac \eta 2 a_N}$ and $|t^*_j-t_j| \geq C a_N^p$.
3. Finally, from the classical bound of $Q'_3$, we obtain,
\begin{eqnarray*}
Q'_3 \leq  2 \cdot  \left(  Q'_2  \right ) ^{1/2} \cdot \left(  Q'_1
\right ) ^{1/2}.
\end{eqnarray*}
But, following a similar method as previously, from (\ref{encaH}
one can find a upper-bound for $Q_2'$, {\em i.e.} for $(p,q) \in
[3/4,1]\times [0,1]$, for all $\varepsilon>0$, for $N$ sufficiently
large, we can also chose $C>0$ such as~:
\begin{eqnarray*}
\Pr  \left ( Q'_2 \leq C^4 \cdot C_3 \cdot a_N^{4p-3}   \right ) &\geq&
1-\varepsilon,
\end{eqnarray*}
with $C_3>0$ a real number not depending on $C$, $N$ and
$\varepsilon$. Thus, for $(p,q) \in [3/4,1]\times [0,1]$, for all
$\varepsilon>0$, we can also chose $C>0$ such as~:
\begin{eqnarray}
\label{Q'3} \Pr  \left ( Q'_3 \leq C^2 \cdot C_4 \cdot
\frac{a_N^{2p-2}}{\sqrt{N \Delta_N}} \right ) &\geq&
1-\varepsilon,
\end{eqnarray}
with $C_4>0$ a real number not depending on $C$ and $N$.\\ \\
Now, from (\ref{Q'1}), (\ref{Q'2}) and (\ref{Q'3}), one deduces
that for $(p,q) \in [3/4,1]\times [0,1]$, for all $\varepsilon>0$, for
$N$ sufficiently large, we can chose $C>0$ sufficiently large such
as~:
$$
\Pr \left (\min _{T \in W^\eta _{Ca_N^p}}
Q^{(N)}(T,\widehat{\Lambda}(T)) \geq C^4 \cdot \frac {C_2} 2 \cdot
a_N^{4p-3} \right ) \geq 1- \varepsilon.
$$
and thus like $\displaystyle{Q^{(N)}_* = O_P \left ( \frac {a_N}{N
\Delta_N}\right )}$ from (\ref{Q*}),
$$
\Pr \left (\min _{T \in W^\eta _{Ca_N^p}}
Q^{(N)}(T,\widehat{\Lambda}(T)) \leq  Q_*^{(N)} \right ) \leq
\varepsilon,
$$
that leads to $\displaystyle{\Pr \left (a_N^{1-p} \left |
\widehat{\omega}_j^{(N)} -\omega_j^*\right |\geq C\right ) \leq
\varepsilon}$ for sufficiently large $C$ and $N$.\end{dem}
\begin{dem}[{\bf Proposition \ref{HetK}}]
From Proposition \ref{omegas}, we deduce that $\forall j=0,\cdots,K$,
$$
\Pr \left (  [\tilde{U}_j^{(N)},\tilde{V}_j^{(N)} ] \subset [t_j^*,
t_{j+1}^* -\tau_N] \right ) \limiteN  1.
$$
Denote $A_j^{(N)}$ the event $
[\tilde{U}_j^{(N)},\tilde{V}_j^{(N)} ] \subset [t_j^*, t_{j+1}^*
-\tau_N]$. Then, $\forall j=0,\cdots,K$ and $\forall (x,y) \in \R^2$,
\begin{multline*}
\Pr \left ( \sqrt{N \Delta_N} \left
(\tilde{\lambda}_j^{(N)}-\lambda_j^*\right ) \in ]-\infty,x]\times
]-\infty,y] \right )
\\
=\Pr \left (A_j^{(N)} \right )\times \Pr \left (\sqrt{N \Delta_N} \left
(\tilde{\lambda}_j^{(N)}-\lambda_j^*\right ) \in ]-\infty,x]\times
]-\infty,y]~ |~A_j^{(N)} \right )+
\\
+ \Pr \left ( \overline{A_j^{(N)}}\right )\times \Pr \left ( \sqrt{N
\Delta_N} \left (\tilde{\lambda}_j^{(N)}-\lambda_j^*\right ) \in
]-\infty,x]\times ]-\infty,y]~ |~\overline{A_j^{(N)}}\right ).
\end{multline*}
Now, since $\Pr \left ( \sqrt{N \Delta_N} \left
(\tilde{\lambda}_j^{(N)}-\lambda_j^*\right ) \in ]-\infty,x]\times
]-\infty,y]~ |~\overline{A_j^{(N)}}\right ) \leq 1$ and $\Pr \left (
\overline{A_j^{(N)}}\right )=1 - \Pr\left (A_j^{(N)} \right )$, we
obtain~:
\begin{multline}
\Pr\left (A_j^{(N)} \right ) \cdot \Pr \left ( \sqrt{N \Delta_N} \left
(\tilde{\lambda}_j^{(N)}-\lambda_j^*\right ) \in ]-\infty,x]\times
]-\infty,y] ~ |~ A_j^{(N)}  \right )
\\ \hspace{-2.5 cm}
\leq \Pr \left ( \sqrt{N \Delta_N} \left
(\tilde{\lambda}_j^{(N)}-\lambda_j^*\right )
\in ]-\infty,x]\times ]-\infty,y]   \right ) \\
\label{probcond}  \leq   \Pr \left ( \sqrt{N \Delta_N} \left
(\tilde{\lambda}_j^{(N)}-\lambda_j^*\right ) \in ]-\infty,x]\times
]-\infty,y] ~ |~ A_j^{(N)}  \right ) + 1-\Pr\left (A_j^{(N)} \right
).
\end{multline}
Since $~\widehat{\omega}_j^{(N)}\limiteproba \omega_j^*\;~$ and
$~\;\widehat{\omega}_{j+1}^{(N)}\limiteproba \omega_{j+1}^*,\;$
therefore
$~\;(\widehat{\omega}_j^{(N)},\widehat{\omega}_{j+1}^{(N)})
\limiteproba (\omega_j^*,\omega_{j+1}^*),\;$ we have\\ $~(f_k)_{k \in
\{\tilde{U}_j^{(N)},\cdots,\tilde{V}_j^{(N)} \}} \limiteproba
(g_j^*(k))_{1\leq k \leq m}$ and $~\tilde{X}_j ^{(N)}\limiteproba
X^{*}_j.\;$ Thus, from Proposition \ref{propTLC} and central limit
theorem (\ref{TLCdis}), for all $(x_k)_{1\leq k \leq m} \in
\R^m$, we get
\begin{multline*}
\Pr \left (\sqrt {N\Delta_N} \left ( \tilde{Y}_j^{(N)}-
\tilde{X}^{(N)}_j\lambda_j^* \right )\in \prod _{k=1}^m
]-\infty,x_k]~~ |~ A_j^{(N)} \right ) - \Pr \left (\tilde{Z}_j \in
\prod _{k=1}^m ]-\infty,x_k]~~ |~ A_j^{(N)} \right )\limiteN 0,
\end{multline*}
with $\tilde{Z}_j \egaleloi {\cal N}_m(0,\Sigma_j^*)$ and
$\displaystyle{\Sigma_j^*=\left ( \cov \Big( Z \Big ( \frac 1
{g_j^*(k)} \Big ), Z \Big ( \frac 1 {g_j^*(l)} \Big ) \Big )
\right )_{1\leq k,l \leq m}}$ (it explains the expression
(\ref{Sigmaj}) of $\Sigma_j^*$). From the equality
$\displaystyle{\tilde{\lambda}_j^{(N)}=\Big ( (\tilde{X}_j^{(N)})'
\tilde{X}_j^{(N)}\Big )^{-1}(\tilde{X}_j^{(N)})'
\tilde{Y}_j^{(N)}}$, we deduce that for all $(x,y) \in \R^2$, with
$\tilde{\xi}_j \egaleloi {\cal N}_2 (0,\Gamma_1^{\lambda_j^*})$
and $\displaystyle{\Gamma_1^{\lambda_j^*}=\left (  X^{*'}_j
X^{*}_j \right )^{-1}X^{*}_j\Sigma_j^* X^{*'}_j\left (  X^{*'}_j
X^{*}_j \right )^{-1}}$,
\begin{eqnarray}
\label{TLCxi} \Pr \left ( \sqrt{N \Delta_N} \left
(\tilde{\lambda}_j^{(N)}-\lambda_j^*\right ) \in ]-\infty,x]\times
]-\infty,y]  ~ |~ A_j^{(N)}  \right ) -\Pr \left ( \tilde{\xi}_j
\in ]-\infty,x]\times ]-\infty,y]~ |~ A_j^{(N)}\right )\limiteN 0.
\end{eqnarray}
We also have~:
\begin{eqnarray}
\nonumber \Pr \left ( \tilde{\xi}_j \in ]-\infty,x]\times
]-\infty,y]\right )+\Pr \left (A_j^{(N)}\right )-1  \leq
&&  \\
\label{xi}&& \hspace{-5cm}\leq \Pr \left ( \tilde{\xi}_j \in
]-\infty,x]\times ]-\infty,y]~ |~ A_j^{(N)}\right ) \leq \frac {\Pr
\left ( \tilde{\xi}_j \in ]-\infty,x]\times ]-\infty,y]\right )}{\Pr
\left (A_j^{(N)}\right )}.
\end{eqnarray}
Now, as $\displaystyle{\Pr \left (A_j^{(N)}\right ) \limiteN 1}$,
from (\ref{probcond}), (\ref{TLCxi}) and (\ref{xi}), we deduce
that for all $(x,y)\in \R^2$~:
$$
\Pr \left ( \sqrt{N \Delta_N} \left
(\tilde{\lambda}_j^{(N)}-\lambda_j^*\right ) \in ]-\infty,x]\times
]-\infty,y]   \right ) \limiteN \Pr \left ( \tilde{\xi}_j \in
]-\infty,x]\times ]-\infty,y]\right ),
$$
that achieves the proof.
\end{dem}
\begin{dem}[{\bf Proposition \ref{HetK2}}]
First, from the expression of each $s_{kl}$ given in
(\ref{Sigmaj}) and with ${\cal M}_m(\R)$ the set of real
$m$-by-$m$ matrix, the function $\Sigma:(H,u,v)\mapsto
\Sigma(H,u,v)\in {\cal M}_m(\R)$ is a continuous (and therefore
measurable) function of $(H,u,v)$ for $H$ in a compact set
included in $]0,1[$ and $(u,v)\in ]f_{min},f_{max}[^2$. For all
$j=0,\cdots,K$, we have~:
\begin{enumerate}
\item from Assumptions ($B_K)$ and ($C$),
$(\tilde{H_j}^{(N)},\tilde{\sigma}_j^{(N)} )\in {\cal K}$ and
$(\widehat{\omega}_j^{(N)},\widehat{\omega}_{j+1}^{(N)}) \in
]f_{min},f_{max}[^2$;
\item from (\ref{conv_omega}) and (\ref{conv_lambda}),
$\tilde{H_j}^{(N)} \limiteproba H_j^*$,
$~\widehat{\omega}_j^{(N)}\limiteproba \omega_j^*$,
$~\widehat{\omega}_{j+1}^{(N)}\limiteproba \omega_{j+1}^*$ and
therefore
$$
(\tilde{H_j}^{(N)},\widehat{\omega}_j^{(N)},
\widehat{\omega}_{j+1}^{(N)}) \limiteproba
(H_j^*,\omega_j^*,\omega_{j+1}^*).
$$
\end{enumerate}
As a consequence,
$\widehat{\Sigma}_j^{(N)}=\Sigma(\tilde{H_j}^{(N)},
\widehat{\omega}_j^{(N)},\widehat{\omega}_{j+1}^{(N)}) \limiteproba
\Sigma_j^*$, for all $j=0,\cdots,K$, and since $\Sigma(H,u,v)$ is an
invertible covariance matrix for all $(H,u,v)\in ]0,1[ \times
]f_{min},f_{max}[^2$,
\begin{eqnarray}\label{fgls1}
\Big(\widehat{\Sigma}_j^{(N)}\Big ) ^{-1} \limiteproba \Big (
\Sigma_j^*\Big )^{-1},~~\mbox{for all}~~j=0,\cdots,K.
\end{eqnarray}
Secondly, denote $\displaystyle{\left \{  \begin{array}{l}
\tilde{M}^{(N)}_j=\Big ( (\tilde{X}^{(N)}_j)' \Big ( \Sigma_j^*\Big
) ^{-1} \tilde{X}^{(N)}_j \Big )^{-1}( \tilde{X}^{(N)}_j)'  \Big (
\Sigma_j^*\Big ) ^{-1}\\
\widehat{M}_j^{(N)}=\Big ( (\tilde{X}^{(N)}_j)' \Big (
\widehat{\Sigma}_j^{(N)}\Big ) ^{-1} \tilde{X}^{(N)}_j  \Big
)^{-1}( \tilde{X}^{(N)}_j)'  \Big ( \widehat{\Sigma}_j^{(N)}\Big )
^{-1} \end{array} \right . }$. \\
The $2$-by-$m$ matrix $\tilde{M}^{(N)}_j$ verifies~:
$$
\tilde{\lambda}^{(N)}_j=\tilde{M}^{(N)}_j\tilde{Y}^{(N)}_j=
\lambda_j^*+\frac 1 {\sqrt{N \Delta_N}}\tilde{M}^{(N)}_j
\tilde{Z}^{(N)}_j
$$
with $\tilde{Z}^{(N)}_j=\displaystyle{\Big ( Z^{(N)} (1/f_i)\Big
)_{i \in \{\tilde{U}_j^{(N)},\cdots,\tilde{V}_j^{(N)} \}}}$ and
$\tilde{Z}^{(N)}_j \limiteloi \tilde{Z}_j=\Big (Z(1/g_j^*(k)) \Big
)_{1\leq k \leq m}$ from the central limit theorem (\ref{TLCdis}).
In the same way,
$$
\underline{\lambda}^{(N)}_j=\widehat{M}^{(N)}_j\tilde{Y}^{(N)}_j=
\lambda_j^*+\frac 1 {\sqrt{N \Delta_N}}\widehat{M}^{(N)}_j
\tilde{Z}^{(N)}_j.
$$
From (\ref{fgls1}), we obtain $\widehat{M}^{(N)}_j
-\tilde{M}^{(N)}_j\limiteproba 0$, and thus,
$$
{\sqrt{N \Delta_N}} \left ( \underline{\lambda}^{(N)}_j-
\lambda_j^* \right)  - \tilde{M}^{(N)}_j \tilde{Z}^{(N)}_j
\limiteproba 0,
$$
with $\tilde{M}^{(N)}_j \tilde{Z}^{(N)}_j \limiteloi {\cal
N}_2(0,\Gamma_2^{\lambda_j^*})$ (the same covariance matrix as
that obtained with a generalized least squares estimation), and
this implies Proposition \ref{HetK2}.
\end{dem}
\begin{dem}[{\bf Proposition \ref{test}}]
For each $j=0,\cdots,K$, one first show that
\begin{eqnarray}\label{cochran}
N\Delta_N \cdot  \parallel  \tilde{Y}_j^{(N)} -\tilde{X}_j
^{(N)}\underline{\lambda}_j^{(N)}
\parallel ^2_{\widehat{\Sigma}_j^{(N)}} \limiteloi \chi ^2 (m-2).
\end{eqnarray}
Indeed, $\displaystyle{\parallel \tilde{Y}_j^{(N)} -\tilde{X}_j
^{(N)}\underline{\lambda}_j^{(N)}
\parallel ^2_{\widehat{\Sigma}_j^{(N)}}=\parallel
\widehat{P}_{j\bot}^{(N)} \tilde{Y}_j^{(N)}  \parallel
^2_{\widehat{\Sigma}_j^{(N)}}=\frac 1 {N\Delta_N} \parallel
\widehat{P}_{j\bot}^{(N)} \tilde{Z}_j^{(N)}  \parallel
^2_{\widehat{\Sigma}_j^{(N)}}}$ where
$\widehat{P}_{j\bot}^{(N)}=I_m -\tilde{X}^{(N)}_j
\widehat{M}_{j}^{(N)} $ is the matrix of the orthogonal projector
in $\R^m$ on the orthogonal of $V_j$, where
$V_j=\{\tilde{X}^{(N)}_j \lambda,~\lambda \in \R^2\}$ is the
$2$-dimensional subspace of $\R^m$ generated by
$\tilde{X}^{(N)}_j$ (here the notion of orthogonality is based on
the inner product $<u,v>_{\widehat{\Sigma}_j^{(N)}}=u'\cdot \Big (
\widehat{\Sigma}_j^{(N)}\Big ) ^{-1}\cdot v$ for $u,v\in \R^m$). From
the previous proofs, we know~:
\begin{itemize}
\item $\widehat{\Sigma}_j^{(N)} \limiteproba
\Sigma_j^*$, $~\tilde{X}_j ^{(N)}\limiteproba X^{*}_j$ and
therefore $\displaystyle{\widehat{P}_{j\bot}^{(N)} \limiteproba
P^*_{j\bot}}$ where \\
$\displaystyle{P^*_{j\bot}=\left (I_m-X_j^{*}\left ( X_j^{*'}
(\Sigma_j^*)^{-1}X_j^{*} \right) ^{-1}X_j^{*'}(\Sigma_j^*)^{-1}
\right)}$ is the matrix of an orthogonal projector on a
$(m-2)$-dimensional subspace of $\R^m$;
\item $<u,v>_{\widehat{\Sigma}_j^{(N)}} \limiteproba
<u,v>_{\Sigma_j^*}$ for $u,v\in \R^m$;
\item $\tilde{Z}_j^{(N)} \limiteproba \tilde{Z}_j $ with
$\tilde{Z}_j \egaleloi {\cal N}_m(0,\Sigma_j^*)$.
\end{itemize}
Consequently, $\parallel \widehat{P}_{j\bot}^{(N)}
\tilde{Z}_j^{(N)}  \parallel ^2_{\widehat{\Sigma}_j^{(N)}}
\limiteloi \parallel P_{j\bot}^* \tilde{Z}_j
\parallel ^2_{\Sigma_j^*}$. From Cochran's Theorem, we know
$\parallel P_{j\bot}^* \tilde{Z}_j \parallel ^2_{\Sigma_j^*}
\egaleloi \chi^2(m-2)$ and therefore (\ref{cochran}) is proved. \\ \\
Moreover, with the notations of Proposition \ref{propTLC}, if
$\log f \geq \log f' + \log \beta /\alpha$ then $\cov
(Z(1/f),Z(1/f'))=0$. But for all $(i,j) \in \{0,\cdots,K \}^2$, $i\neq
j$, $\forall k \in \{\tilde{U}_i^{(N)},\cdots,\tilde{V}_i^{(N)} \}$
and  $\forall k' \in \{\tilde{U}_j^{(N)},\cdots,\tilde{V}_j^{(N)} \}$,
$|\log f_k - \log f_{k'}|\geq  \log \beta /\alpha$. Thus, we
deduce that the different $\underline{\lambda}_j^{(N)}$ are
asymptotically Gaussian and independent. It provides the end of
the proof of the Proposition \ref{test}.
\end{dem}
~\\
{\bf Acknowledgments.} The authors are extremely grateful to the
anonymous referee for a very careful reading and many relevant
suggestions and corrections that strongly improve the content and
the form of the paper. From the important work that he has
generously made, it should be fair to consider him as a co-author
of the article.
%%%%%%%%%%%%%%%%%%%%%%%%%%%%%%%%%%%%%%%%%%%%%%%%%%%%%%%%%%%%%%%%


\begin{thebibliography}{99}

\bibitem{Abry} Abry, P., Flandrin, P., Taqqu, M.S. and Veitch, D. (2002).
Self-similarity and long-range dependence through the wavelet
lens, in  {\it Long-range Dependence: Theory and Applications}, P.
Doukhan, G. Oppenheim and M.S. Taqqu editors, Birkh{\"a}user.

\bibitem{Amemiya} Amemiya,T. (1985). {\it Advanced Econometrics.}
Cambridge : Harvard University Press.

\bibitem{AyLV:99} Ayache, A. and L\'evy~V\'ehel, J. (1999):  Generalized multifractional Brownian
motion: definition and preliminary results, {\it  in} M.~Dekking,
J.~L\'evy~V\'ehel, E.~Lutton \& C.~Tricot, eds, "Fractal theory an
applications in engineering", Springer Verlag.

\bibitem{Bai} Bai J. (1998). Least squares estimation of a shift in linear processes.
{\it J. of Time Series Anal.} 5, p. 453-472.

\bibitem{Bai-Per} Bai J. and Perron P. (1998). Estimating and testing linear models
with multiple structural changes. {\it Econometrica} 66, p. 47-78.

\bibitem{Bardb} Bardet J.M. (2000). Testing for the presence of self-similarity of
Gaussian time series having stationary increments. {\it J. of Time
Series Anal.} 21, p. 497-516.

\bibitem{Barda} Bardet J.M. (2002). Statistical study of the wavelet analysis of
fractional Brownian motion.  {\it IEEE Trans. Inform. Theory.} 48,
p. 991-999.

\bibitem{BaBe} Bardet J.M. and Bertrand, P. (2003). Definition,
properties and wavelet analysis of multiscale fractional Brownian
motion. Preprint LSP, Toulouse III.

\bibitem{Blms} Bardet J.M., Lang G., Moulines E. and Soulier P. (2000).
Wavelet estimator of long-range dependent processes. {\it
Statistical Inference for Stochastic Processes} 3, p. 85-99.

\bibitem{Benassi98} Benassi, A., Cohen, S. and Istas, J. (1998).
Identifying the multifractional function of a Gaussian process.
{\it Statistics and probability letters} 39, p. 337-345.

\bibitem{BDeguy:1999} Benassi, A. and  Deguy, S. (1999). Multi-scale fractional
Brownian motion : definition and identification. Preprint LAIC.

\bibitem{BJR:Rev.Math.Ibero}  Benassi, A., Jaffard, S. and Roux, D. (1997). Elliptic Gaussian random processes.
{\it Rev. Mathem\`atica Iberoamericana} 13 (1), p. 19-90.

%% \bibitem{BP:scie} Bertrand, P. (2000). A local method for estimating change points : the hat-function. {\it Statistics} 34, p. 215-235.

\bibitem{B2DMV:IEEE} Bertrand, P., Bardet, J.M, Dabonneville, M., Mouzat, A.
and Vaslin, P. (2001). Automatic determination of the different
control mechanisms in upright position by a wavelet method, {\it
IEEE Engineering in Medicine and Biology Society} p. 25-28.

\bibitem{Bil} Billingsley, P. (1968). {\it Convergence of Probability Measures}.
New-York. Wiley.

\bibitem{Ch:2000}  Cheridito, P. (2003). Arbitrage in Fractional Brownian
Motion Models. {\it Finance and Stochastics} 7, p. 533-553.

\bibitem{Cohen:2000} Cohen, S. (2000). Champs localement auto-similaires.
In {\it Fractals et loi d'\'echelles},  P.Abry, P. Goncalv\`es and  J.
L\'evy-Vehel editors, Herm\`es, Paris

\bibitem{Collins:DeLuca:93} Collins J.J. and de Luca C.J. (1993).
Open-loop and closed-loop control of posture : A random walk
analysis of center-of-pressure trajectories. {\it  Experimental
Brain Research} 9, p. 308-318.

\bibitem{dahl89} Dahlhaus R. (1989). Efficient parameter estimation
for self-similar processes. {\em Ann. Statist.} 17, p. 1749-1766

\bibitem{flan2} Flandrin, P. (1992). Wavelet analysis and synthesis of
fractional Brownian motion. {\em IEEE Trans. on Inform. Theory}
38, p. 910-917.

\bibitem{fox-taq1} Fox, R. and Taqqu, M.S. (1986). Large-sample properties
of parameter estimates for strongly dependent Gaussian time
series. {\em Ann. Statist.} 14, p. 517-532.

\bibitem{gira90} Giraitis L. and  Surgailis D. (1990). A central limit
theorem for quadratic forms in strongly dependent linear variables
and its applications to the asymptotic normality of Whittle
estimate. {\em Prob. Th. and Rel. Fields.} 86, p. 87-104.

%% \bibitem{Hunt :51} Hunt, G.A. (1951). Random Fourier transforms. {\it Trans. of the American Society of Civil Engineers} 116, p. 770-808.

\bibitem{IW} Ikeda, N. and Watanabe, S. (1989). {\it Stochastic Differential Equations and Diffusion
Processes}, North Holland Publishing Co.

\bibitem{is-lan} Istas, J. and Lang, G. (1997). Quadratic variations
and estimation of the local H\"older index of a Gaussian process.
{\em Ann. Inst. Poincar\'e} 33, 407-436.

\bibitem{JS} Jacod, J. and Shiryaev, A.N. (1987). {\it Limit theorems for stochastic
processes}, Springer Verlag.

\bibitem{Lav} Lavielle, M. (1999). Detection of multiple changes in a
sequence of dependent variables. {\it Stoc. Proc. Appl.} 83, p.
79-102.

\bibitem{Lav-Mou} Lavielle, M. and Moulines, E. (2000). Least-squares
estimation of an unknown number of shifts in a time series. {\it
J. of Time Series Anal.} 21, p. 33-59.

\bibitem{MvN:68} Mandelbrot, B. and  Van Ness J. (1968). Fractional
Brownian motion, fractional noises and applications. {\it SIAM
review} 10, p.422-437.

\bibitem{PLV:96} Peltier, R. and L\'evy~V\'ehel, J. (1995). Multifractional Brownian motion : definition and preliminary
results, Technical Report 2645, INRIA, Le Chesnay, France.

\bibitem{rogers:97} Rogers, L.C.G. (1997). Arbitrage with fractional
Brownian motion. {\it Mathematical Finance} 7, p. 95-105.

\bibitem{TaSa:1994} Samorodnitsky, G. and Taqqu M.S. (1994). {\it Stable
non-Gaussian Random Processes}, Chapman and Hall.

\bibitem{vv} Van der Vaart, A. (1998). {\it Asymptotic statistics}, Cambridge Series in
Statistical and Probabilistic Mathematics, Cambridge.
\end{thebibliography}
\end{document}